\documentclass[11pt, DIV10,a4paper]{article}

\usepackage{natbib}
\newcommand {\ctn}{\citet} 
\newcommand {\ctp}{\citep}       
\usepackage{graphicx,subfigure,amsmath,latexsym,amssymb}
\usepackage{float,epsfig,multirow,rotating,times}
\usepackage{upgreek,wrapfig}
\usepackage{comment}
\usepackage[colorlinks,citecolor=blue,urlcolor=blue,filecolor=blue,backref=page]{hyperref}
\usepackage{enumitem}
\usepackage{physics}

\newcommand{\mylabel}[2]{#2\def\@currentlabel{#2}\label{#1}}

\newcommand{\bd}{\boldsymbol{d}}

\newcommand{\btheta}{\boldsymbol{\theta}}

\newcommand{\bbeta}{\boldsymbol{\beta}}

\newcommand{\bTheta}{\boldsymbol{\Theta}}
\newcommand{\bgamma}{\boldsymbol{\gamma}}
\newcommand{\bSigma}{\boldsymbol{\Sigma}}

\newcommand{\bepsilon}{\boldsymbol{\epsilon}}

\newcommand{\bmu}{\boldsymbol{\mu}}

\newcommand{\bC}{\boldsymbol{C}}

\newcommand{\bI}{\boldsymbol{I}}

\newcommand{\bA}{\boldsymbol{A}}

\newcommand{\bu}{\boldsymbol{u}}

\newcommand{\bx}{\bm{x}}
\newcommand{\bX}{\boldsymbol{X}}
\newcommand{\by}{\boldsymbol{y}}

\newcommand{\bz}{\boldsymbol{z}}

\newcommand{\bzero}{\boldsymbol{0}}

\newcommand{\bone}{\boldsymbol{1}}

\newcommand{\postp}{P_{\btheta|\bX_n}}

\newcommand{\pexp}{E_{\btheta|\bX_n}}

\newcommand{\bThetainf}{\mathbf{\Theta}^\infty} 
\newcommand{\mm}{\mathcal{M}}
\newcommand{\dnm}{\delta_{\mathcal{NM}}}
\newcommand{\fdrx}{FDR_{\bX_n}}
\newcommand{\fnrx}{FNR_{\bX_n}}
\newcommand{\mfdr}{mFDR_{\bX_n}}
\newcommand{\mfnr}{mFNR_{\bX_n}}

\DeclareMathOperator*{\argmax}{argmax}

\newtheorem{theorem}{Theorem}

\newtheorem{corollary}[theorem]{Corollary}

\newtheorem{definition}[theorem]{Definition}

\newtheorem{lemma}[theorem]{Lemma}

\newtheorem{remark}[theorem]{Remark}

\newenvironment{proof}[1][Proof]{\textbf{#1.} }{\ \rule{0.5em}{0.5em}}

\newcommand{\bm}{\mathbf}

\numberwithin{equation}{section}
\numberwithin{algo}{section}
\numberwithin{table}{section}
\numberwithin{figure}{section}

\usepackage{setspace}
\usepackage[mathscr]{euscript}
\usepackage[margin=1in]{geometry}
\begin{document}

\normalsize

\title{\vspace{-0.8in}
{\bf High-dimensional Asymptotic Theory of Bayesian Multiple Testing Procedures Under General Dependent Setup and Possible Misspecification}}
\author{Noirrit Kiran Chandra and Sourabh Bhattacharya\thanks{
Noirrit Kiran Chandra is a postdoctoral researcher at Department of Statistical Science, Duke University, USA, and Sourabh Bhattacharya 
is an Associate Professor in Interdisciplinary Statistical Research Unit, Indian Statistical
Institute, 203, B. T. Road, Kolkata 700108.
Corresponding e-mail: noirritchandra@gmail.com.}}
\date{\vspace{-0.5in}}
\maketitle%

\begin{abstract}
In this article, we investigate the asymptotic properties of Bayesian multiple testing procedures under general dependent setup, when
the sample size and the number of hypotheses both tend to infinity. Specifically, we investigate strong consistency of the procedures
and asymptotic properties of different versions of false discovery and false non-discovery rates under the high dimensional setup.
We particularly focus on a novel Bayesian non-marginal multiple testing procedure and its associated error rates in this regard.
Our results show that the asymptotic convergence rates of the error rates are directly associated with the Kullback-Leibler divergence from
the true model, and the results hold even when the postulated class of models is misspecified. 

For illustration of our high-dimensional asymptotic theory, we consider a Bayesian variable selection problem in a time-varying covariate selection framework,
with autoregressive response variables. We particularly focus on the setup where the number of hypotheses
increases at a faster rate compared to the sample size, which is the so-called ultra-high dimensional situation.
\\[2mm]
\\[2mm]
\end{abstract}

\section{Introduction}
\label{sec:introduction}

The area of multiple hypotheses testing has gained much importance and popularity, particularly in this era of big data, where often very large number of hypotheses
need to be tested simultaneously. There are applications abound in the fields of statistical genetics, spatio-temporal statistics, brain imaging, to name a few.
On the theoretical side, it is important to establish validity of the multiple testing procedure in the sense that the procedure controls the false discovery rate $(FDR)$ at
some pre-specified level or attains oracle, as the number of tests grows to infinity. 

Although there is considerable literature addressing these issues, 
the important factor of dependence among the tests seem to have attained less attention. Indeed, realistically, the test statistics or the parameters can not be expected
to be independent. In this regard, \ctn{chandra2019} introduced a novel Bayesian multiple testing procedure that coherently accounts for such dependence
and yields joint decision rules that are functions of appropriate joint posterior probabilities. As demonstrated in \ctn{chandra2019} and \ctn{chandra2020}, 
the new Bayesian method significantly outperforms existing popular multiple testing methods by proper utilization of the dependence structures.
Since in the new method the decisions are taken jointly, the method is referred to as Bayesian non-marginal multiple testing procedure.

\ctn{chandra2020} investigated in detail the asymptotic theory of the non-marginal procedure, and indeed general Bayesian multiple testing
methods under additive loss, for fixed number of hypotheses, when the sample size tends to
infinity. In particular, they provided sufficient conditions for strong consistency of such procedures and also showed that the asymptotic convergence rates
of the versions of $FDR$ and false non-discovery rate ($FNR$) are directly related to the Kullback-Leibler (KL) divergence from the true model.
Interestingly, their results continue to hold even under misspecifications, that is, if the class of postulated models does not include the true model.
In this work, we investigate the asymptotic properties of the Bayesian non-marginal procedure in particular, and Bayesian multiple testing methods
under additive loss in general, when the sample size, as well as the number of hypotheses, tend
to infinity.

As mentioned earlier, asymptotic works in multiple testing when the number of hypotheses grows to infinity, are not rare. 
\ctn{caioptim11} have proposed an asymptotic optimal decision rule for short range dependent data with dependent test statistics. \ctn{bogdan2011} studied the oracle properties and Bayes risk of several multiple testing methods under sparsity in Bayesian decision-theoretic setup. \ctn{datta2013} studied oracle properties for horse-shoe prior when the number of tests grows to infinity. However, in the aforementioned works, the test statistics are independent and follow Gaussian distribution. \ctn{fan2012} proposed a method of dealing with correlated test statistics where the covariance structure is known. Their method is based on principal eigen values of the covariance matrix, which they termed as principal factors. Using those principal factors their method dilutes the association between correlated statistics to deal with an arbitrary dependence structure. They also derived an approximate consistent estimator for the false discovery proportion (FDP) in large-scale multiple testing. \ctn{fan2017} extended this work when the dependence structure is unknown. In these approaches, the decision rules are marginal and the test statistics jointly follow multivariate Gaussian distribution.
\ctn{chandra2019} argue that when the decision rules corresponding to different hypotheses are marginal, the full potential of the dependence structure is not properly
exploited. Results of extensive simulation studies reported in \ctn{chandra2019} and \ctn{chandra2020}, demonstrating superior performance of the Bayesian
non-marginal method compared to popular marginal methods, even for large number of hypotheses, seem to vindicate this issue. This makes asymptotic analysis 
of the Bayesian non-marginal method with increasing number of hypotheses all the more important.

To be more specific, we investigate the asymptotic theory of the Bayesian non-marginal procedure in the general dependence setup, without any particular model assumption,
when the sample size ($n$) and the number of hypotheses ($m_n$, which may be a function of $n$), both tend to infinity. 
We establish strong consistency of the procedure and show that even in this setup, the convergence rates of versions of $FDR$ and $FNR$ are directly related to the KL-divergence
from the true model. We show that our results continue to hold for general Bayesian procedures under the additive loss function. 
In the Bayesian non-marginal context we illustrate the theory with the time-varying covariate 
selection problem, where the number of covariates tends to infinity with the sample size $n$. We distinguish between the two setups: ultra high-dimensional case, that is, where
$\frac{m_n}{n}\rightarrow \infty$ (or some constant), as $n\rightarrow\infty$, and the high-dimensional but not ultra high-dimensional case, that is, 
$m_n\rightarrow\infty$ and $\frac{m_n}{n}\rightarrow 0$, as $n\rightarrow\infty$. We particularly focus on the ultra high-dimensional setup because of its much more challenging nature.

\section{A brief overview of the Bayesian non-marginal procedure}
\label{sec:overview_non_marginal}

Let $\bX_n=\{X_1,\ldots,X_n\}$ denote the available data set. Suppose the data is modelled by the family of distributions $P_{\bX_n|\btheta}$ (which may also be non-parametric). 
For $M>1$, let us denote by $\bTheta=\Theta_1\times\cdots\times\Theta_M$ 
the relevant parameter space associated with $\btheta=(\theta_1,\ldots,\theta_M)$, where we allow $M$ to be infinity as well. 
Let $\postp(\cdot)$ and $\pexp(\cdot)$ denote the posterior distribution and expectation respectively of $\btheta$ given $\bX_n$ and let $P_{\bX_n}(\cdot)$ and $E_{\bX_n} (\cdot)$ denote the marginal distribution and expectation of $\bX_n$ respectively. Let us consider the problem of testing $m$ hypotheses simultaneously corresponding to the actual parameters of interest, where
$1<m\leq M$.  

Without loss of generality, let us consider testing the parameters associated with $\Theta_i$; $i=1,\ldots,m$,
formalized as:
$$ H_{0i}:\theta_i \in \Theta_{0i}  \hbox{ versus } H_{1i} : \theta_i \in \Theta_{1i},$$ 
where $\Theta_{0i} \bigcap \Theta_{1i}=\emptyset \mbox{ and } \Theta_{0i} \bigcup \Theta_{1i} 
= \Theta_{i},\mbox{ for $i=1,\cdots,m$}.$

Let
\begin{align}
d_i=&\begin{cases}
1&\text{if the $i$-th hypothesis is rejected;}\\
0&\text{otherwise;}
\end{cases}\\
r_i=&\begin{cases}
1&\text{if $H_{1i}$ is true;}\\
0&\text{if $H_{0i}$ is true.} 
\end{cases}
\end{align}
Following \ctn{chandra2019} we define $G_i$ to be the set of hypotheses, including the $i$-th one, which are highly dependent, and define
\begin{equation}
z_i=\begin{cases}
1&\mbox{if $H_{d_j,j}$ is true for all $j\in G_i\setminus\{i\}$;}\\
0&\mbox{otherwise.}
\end{cases}\label{eq:z}
\end{equation}

If, for any $i\in\{1,\ldots,m\}$, $G_i=\{i\}$, a singleton, then we define $z_i=1$.
\ctn{chandra2019} maximize the posterior expectation of the number of true positives
\begin{equation}
TP=\sum_{i=1}^md_ir_iz_i,
\label{eq:tp}
\end{equation}
subject to controlling the posterior expectation of the error term
\begin{equation}
E=\sum_{i=1}^md_i(1-r_iz_i),
\label{eq:e}
\end{equation}
which is actually the posterior mean of the sum of three error terms
$E_1=\sum_{i=1}^md_i(1-r_i)z_i$, $E_2=\sum_{i=1}^md_i(1-r_i)(1-z_i)$ and $E_3=\sum_{i=1}^md_ir_i(1-z_i)$.
For detailed discussion regarding these, see \ctn{chandra2019}.

It follows that the decision configuration can be obtained by minimizing the function
\begin{align}
\xi(\bd)&=-\sum_{i=1}^md_i\pexp(r_iz_i)+\lambda_n\sum_{i=1}^md_i \pexp (1-r_iz_i)\notag\\
&= -(1+\lambda_n)\sum_{i=1}^md_i\left(w_{in}(\bd)-\frac{\lambda_n}{1+\lambda_n}\right),\notag
\end{align}
with respect to all possible decision configurations of the form $\bd=\{d_1,\ldots,d_m\}$, where
$\lambda_n>0$,
and
\begin{equation}
w_{in}(\bd)=\pexp(r_iz_i)= \postp\left(H_{1i}\cap\left\{\cap_{j\neq i,j\in G_i}H_{d_j,j}\right\}\right),
\label{eq:w}
\end{equation}
is the posterior probability of the decision configuration $\{d_1,\ldots,d_{i-1},1,d_{i+1},\ldots,d_m\}$
being correct.
Letting $\beta_n=\lambda_n/(1+\lambda_n)$, one can equivalently maximize
\begin{equation}
f_{\beta_n}(\bd)=\sum_{i=1}^m d_i\left(w_{in}(\bd)-\beta_n\right)\label{eq:beta1}
\end{equation}
with respect to $\bd$ and obtain the optimal decision configuration.

\begin{definition}
	Let $\mathbb D$ be the set of all $m$-dimensional binary vectors denoting all possible decision configurations. Define $$\widehat{\bd}=\argmax_{\bd\in\mathbb{D}} f_\beta(\bd)$$ where $0<\beta<1$. Then $\widehat{\bd}$ is the \textit{optimal decision configuration} obtained as the solution of the non-marginal multiple testing method.
	\label{def:nmd}
\end{definition}

For detailed discussion regarding the choice of $G_i$s in (\ref{eq:z}), see \ctn{chandra2019} and \ctn{chandra2020}. In particular, \ctn{chandra2020} show that
asymptotically, the Bayesian non-marginal method is robust with respect to $G_i$s in the sense that it is consistent with respect to any choice of the grouping structure.
As will be shown in this article, the same holds even in the high-dimensional asymptotic setup.

\subsection{Error measures in multiple testing}
\label{subsec:Bayesian_errors}

\ctn{storey03} advocated \textit{positive False Discovery Rate} $(pFDR)$ as a measure of type-I error in multiple testing. Let $\delta_{\mm}(\bd|\bX_n)$ be the probability of choosing $\bd$ as the optimal decision configuration given data $\bX_n$ when a multiple testing method $\mm$ is employed. Then $pFDR$ is defined as:
\begin{equation}
pFDR=E_{\bX_n} \left[ \sum_{\bd\in\mathbb{D}}  
\frac{\sum_{i=1}^{m}d_i(1-r_i)}{\sum_{i=1}^{m}d_i}\delta_\mm(\bd|\bX_n)\bigg{|}\delta_\mm(\bd=\mathbf{0}|\bX_n)=0 \right].
\label{eq:pfdr}
\end{equation}

Analogous to type-II error, the \textit{positive False Non-discovery Rate} $(pFNR)$ is defined as
\begin{align}
pFNR= E_{\bX_n}\left[\sum_{\bd\in\mathbb D} \frac{\sum_{i=1}^m(1-d_i)r_i} {\sum_{i=1}^m(1-d_i)} \delta_{\mathcal M}\left(\bd|\bX_n\right)
\bigg | \delta_{\mathcal M}\left(\bd=\bone|\bX_n\right)=0\right].
\label{eq:pFNR}
\end{align}

Under prior $\pi(\cdot)$, \ctn{SanatGhosh08} defined posterior $FDR$ and $FNR$. The measures are given as following:
\begin{align}
posterior~FDR
&= \pexp\left[\sum_{\bd\in\mathbb D}\frac{\sum_{i=1}^md_i(1-r_i)}{\sum_{i=1}^md_i \vee 1}\delta_{\mathcal M}\left(\bd|\bX_n\right) \right]= \sum_{\bd\in\mathbb{D}} 
\frac{\sum_{i=1}^{m}d_i(1-v_{in})}{\sum_{i=1}^{m}d_i \vee 1}\delta_{\mathcal M}(\bd|\bX_n); \label{eq:pBFDR}\\
posterior~FNR
&= \pexp\left[\sum_{\bd\in\mathbb D} \frac{\sum_{i=1}^m(1-d_i)r_i} {\sum_{i=1}^m(1-d_i)\vee 1} \delta_{\mathcal M}\left(\bd|\bX_n\right)
\right]=\sum_{\bd\in\mathbb D}
\frac{\sum_{i=1}^{m}(1-d_i)v_{in}}{\sum_{i=1}^{m}(1-d_i)\vee 1}\delta_{\mathcal M}(\bd|\bX_n),\label{eq:pBFNR}
\end{align}

where $v_{in}=\postp(\varTheta_{1i})$. Also under any non-randomized decision rule $\mathcal{M}$, $\delta_{\mathcal{M}}(\bd|\bX_n)$ is either 1 or 0 depending on data $\bX_n$. Given $\bX_n$, we denote these posterior error measures by $\fdrx$ and $\fnrx$ respectively.

With respect to the new notions of errors in (\ref{eq:tp}) and (\ref{eq:e}),
\ctn{chandra2019} modified $\fdrx$ as
\begin{align}
modified~\fdrx &=\pexp\left[ \sum_{\bd\in\mathbb D}\frac{\sum_{i=1}^md_i(1-r_iz_i)}{\sum_{i=1}^md_i\vee1}
\delta_{\mathcal M}\left(\bd|\bX_n\right)\right] 
\notag\\
&=  \sum_{\bd\in\mathbb{D}}  
\frac{\sum_{i=1}^{m}d_i(1-w_{in} (\bd))}{\sum_{i=1}^{m}d_i\vee1}\delta_{\mathcal M} (\bd|\bX_n),
\label{eq:mpBFDR}
\end{align}
and $\fnrx$ as
\begin{align}
modified~\fnrx &= \pexp\left[\sum_{\bd\in\mathbb D} \frac{\sum_{i=1}^m(1-d_i)r_iz_i}{\sum_{i=1}^m(1-d_i)\vee1}
\delta_{\mathcal M}\left(\bd|\bX_n\right)
\right]\notag\\
&= \sum_{\bd\in\mathbb D} 
\frac{\sum_{i=1}^{m}(1-d_i)w_{in}(\bd)}{\sum_{i=1}^{m}(1-d_i)\vee1}
\delta_{\mathcal M}(\bd|\bX_n).
\label{eq:mpBFNR}
\end{align}

We denote $modified~\fdrx$ and $\fnrx$ by $\mfdr$ and $\mfnr$ respectively. Notably, the expectations of $\fdrx$ and $\fnrx$ with respect to $\bX_n$, conditioned on the fact that their respective denominators are positive, yields the \textit{positive Bayesian} $FDR~(pBFDR)$ and $FNR$ $(pBFNR)$ respectively. The same expectation over $\mfdr$ and $\mfnr$ yields \textit{modified positive} $BFDR~(mpBFDR)$ and \textit{modified positive} $BFNR~(mpBFNR)$ respectively.



\ctn{muller04} (see also \ctp{sun2009, caioptim11}) considered the following additive loss function
\begin{equation}
L(\bd,\btheta)= c\sum_{i=1}^m d_i(1-r_i)+ \sum_{i=1}^m (1-d_i)r_i,
\label{eq:loss_mul}
\end{equation}
where $c$ is a positive constant.
The decision rule that minimizes the posterior risk of the above loss is 
$d_i=I\left(v_i>\frac{c}{1+c}\right)$ for all $i=1,\cdots,m$, where $I(\cdot)$ is the indicator function. 
%
Observe that the non-marginal method boils down to this additive loss function based approach when $G_i=\{i\}$, that is, when the information regarding dependence between hypotheses is not available or overlooked. Hence, the convergence properties of the additive loss function based methods can be easily derived from our theories. 

Note that multiple testing problems can be regarded as model selection problems where the task is to choose the correct specification for the parameters under consideration. 
The model is misspecified even if one decision is taken incorrectly. Under quite general conditions, \ctn{Shalizi09} investigated asymptotic behaviour of misspecified models. 
We adopt his basic assumptions and some of his convergence results to build a general asymptotic theory for our Bayesian non-marginal multiple testing method in high dimensions. 

In Section \ref{sec:shalizi}, we provide the setup, assumptions and the main result of \ctn{Shalizi09} which we adopt for our purpose. 
In Section \ref{ref:infinite_hypotheses} we address consistency of the Bayesian non-marginal method and convergence of the associated error terms
in the high-dimensional setup. High-dimensional asymptotic analyses of versions of $FDR$ and $FNR$ are detailed in Sections \ref{sec:highdim_FDR}
and \ref{sec:highdim_FNR}, respectively.
In Section \ref{sec:asymp_category_a}, we establish the high-dimensional asymptotic theory for $FNR_{\bX_n}$ and $BFNR$ when versions of $BFDR$ are 
$\alpha$-controlled asymptotically.
We illustrate the asymptotic properties of the non-marginal method in a multiple testing setup associated with an autoregressive model
involving time-varying covariates in Section \ref{sec:ar1_inf}, in high-dimensional contexts.
Finally, in Section \ref{sec:conclusion_inf} we summarize our contributions and provide concluding remarks.


\section{Preliminaries for ensuring posterior convergence under general setup}
\label{sec:shalizi}

Following \ctn{Shalizi09} we consider a probability space $(\Omega,\mathcal F, P)$, 
and a sequence of random variables $X_1,X_2,\ldots$,   
taking values in some measurable space $(\Xi,\mathcal X)$, whose
infinite-dimensional distribution is $P$. The natural filtration of this process is
$\sigma(\bX_n)$.

We denote the distributions of processes adapted to $\sigma(\bX_n)$ 
by $P_{\bX_n|\btheta}$, where $\btheta$ is associated with a measurable
space $(\bTheta,\mathcal T)$, and is generally infinite-dimensional. 
For the sake of convenience, we assume, as in \ctn{Shalizi09}, that $P$
and all the $P_{\bX_n|\btheta}$ are dominated by a common reference measure, with respective
densities $p$ and $f_{\btheta}$. The usual assumptions that $P\in\bTheta$ or even $P$ lies in the support 
of the prior on $\bTheta$, are not required for Shalizi's result, rendering it very general indeed. We put the prior distribution $\pi(\cdot)$ on the parameter space $\bTheta$.

\subsection{Assumptions and theorem of Shalizi}
\label{subsec:assumptions_shalizi}

\begin{enumerate}[label={(S\arabic*)}]	
	\item  \label{shalizi1} Consider the following likelihood ratio:
	\begin{equation}
	R_n(\btheta)=\frac{f_{\btheta}(\bX_n)}{p(\bX_n)}.
	\label{eq:R_n}
	\end{equation}
	Assume that $R_n(\btheta)$ is $\sigma(\bX_n)\times \mathcal T$-measurable for all $n>0$.
	
	\item \label{s2} For each $\btheta\in\Theta$, the generalized or relative asymptotic equipartition property holds, and so,
	almost surely,
	\begin{equation*}
	\underset{n\rightarrow\infty}{\lim}~\frac{1}{n}\log R_n(\btheta)=-h(\btheta),
	\end{equation*}
	where $h(\btheta)$ is given in (S3) below.
	
	\item \label{s3} For every $\btheta\in\Theta$, the KL-divergence rate
	\begin{equation}
	h(\btheta)=\underset{n\rightarrow\infty}{\lim}~\frac{1}{n}E\left(\log\frac{p(\bX_n)}{f_{\btheta}(\bX_n)}\right).
	\label{eq:S3}
	\end{equation}
	exists (possibly being infinite) and is $\mathcal T$-measurable.
	
	\item \label{s4}
	Let $I=\left\{\btheta:h(\btheta)=\infty\right\}$. 
	The prior $\pi$ satisfies $\pi(I)<1$.
	
	Following the notation of \ctn{Shalizi09}, for $A\subseteq\Theta$, let
	\begin{align}
	h\left(A\right)&=\underset{\btheta\in A}{\mbox{ess~inf}}~h(\btheta);\label{eq:h2}\\
	J(\btheta)&=h(\btheta)-h(\Theta);\label{eq:J}\\
	J(A)&=\underset{\btheta\in A}{\mbox{ess~inf}}~J(\btheta).\label{eq:J2}
	\end{align}
	
	\item \label{s5} There exists a sequence of sets $\mathcal G_n\rightarrow\Theta$ as $n\rightarrow\infty$ 
	such that: 
	\begin{enumerate}
		\item[(1)]
		\begin{equation}
		\pi\left(\mathcal G_n\right)\geq 1-\alpha\exp\left(-\varsigma n\right),~\mbox{for some}~\alpha>0,~\varsigma>2h(\Theta);
		\label{eq:S5_1}
		\end{equation}
		\item[(2)]The convergence in (S3) is uniform in $\theta$ over $\mathcal G_n\setminus I$.
		\item[(3)] $h\left(\mathcal G_n\right)\rightarrow h\left(\Theta\right)$, as $n\rightarrow\infty$.
	\end{enumerate}
	
	For each measurable $A\subseteq\Theta$, for every $\delta>0$, there exists a random natural number $\tau(A,\delta)$
	such that
	\begin{equation}
	n^{-1}\log\int_{A}R_n(\btheta)\pi(\btheta)d\btheta
	\leq \delta+\underset{n\rightarrow\infty}{\limsup}~n^{-1}
	\log\int_{A}R_n(\btheta)\pi(\btheta)d\btheta,
	\label{eq:limsup_2}
	\end{equation}
	for all $n>\tau(A,\delta)$, provided 
	$\underset{n\rightarrow\infty}{\lim\sup}~n^{-1}\log\pi\left(\mathbb I_A R_n\right)<\infty$.
	Regarding this, the following assumption has been made by Shalizi:
	
	\item\label{s6} The sets $\mathcal G_n$ of (S5) can be chosen such that for every $\delta>0$, the inequality
	$n>\tau(\mathcal G_n,\delta)$ holds almost surely for all sufficiently large $n$.
	
	\item \label{shalizi7} The sets $\mathcal G_n$ of (S5) and (S6) can be chosen such that for any set $A$ with $\pi(A)>0$, 
	\begin{equation}
	h\left(\mathcal G_n\cap A\right)\rightarrow h\left(A\right),
	\label{eq:S7}
	\end{equation}
	as $n\rightarrow\infty$.
\end{enumerate}
Under the above assumptions, the following version of the theorem of \ctn{Shalizi09} can be seen to hold.

\begin{theorem}[\ctp{Shalizi09}]
	\label{th:shalizi}
	Consider assumptions \ref{shalizi1}--\ref{shalizi7} and any set $A\in\mathcal T$ with $\pi(A)>0$. If $\varsigma>2h(A)$, where
	$\varsigma$ is given in (\ref{eq:S5_1}) under assumption (S5), then
	\begin{equation}
	\underset{n\rightarrow\infty}{\lim}~\frac{1}{n}\log\postp(A|\bX_n)=-J(A).
	\label{eq:post_conv1}
	\end{equation}
\end{theorem}
We shall frequently make use of this theorem for our purpose. Also throughout this article, we show consistency results for general models which satisfy \ref{shalizi1}--\ref{shalizi7}. For all our results, we assume these conditions.

\section{Consistency of multiple testing procedures when the number of hypotheses tends to infinity}
\label{ref:infinite_hypotheses}

In this section we show that the non-marginal procedure is asymptotically consistent under any general dependency model satisfying the conditions 
in Section \ref{subsec:assumptions_shalizi}. Since one of our main goals is to allow for misspecification, we must define consistency of multiple testing
methods encompassing misspecification, while also allowing for $m_n$ hypotheses where $m_n/n\rightarrow c$, where $c\geq 0$ or $\infty$. 
We formalize this below by introducing appropriate notions.

\subsection{Consistency of multiple testing procedures under misspecification}
\label{cons_miss}

Let $\bThetainf$ be the infinite dimensional parameter space of the countably infinite set of parameters $\{\theta_1,\theta_2,\ldots\}$. In this case, 
any decision configuration $\bd$ is also an infinite dimensional vector of 0's and 1's. 
Define $\bTheta^t=\otimes_{i=1}^\infty \bTheta_{d_i^t,i}$, where $``\otimes"$ denotes cartesian product, and $\bd^t=(d^t_1,d^t_2,\ldots)$ denotes the actual infinite dimensional
decision configuration satisfying $J\left(\bTheta^t\right)=J\left(\bThetainf\right)$. This definition of $\bd^t$ accounts for misspecification in the sense
that $\bd^t$ is the minimizer of the KL-divergence from the true data-generating model.
For any decision $\bd$, let $\bd(m_n)$ denote the first $m_n$ components of $\bd$. Let $\mathbb D_{m_n}$ denote the set of all possible decision configurations
corresponding to $m_n$ hypotheses.
With the aforementioned notions, we now define consistency of multiple testing procedures.

\begin{definition}
	Let $\bd^t(m_n)$ be the true decision configuration among all possible decision configurations in $\mathbb D_{m_n}$. 
	Then a multiple testing method $\mm$ is said to be asymptotically consistent if almost surely
	\begin{equation}
		\lim_{n\rightarrow\infty} \delta_\mm(\bd^t(m_n)|\bX_n) =1.
	\end{equation}
	\label{def:consistent}
\end{definition}
Recall the constant $\beta_n$ in (\ref{eq:beta1}), which is the penalizing constant between the error $E$ and true positives $TP$. 
For consistency of the non-marginal procedure, we need certain conditions on $\beta_n$, which we state below. These conditions will also play important roles in the asymptotic studies
of the different versions of $FDR$ and $FNR$ that we consider.
\begin{enumerate}[label={(A\arabic*)}]	
	\item \label{A1} We assume that the sequence $\beta_n$ is neither too small nor too large, that is,
	\begin{align}
	\underline\beta&=\underset{n\geq 1}{\liminf}\beta_n>0;\label{eq:liminf_beta}\\
	\overline\beta&=\underset{n\geq 1}{\limsup}\beta_n<1.\label{eq:limsup_beta}
	\end{align}
	\item \label{A2} We assume that neither all the null hypotheses are true and nor all of then are false for $m_n$ hypotheses being considered, 
		that is, $\bd^t(m_n)\neq\bzero$ and $\bd^t(m_n)\neq\bone$, where
	$\bzero$ and $\bone$ are vectors of 0's and 1's respectively.
\end{enumerate}
Condition \ref{A1} is necessary for the asymptotic consistency of both the non-marginal method and additive loss function based method. This ensures that the penalizing constant is asymptotically bounded away from 0 and 1, that is, it is neither too small nor too large. Notably, \ref{A2} is not required for the consistency results. The role of \ref{A2} is to ensure that the denominator terms in the multiple testing error measures (defined in Section \ref{subsec:Bayesian_errors}) do not become 0. 

\subsection{Main results on consistency in the infinite-dimensional setup}
\label{subsec:main_consistency}

In this section we investigate the asymptotic properties of the Bayesian non-marginal method and \ctn{muller04} when $m_n/n$ tends to infinity or some
positive constant. It is to be noted that result (\ref{eq:post_conv1}) of \ctn{Shalizi09} 
holds even for infinite-dimensional parameter space. Exploiting this fact we derive the results in this section. 

Note that if there exists a value $\btheta^t$ of $\btheta$ that minimizes the KL-divergence, then $\btheta^t$ is in the set $\bTheta^t$. 
Let us denote by $\bTheta^{tc}$ the complement of $\bTheta^t$. Observe that if $\btheta^t$ lies in the interior of $\bTheta^t$, then $J\left(\bTheta^{tc}\right)>0$.
It then holds that
\begin{equation}
	\lim_{n\rightarrow\infty} \frac{1}{n}\log \postp \left( \bTheta^{tc}\right) = -J\left( \bTheta^{tc}\right),
\end{equation}
which implies that for any $\epsilon>0$, there exists a $n_0(\epsilon)$ such that for all $n>n_0(\epsilon)$
\begin{align}
	&\exp\left[-n\left(J\left(\bTheta^{tc}\right)+\epsilon\right)\right]
	< \postp \left( \bTheta^{tc}\right) < \exp \left [-n\left(J\left(\bTheta^{tc}\right)-\epsilon\right) \right]	\label{eq:n01}\\
\Rightarrow& 1-\exp \left [-n\left(J\left(\bTheta^{tc}\right)-\epsilon\right) \right] < \postp \left( \bTheta^t\right) <1- \exp\left[-n\left(J\left(\bTheta^{tc}\right)+\epsilon\right)\right].\label{eq:n02}
\end{align}
For notational convenience, we shall henceforth denote $J\left(\bTheta^{tc}\right)$ by $J$.


Note that the groups $G_i$ also depend upon $m_n$ in our setup; hence, we denote them by $G_{i,m_n}$.
For any decision configuration $\bd(m_n)$ and group $G_{m_n}$ let $\bd_{G_{m_n}}=\{d_j:j\in G_{m_n}\}$. Define
$$\mathbb D_{i,m_n}=\left\{\bd(m_n):~\mbox{all decisions in}~\bd_{G_{i,m_n}}~\mbox{are correct}\right\}.$$
Here $\mathbb D_{i,m_n}$ is the set of all decision configurations where the decisions corresponding to the hypotheses in $G_{i,m_n}$ are at least correct. 
Clearly $\mathbb D_{i,m_n}$ contains $\bd^t(m_n)$ for all $i=1,2,\ldots,m_n$.

Hence, $\mathbb D^c_{i,m_n}=\left\{\bd(m_n):~\mbox{at least one decision in} ~\bd_{G_{i,m_n}}~ \mbox{is incorrect}\right\}$.
Observe that if $\bd(m_n)\in\mathbb D^c_{i,m_n}$, at least one decision is wrong corresponding to some parameter in $G_{i,m_n}$. 
As $\postp \left( \bTheta^{tc} \right)$ is the posterior 
probability of at least one wrong decision in the infinite dimensional parameter space, we have
\begin{equation}
	w_{in}(\bd(m_n))\leq w_{in}(\bd)<\postp \left( \bTheta^{tc}\right) < \exp \left [-n\left(J-\epsilon\right) \right].
	\label{eq:ineq_start}
\end{equation}
Also if $H_{0i}$ is true, then
\begin{equation}
	v_{in}\leq w_{in}(\bd)<\postp \left( \bTheta^{tc}\right) < \exp \left [-n\left(J-\epsilon\right) \right].
	\label{eq:ineq_end0}
\end{equation}
Similarly for $\bd(m_n)\in\mathbb D_{i,m_n}$ and for false $H_{0i}$
\begin{equation}
	w_{in}(\bd(m_n))\geq w_{in}(\bd^t)>\postp \left( \bTheta^t\right) >1- \exp \left [-n\left(J-\epsilon\right) \right];
\label{eq:ineq2}
\end{equation}
\begin{equation}
	v_{in}\geq w_{in}(\bd^t)>\postp \left( \bTheta^t\right) >1- \exp \left [-n\left(J-\epsilon\right) \right].
	\label{eq:ineq_end}
\end{equation}
It is important to note that the inequalities (\ref{eq:ineq_start})-(\ref{eq:ineq_end}) hold for all $n>n_0$ and this $n_0$ is the same for all $i$, thanks to validity of 
Shalizi's result in infinite dimensional parameter space.
Exploiting the properties of Shalizi's theorem we will now establish consistency of the Bayesian non-marginal method for increasing number of hypotheses.
\begin{theorem}
	Let $\delta_{\mathcal {NM}}$ denote the decision rule corresponding to the Bayesian non-marginal procedure for $m_n$ hypotheses being tested
	using samples of size $n$, where $m_n\rightarrow\infty$ as $n\rightarrow\infty$. Assume Shalizi's conditions and assumption \ref{A1}.
	Also assume that 
	$J\left(\bTheta^t\right)>0$.
	Then, 
	\begin{align}
		&\lim_{n\rightarrow\infty} \delta_{\mathcal {NM}}(\bd^t(m_n)|\bX_n)=1,~\mbox{almost surely, and}
		\label{eq:inf_dim_cons1}\\
		& \lim_{n\rightarrow\infty} E\left[\delta_{\mathcal {NM}}(\bd^t(m_n)|\bX_n)\right]=1.
		\label{eq:inf_dim_cons2}
	\end{align}
	\label{theorem:inf_dim_cons}
\end{theorem}

\begin{corollary}
	\label{corr:corr1}
	Assuming condition \ref{A1}, the optimal decision rule corresponding to the additive loss function (\ref{eq:loss_mul}) is asymptotically consistent.
	The proof follows in the same way as that of Theorem \ref{theorem:inf_dim_cons} using (\ref{eq:ineq_end0}) and (\ref{eq:ineq_end}).
\end{corollary}

\begin{remark}
Note that Theorem \ref{theorem:inf_dim_cons} does not require any condition regarding the growth of $m_n$ with respect to $n$, and holds
if $m_n/n\rightarrow c$ as $n\rightarrow\infty$, where $c\geq 0$ is some constant, or infinity. Thus, the result seems to be extremely satisfactory.
However, restrictions on the growth of $m_n$ needs to be generally imposed to satisfy the conditions of Shalizi.
An illustration in this regard is provided in Section \ref{sec:ar1_inf}.
\end{remark}

\section{High-dimensional asymptotic analyses of versions of $FDR$}
\label{sec:highdim_FDR}
For a fixed number of hypotheses $m$, \ctn{chandra2020} investigated convergence of different versions of $FDR$ as the sample size $n$ tends to infinity.
They show that show that the convergence rates of the posterior error measures $ mFDR_{\bX_n}$ and $FDR_{\bX_n}$ are directly associated with the KL-divergence
from the true model. Indeed, they were able to obtain the exact limits of $\frac{1}{n}\log mFDR_{\bX_n}$ and $\frac{1}{n}\log FDR_{\bX_n}$ in terms of
the relevant $m$-dimensional KL-divergence rate. 

In the current high-dimensional setup, however, such exact
KL-divergence rate can not be expected to be available since the number of hypotheses $m_n$ is not fixed. As $m_n\rightarrow\infty$, it is plausible 
to expect that the convergence rates depend upon the infinite-dimensional KL-divergence $J$. We show that this is indeed the case, but the exact limit
is not available, which is again to be expected, since $m_n$ approaches infinity, not equal to infinity. Here, in the high-dimensional setup we obtain $-J$ 
as an upper bound of the limit supremums. It is easy to observe that the limits in the finite-dimensional setup are bounded above by $-J$, thus providing evidence
of internal consistency as we move from fixed-dimensional to high-dimensional setup.

We also show that $mpBFDR$ and $pBFDR$ approach zero, even though the rates of convergence are not available. Recall that even in the fixed-dimensional setup,
the convergence rates of $mpBFDR$ and $pBFDR$ were not available. As in the consistency result, these results too do not require any restriction on the growth rate
of $m_n$, except that required for Shalizi's conditions to hold.

We present our results below, the proofs of which are presented in the supplement.
\begin{theorem}
	\label{theorem:mpBFDR_bound_inf}
	Assume the setup and conditions of Theorem \ref{theorem:inf_dim_cons}.
	Then, for any $\epsilon>0$, there exists $n_0(\epsilon)\geq 1$ such that for 
	$n\geq n_0(\epsilon)$, the following hold almost surely: 
	\begin{align}
		mFDR_{\bX_n}
		&\leq e^ {-n(J-\epsilon ) };	
	\label{eq:inf_dim_mpBFDR}\\
		FDR_{\bX_n}
		&\leq e^ {-n(J-\epsilon ) }.	
	\label{eq:inf_dim_pBFDR}
	\end{align}	
\end{theorem}

The above theorem shows that the convergence rate of $\mfdr$ and $FDR$ to 0 for arbitrarily large number of hypotheses is at exponential rate, for arbitrary growth
rate of $m_n$ with respect to $n$. However, again Shalizi's conditions would require restriction on the growth rate of $m_n$.

\begin{corollary}
	Under the setup and assumptions of Theorem \ref{theorem:inf_dim_cons},
	\begin{align}
		\limsup_{n\rightarrow\infty} \frac{1}{n}\log mFDR_{\bX_n}&\leq -J;
		\label{eq:loginf_dim_mpBFDR}\\
		\limsup_{n\rightarrow\infty} \frac{1}{n}\log FDR_{\bX_n}&\leq -J;
		\label{eq:log_inf_dim_pBFDR}
	\end{align}
\end{corollary}

\begin{theorem}
	\label{theorem:lim_BFDR}
		Assume the setup and conditions of Theorem \ref{theorem:inf_dim_cons}, along with assumption \ref{A2}. Then
	\begin{align}
		\lim_{n\rightarrow\infty}mpBFDR &= 0;
		\label{eq:lim_mpBFDR}\\
		\lim_{n\rightarrow\infty}pBFDR &= 0.
		\label{eq:lim_pBFDR}
	\end{align}
\end{theorem}

\section{High-dimensional asymptotic analyses of versions of $FNR$}
\label{sec:highdim_FNR}
High-dimensional asymptotic treatments of versions of $FNR$ are similar to those for versions of $FDR$. In particular, limit supremums of both 
$\frac{1}{n}\log mFNR_{\bX_n}$ and $\frac{1}{n}\log FNR_{\bX_n}$ are bounded above by $-J$, and that both $mpBFNR$ and $pBFNR$ converge to zero.  
The proofs of these results are also similar to those for the respective $FDR$ versions. Internal consistency of these results is again evident
as the limits of $\frac{1}{n}\log mFNR_{\bX_n}$ and $\frac{1}{n}\log FNR_{\bX_n}$ in the finite dimensional setups are bounded above by $-J$ and
$mpBFNR$ and $pBFNR$ converge to zero for fixed number of hypotheses. In the latter cases, convergence rates are not available for either fixed or
high-dimensional cases. Below we provide the relevant results on versions of $FNR$, with proofs in the supplement. 
\begin{theorem}
	\label{theorem:mpBFNR_bound_inf}
	Assume the setup and conditions of Theorem \ref{theorem:inf_dim_cons}.
	Then, for any $\epsilon>0$, there exists $n_0(\epsilon)\geq 1$ such that for 
	$n\geq n_0(\epsilon)$, the following hold almost surely: 
	\begin{align}
		mFNR_{\bX_n}
		&\leq e^ {-n(J-\epsilon ) };	
	\label{eq:inf_dim_mpBFNR}\\
		FNR_{\bX_n}
		&\leq e^ {-n(J-\epsilon ) }.	
	\label{eq:inf_dim_pBFNR}
	\end{align}	
\end{theorem}

The above theorem shows that the convergence rate of $\mfnr$ and $FNR$ to 0 for arbitrarily large number of hypotheses is at exponential rate, for arbitrary growth
rate of $m_n$ with respect to $n$. However, again Shalizi's conditions would require restriction on the growth rate of $m_n$.

\begin{corollary}
	Under the setup and assumptions of Theorem \ref{theorem:inf_dim_cons},
	\begin{align}
		\limsup_{n\rightarrow\infty} \frac{1}{n}\log mFNR_{\bX_n}&\leq -J;
		\label{eq:loginf_dim_mpBFNR}\\
		\limsup_{n\rightarrow\infty} \frac{1}{n}\log FNR_{\bX_n}&\leq -J;
		\label{eq:log_inf_dim_pBFNR}
	\end{align}
\end{corollary}

\begin{theorem}
	\label{theorem:lim_BFNR}
		Assume the setup and conditions of Theorem \ref{theorem:inf_dim_cons}, along with assumption \ref{A2}. Then
	\begin{align}
		\lim_{n\rightarrow\infty}mpBFNR &= 0;
		\label{eq:lim_mpBFNR}\\
		\lim_{n\rightarrow\infty}pBFNR &= 0.
		\label{eq:lim_pBFNR}
	\end{align}
\end{theorem}


\section{High-dimensional asymptotics for $FNR_{\bX_n}$ and $BFNR$ when versions of $BFDR$ are $\alpha$-controlled}
\label{sec:asymp_category_a}

It has been proved in \ctn{chandra2019} for the non-marginal multiple testing procedure and additive loss-function based methods, 
$mpBFDR$ and $pBFDR$ are continuous and non-increasing in $\beta$. Consequently, for suitable values of $\beta$ any $\alpha\in (0,1)$
can be achieved by these errors. For suitably chosen positive values of $\alpha$, one can hope to reduce the corresponding $BFNR$. 
This is standard practice even in the single hypothesis testing literature, where the Type-I error is controlled at some positive value so that a reduced Type-II error
may be incurred. However, as shown in \ctn{chandra2020} in the fixed-dimensional setup, 
for the non-marginal multiple testing procedure and additive loss-function based methods, values of $\alpha$ that are as close to 1 as desired, can not be attained
by versions of $FDR$ as the sample size $n$ tends to infinity. This is not surprising, however, since consistent procedures are not expected to incur large errors
asymptotically, at least when the number of hypothesis is fixed. Indeed, in the fixed-dimensional setup, 
\ctn{chandra2020} provided an interval of the form $(a,b)$ where $0<a<b<1$, in which maximum values of the versions of $FDR$ can lie asymptotically
and obtained asymptotic results for $FNR$ for such $\alpha$-controlled versions of $FDR$.

In this section we investigate the asymptotic theory for $\alpha$-control in the high-dimensional context, that is, when $m_n\rightarrow\infty$ as $n\rightarrow\infty$.
Although none of our previous high-dimensional results did not require any explicit restrictions on the growth rate of $m_n$ given that the posterior convergence
result of Shalizi holds, here we need a very mild condition on $m_n$ that it grows slower than the exponential rate in $n$. We also need to fix the proportion ($p$) of
true alternatives as $m_n\rightarrow\infty$, and the proportion ($q$) of groups associated with at least one false null hypothesis. As we show, these two proportions
define an interval of the form $(0,b)$, with $b=\frac{1-q}{1+p-q}<1$, in which the maximum of the versions of $FDR$ lie, as $m_n\rightarrow\infty$ with $n$. In contrast with the
fixed-dimensional asymptotics of \ctn{chandra2020}, the lower bound of the interval is zero for high dimensions, not strictly positive. To explain, for fixed dimension $m$,
the lower bound was $a=\frac{1}{\sum_{i=1}^md^t_i+1}$. Intuitively, replacing $a$ and $m$ with $a_{m_n}$ and $m_n$ respectively, 
dividing both numerator and denominator of $a$ by $m_n$, taking the limit, replacing the denominator with $p$, we obtain $a_{m_n}\rightarrow 0$, as $n\rightarrow\infty$.
Similar intuition can be used to verify that the upper bound $b$ in the fixed dimensional case converges to $\frac{1-q}{1+p-q}$ in the high-dimensional setup.
As in our previous results, these provide a verification of internal consistency in the case of transition from fixed-dimensional to high-dimensional situations.

Our results regarding asymptotic $\alpha$ control of versions of $FDR$ and corresponding convergence of versions of $FNR$ are detailed in Sections \ref{subsec:constant_BFDR}
and \ref{subsec:highdim_fnr}.

\subsection{High-dimensional $\alpha$-control of $mpBFDR$ and $pBFDR$ for the non-marginal method}
\label{subsec:constant_BFDR}

The following theorem provides the interval for the maximum $mpBFDR$ that can be incurred asymptotically in the high-dimensional setup. 
\begin{theorem}
	\label{theorem:mpBFDR_alpha}
	In addition to \ref{A1}-\ref{A2}, assume the following: 
	\begin{enumerate}[label={(B)}]	
		\item \label{ass_groups} 
			For each $n>1$, let each group of a particular set of $m_{1n}~(<m_n)$ groups out of the total $m_n$ groups be associated with at least one false null hypothesis, and that
			all the null hypotheses associated with the remaining $m_n-m_{1n}$ groups be true. Let us further assume that the latter $m_n-m_{1n}$ groups do not have any overlap with the remaining $m_{1n}$ groups. 
			Without loss of generality assume that $G_{1n},\ldots,G_{m_{1n}}$ are the groups each consisting of at least one false null and 
			$G_{m_{1n}+1},G_{m_{1n}+2},\cdots,G_{m_n}$ are the groups
		where all the null hypotheses are true. 	
			Assume further, the following limits:
			\begin{align}
				&\lim_{n\rightarrow\infty}\frac{m_{1n}}{m_n}=q\in(0,1);\label{eq:q}\\
				&\lim_{n\rightarrow\infty}\frac{\sum_{i=1}^{m_{n}}d^t_i}{m_n}=p\in(0,1);\label{eq:p}\\
				&\lim_{n\rightarrow\infty}m_ne^{-nc}=0~\mbox{for all}~c>0.\label{eq:infdim_limit}
			\end{align}
	\end{enumerate}
	Then the maximum $mpBFDR$ that can be incurred, asymptotically lies in $\left(0, \frac{1-q}{1+p-q} \right)$.
\end{theorem}

\begin{remark}
	If $p$ is close to zero, that is, if all but a finite number of null hypotheses are true, 
	then $\frac{1-q}{1+p-q}\approx 1$, showing that in such cases, better $\alpha$-control can be exercised.
	Indeed, as the proof of the theorem shows, the optimal decision in this case will be given by all but a finite set of one's, so that all but a finite number
	of decisions are correct. Hence, maximum error occurs in this case.
	Also, if $q$ is close to $1$, then $\frac{1-q}{1+p-q}\approx 0$. In other words, if all but a finite number of groups are associated with at least one
	false null hypothesis, then almost no error can be incurred. As the proof Theorem \ref{theorem:mpBFDR_alpha} shows, this is the case where
	all but a finite number of decisions are correct, and hence, it is not surprising that almost no error can be incurred in this case.
\end{remark}

\begin{remark}
	\label{remark:single_G}
	Also, as in the fixed-dimensional case, Theorem \ref{theorem:mpBFDR_alpha} holds, if for at least one $i\in\{1,\ldots,m_n\}$, 
	$G_i\subset\{1,\ldots,m_n\}$. But if
	$G_i=\{1,\ldots,m_n\}$ for $i=1,\ldots,m_n$, then $mpBFDR\rightarrow 0$ as $n\rightarrow\infty$, for any
	sequence $\beta_n\in[0,1]$. 
\end{remark}

\begin{remark}
	Note that in the same way as in the fixed-dimensional setup, Theorem \ref{theorem:mpBFDR_alpha} remains valid even for $mFDR_{\bX_n}$ 
	thanks to its monotonicity with respect to $\beta$, the property
	crucially used to prove Theorem \ref{theorem:mpBFDR_alpha}.
\end{remark}
The following theorem shows that for feasible values of $\alpha$ attained asymptotically by the maximum of $mpBFDR$, for appropriate sequences of penalizing
constants $\beta_n$, it is possible to asymptotically approach such $\alpha$ through $mpBFDR_{\beta_n}$, where 
$mpBFDR_\beta$ denotes $mpBFDR$ for the non-marginal procedure where the penalizing constant is $\beta$.
\begin{theorem}
	\label{corollary:beta_n0}
	
	Suppose that
	\begin{align}
	\lim_{n\rightarrow\infty} mpBFDR_{\beta=0}= E.\label{eq:lim_mpbfdr}
	\end{align}
	Then, for any $\alpha<E$ and $\alpha\in\left(0,\frac{1-q}{1+p-q} \right)$, under condition \ref{ass_groups},
	there exists a sequence $\beta_n\rightarrow 0$ such that $mpBFDR_{\beta_n}\rightarrow\alpha$ as $n\rightarrow\infty$.
\end{theorem}

From the proofs of Theorem \ref{theorem:mpBFDR_alpha} and \ref{corollary:beta_n0}, it can be seen that replacing $w_{in}(\hat \bd(m_n))$ 
by $v_{in}$ does not affect the results. Hence we state the following corollary.
\begin{corollary}
	\label{theorem:pBFDR_alpha}
	
	Let $pBFDR_\beta$ denote the $pBFDR$ corresponding to the non-marginal procedure where the penalizing constant is $\beta$. Suppose that
	\begin{align*}
	\lim_{n\rightarrow\infty} pBFDR_{\beta=0}= E',
	\end{align*}
	Then, for any $\alpha<E'$ and $\alpha\in\left(0,\frac{1-q}{1+p-q}\right)$, under condition \ref{ass_groups},
	there exists a sequence $\beta_n\rightarrow 0$ such that $pBFDR_{\beta_n}\rightarrow\alpha$ as $n\rightarrow\infty$.
\end{corollary}

As in the fixed-dimensional setup, we see that for $\alpha$-control we must have $\lim_{n\rightarrow\infty}\beta_n=0$, and that for $\liminf_{n\rightarrow\infty}\beta_n>0$,
$mpBFDR$ tends to zero. In other words, even in the high-dimensional setup, $\alpha$-control requires a sequence $\beta_n$ that is smaller that that 
for which $mpBFDR$ tends to zero.


Since the additive loss function based methods are special cases of the non-marginal procedure where $G_i =\{i\}$ for all $i$ (see \ctn{chandra2019}, \ctn{chandra2020}), 
and that in such cases, $mpBFDR$ reduces to $pBFDR$, it is important to investigate
asymptotic $\alpha$-control of $pBFDR$ in this situation. Our result in this direction is provided in Theorem \ref{theorem:pBFDR_alpha2}.
\begin{theorem}
	\label{theorem:pBFDR_alpha2}
	Let $m_{0n}~(<m_n)$ be the number of true null hypotheses such that $m_{0n}/m_n\rightarrow p_0\in(0,1)$, as $n\rightarrow\infty$. 
	Then for any $0<\alpha<p_0$, there exists a sequence $\beta_n\rightarrow 0$
	as $n\rightarrow\infty$ such that for the additive loss function based methods
	\begin{equation*}
	\underset{n\rightarrow\infty}{\lim}~pBFDR_{\beta_n} =\alpha.
	\end{equation*}
\end{theorem}
The result is similar in spirit to that obtained by \ctn{chandra2020} in the corresponding finite dimensional situation. 
The limit of $m_{0n}/m_n$ in the corresponding high-dimensional setup, instead of $m_0/m$ in the fixed dimensional case, plays the central role here.

\ctn{chandra2019} and \ctn{chandra2020} noted that even for additive loss function based multiple testing procedures, $mpBFDR$ may be a more desirable candidate compared to
$pBFDR$ since it can yield non-marginal decisions even if the multiple testing criterion to be optimized is a simple sum of loss functions designed to yield marginal decisions.
The following theorem shows that the same high-dimensional asymptotic result as 
Theorem \ref{theorem:pBFDR_alpha2} also holds for $mpBFDR$ in the case of additive loss functions, without the requirement of condition \ref{ass_groups}.
Non-requirement of condition \ref{ass_groups} even in the high-dimensional setup can be attributed to the fact that 
$mpBFDR (\mathcal M)\geq pBFDR (\mathcal M)$ for any multiple testing method $\mathcal{M}$, for arbitrary sample size.
\begin{theorem}
	\label{theorem:mpBFDR_alpha2}
	Let $m_{0n}~(<m_n)$ be the number of true null hypotheses such that $m_{0n}/m_n\rightarrow p_0\in(0,1)$, as $n\rightarrow\infty$. 
	Let $\alpha$ be the desired level of significance where $0<\alpha < p_0$. Then there exists a sequence $\beta_n\rightarrow 0$
	as $n\rightarrow\infty$ such that for the additive loss function based method
	\begin{equation*}
	\underset{n\rightarrow\infty}{\lim}~mpBFDR_{\beta_n}=\alpha.
	\end{equation*}
\end{theorem}

Note that Bayesian versions of $FDR$ (conditional on the data) need not be continuous with respect to $\beta$, and so results for such Bayesian versions
similar to Theorem \ref{corollary:beta_n0}, Corollary \ref{theorem:pBFDR_alpha} and Theorems \ref{theorem:pBFDR_alpha2}, \ref{theorem:mpBFDR_alpha2}, which heavily
use such continuity property, could not be established. 


Thus, interestingly, all the asymptotic results for $\alpha$-control of versions of $FDR$ in the fixed dimensional setup admitted simple extensions to the high-dimensional
setup, with minimal assumption regarding the growth rate of $m_n$, given Shalizi's conditions hold. Since Shalizi's conditions are meant for posterior consistency,
from the multiple testing perspective, our high-dimensional results are very interesting in the sense that almost no extra assumptions are required in addition
to Shalizi's conditions for our multiple testing results to carry over from fixed dimension to high dimensions.

\subsection{High-dimensional properties of Type-II errors  when $mpBFDR$ and $pBFDR$ are asymptotically controlled at $\alpha$}
\label{subsec:highdim_fnr}
In this section, we investigate the high-dimensional asymptotic theory for $\fnrx$ and $pBFNR$ associated with $\alpha$-control of versions of $FDR$.
Our results in these regards are provided as Theorem \ref{theorem:mpBFNR_alpha} and Corollary \ref{corr:corr2}.
\begin{theorem}
	\label{theorem:mpBFNR_alpha}
	Assume condition \ref{ass_groups} and that $n^{-1}\log m_n\rightarrow 0$, as $n\rightarrow\infty$. 
	Then for asymptotic $\alpha$-control of $mpBFDR$ in the non-marginal procedure the following holds almost surely:
	\begin{equation*}
	\limsup_{n\rightarrow\infty} \fnrx \leq -J.
	\end{equation*}
\end{theorem}
The above theorem requires the very mild assumption that $n^{-1}\log m_n\rightarrow 0$, as $n\rightarrow\infty$, in addition to \ref{ass_groups}.
The result shows that $FNR_{\bX_n}$ converges to zero at an exponential rate, but again the exact limit of $\fnrx$ is not available in this high-dimensional setup.
This is slightly disconcerting in the sense that we are now unable to compare the rates of convergence of $FNR_{\bX_n}$ for cases where $\alpha$-control is imposed
and not imposed. Indeed, for the fixed-dimensional setup, \ctn{chandra2020} could obtain exact limits and consequently show that 
$\fnrx$ converges to zero at a rate faster than or equal to that compared to the case when $\alpha$ control is not exercised. However, as we already argued
in the context of versions of $FDR$, exact limits are not expected to be available in these cases for high dimensions.

\begin{corollary}
	\label{corr:corr2}
	Assume condition \ref{ass_groups} and that $n^{-1}\log m_n\rightarrow 0$, as $n\rightarrow\infty$. 
	Then for asymptotic $\alpha$-control of $mpBFDR$ in the non-marginal procedure the following holds: 
	\begin{equation*}
	\lim_{n\rightarrow\infty} pBFNR =0.
	\end{equation*}
\end{corollary}

Thus, as in the fixed dimensional setup, Corollary \ref{corr:corr2} shows that corresponding to $\alpha$-control, 
$pBFNR$ converges to zero even in the high-dimensional setup, and that
the rate of convergence to zero is unavailable.

\section{Illustration of consistency of our non-marginal multiple testing procedure in time-varying covariate selection in autoregressive process}
\label{sec:ar1_inf}
Let the true model $P$ stand for the following $AR(1)$ model consisting of time-varying covariates: 
\begin{equation}
x_t=\rho_0 x_{t-1}+\sum_{i=0}^m\beta_{i0}z_{it}+\epsilon_t,~t=1,2,\ldots,n,
\label{eq:true_ar1_inf}
\end{equation}
where $x_0\equiv 0$, $|\rho_0|<1$ and $\epsilon_t\stackrel{iid}{\sim}N(0,\sigma^2_0)$, for $t=1,2,\ldots,n$. In (\ref{eq:true_ar1_inf}), $m\equiv m_n\rightarrow\infty$
as $n\rightarrow\infty$. Here $\left\{z_{it}:t=1,2,\ldots\right\}$ are relevant time-varying covariates.
We set $z_{0t}\equiv 1$ for all $t$.

Now let the data be modeled by the same model as $P$ but with $\rho_0$, $\beta_{i0}$ and $\sigma^2_0$ be replaced with the unknown quantities $\rho$, $\beta_i$
and $\sigma^2$, respectively, that is,
\begin{equation}
x_t=\rho x_{t-1}+\sum_{i=0}^m\beta_iz_{it}+\epsilon_t,~t=1,2,\ldots,n,
\label{eq:modeled_ar1_inf}
\end{equation}
where we set $x_0\equiv 0$, $\epsilon_t\stackrel{iid}{\sim}N(0,\sigma^2)$, for $t=1,2,\ldots,n$.

For notational purposes, we let $\bz_{mt}=(z_{0t},z_{1t},\ldots,z_{mt})^\prime$, $\bz_t=(z_{0t},z_{1t},\ldots)^\prime$, $\bbeta_{m0}=(\beta_{00},\beta_{10},\ldots,\beta_{m0})^\prime$, 
$\bbeta_m=(\beta_0,\beta_1,\ldots,\beta_m)^\prime$ and $\bbeta=(\beta_0,\beta_1,\ldots)^\prime$. 

\subsection{The ultra high-dimensional setup}
\label{subsec:ultra}
Let us first consider the setup where $\frac{m_n}{n}\rightarrow\infty$ as $n\rightarrow\infty$. This is a challenging problem, and we require notions of sparsity
to address such a problem. As will be shown subsequently in Section \ref{subsec:remarks_assumptions}, a precise notion of sparsity is available for our problem
in the context of the equipartition property. Specifically sparsity in our problem entails controlling relevant quadratic forms of $\bbeta$. For such sparsity,
we must devise a prior for $\bbeta$ such that $\|\bbeta\|<\infty$. We also assume that $\|\bbeta_0\|<\infty$.

For appropriate prior structures for $\bbeta$, let us consider the following strategy.
First, let us consider an almost surely continuously differentiable random function $\tilde\eta(\cdot)$ on a compact space $\mathcal X$, such that 
\begin{equation}
\|\tilde\eta\|=\underset{\tilde\bx\in\mathcal X}{\sup}~|\tilde\eta(\tilde\bx)|<\infty,~\mbox{almost surely.} 
\label{eq:tilde_eta}
\end{equation}
We denote the class of such functions as $\mathcal C'(\mathcal X)$. A popular prior for $\mathcal C'(\mathcal X)$ is the Gaussian process prior with sufficiently smooth
covariance function, in which case, both $\tilde\eta$ and $\tilde\eta'$ are Gaussian processes; see, for example, \ctn{Cramer67}.
Let us now consider an arbitrary sequence $\left\{\tilde\bx_i:i=1,2,\ldots\right\}$, and let $\tilde\bbeta=\left(\tilde\beta_1,\tilde\beta_2,\ldots\right)^\prime$,
where, for $i=1,2,\ldots$, $\tilde\beta_i=\tilde\eta(\tilde\bx_i)$. We then define $\beta_i=\gamma_i\tilde\beta_i$, 
where for $i=1,2,\ldots$, $\gamma_i$ are independent (but non-identical) random variables, such that $0<|\gamma_i|<L<\infty$ for $i\geq 1$, and
\begin{equation}
\sum_{i=1}^{\infty}|\gamma_i|<\infty,~\mbox{almost surely.}
\label{eq:gamma_series}
\end{equation}

Also, let $\rho\in\mathbb R$ and $\sigma\in (0,\infty)=\mathbb R^+$. 
Thus, $\btheta=(\tilde\eta,\bgamma,\rho,\sigma)$, where $\bgamma=(\gamma_1,\gamma_2,\ldots)^\prime$, 
and $\bTheta=\mathcal C'(\mathcal X)\times\mathbb R^{\infty}\times\mathbb R\times\mathbb R^+$, is the parameter space. 
For our asymptotic theories regarding the multiple testing methods
that we consider, we must verify the assumptions of Shalizi for the modeling setups (\ref{eq:true_ar1_inf}) and (\ref{eq:modeled_ar1_inf}), with
this parameter space.

With respect to the above ultra high-dimensional setup, we consider the following multiple-testing framework:
\begin{align}
&H_{01}:|\rho|<1\text{ versus }H_{11}:|\rho|\geq 1\text{ and}\nonumber\\
&H_{0,i+2}:\beta_i\in\mathcal N_0\text{ versus }H_{1,i+2}:\beta_i\in\mathcal N^c_0,~\text{ for }~i=0,\ldots,m,
\label{eq:ar1_test_inf}
\end{align}
where $\mathcal N_0$ is some neighborhood of zero and $\mathcal N^c_0$ is the complement of the neighborhood in the relevant parameter space.

Verification of consistency of our non-marginal procedure amounts to verification of assumptions \ref{shalizi1}--\ref{shalizi7} of Shalizi for the above setup. 
In this regard, we make the following assumptions: 
\begin{itemize}
       \item[(B1)] $\underset{t\geq 1}{\sup}~\|z_t\|<\infty$, where, for $t\geq 1$, $\|z_t\|=\underset{i\geq 1}{\sup}~|z_{it}|$.
       \item[(B2)] For $k>1$, let $\tilde\lambda_{nk}$ be the largest eigenvalue of $\frac{\sum_{t=1}^n\bz_{m,t+k}\bz^\prime_{mt}}{n}$. 
	       We assume that $\tilde\lambda_{nk}\rightarrow 0$, as $n\rightarrow\infty$, for $k>1$.
       \item[(B3)] Let $\lambda_n$ be the largest eigenvalue of $\frac{\sum_{t=1}^n\bz_{mt}\bz^\prime_{mt}}{n}$. 
	       We assume that $\underset{n\geq 1}{\sup}~\lambda_n\leq K<\infty$.
	\item[(B4)]
	\begin{align}
	&\frac{1}{n}\sum_{t=1}^n\bbeta^\prime_m\bz_{mt} \rightarrow 0~\mbox{almost surely};~~\frac{1}{n}\sum_{t=1}^n\bbeta^\prime_{m0}\bz_{mt} \rightarrow 0;\label{eq:ass1}\\
	&\frac{1}{n}\sum_{t=1}^n\bbeta^\prime_m\bz_{mt}\bz'_{mt}\bbeta_m\rightarrow c(\bbeta)~\mbox{almost surely};
	~~\frac{1}{n}\sum_{t=1}^n\bbeta^\prime_{m0}\bz_{mt}\bz'_{mt}\bbeta_{m0}\rightarrow c(\bbeta_0)\label{eq:ass3},\\
	&\frac{1}{n}\sum_{t=1}^n\bbeta^\prime_m\bz_{mt}\bz'_{mt}\bbeta_{m0}\rightarrow c_{10}(\bbeta,\bbeta_0)~\mbox{almost surely},\label{eq:ass4}
	\end{align}
	as $n\rightarrow\infty$. In the above, $c(\bbeta_0)~(>0)$ is a finite constant; $c(\bbeta)~(>0)$ and $c_{10}(\bbeta,\bbeta_0)$ are finite quantities that depend upon 
	the choice of the sequence
	$\left\{\bbeta_{m};n=1,2,\ldots\right\}$. 

       \item[(B5)] The limits of the quantities $\bz^\prime_t\bbeta$ for almost all $\bbeta$, $\bz^\prime_t\bbeta_0$ and
       $\hat\varrho_t=\sum_{k=1}^t\rho^{t-k}_0\bz^\prime_k\bbeta_0$ exist as $t\rightarrow\infty$.
       \item[(B6)] 
      There exist positive constants $\alpha$, $c_{\rho}$, $c_{\sigma}$, $c_{\tilde\eta}$, $c_{\tilde\eta'}$ and $c_{\gamma}$ such that the following hold for sufficiently large $n$:
\begin{align}
	\pi\left(|\rho|>\exp(\left(\alpha n\right)^{1/16})\right) &\leq c_{\rho}\exp\left(-\alpha n\right);\notag\\
	\pi\left(\exp(-\left(\alpha n\right)^{1/16})\leq\sigma\leq\exp(\left(\alpha n\right)^{1/16})\right)&\geq 1-c_{\sigma}\exp\left(-\alpha n\right);\notag\\
	\pi\left(\|\tilde\eta\|\geq\exp(\left(\alpha n\right)^{1/16})\right)&\leq c_{\tilde\eta}\exp\left(-\alpha n\right);\notag\\
	\pi\left(\|\tilde\eta'\|\geq\exp(\left(\alpha n\right)^{1/16})\right)&\leq c_{\tilde\eta'}\exp\left(-\alpha n\right);\notag\\
	\pi\left(\sum_{i=1}^{\infty}|\gamma_i|\geq\exp(\left(\alpha n\right)^{1/16})\right)&\leq c_{\gamma}\exp\left(-\alpha n\right),\notag
\end{align}

       \item[(B7)] $L(m_{n+1}-m_n)\leq\exp(\left(\alpha (n+1)\right)^{1/16})-\exp(\left(\alpha n\right)^{1/16})$, for $n\geq n_0$, for some $n_0\geq 1$.
\end{itemize}

\subsection{Discussion of the assumptions in the light of the ultra high-dimensional setup}
\label{subsec:remarks_assumptions}

Condition (B1) holds if the covariates 
$\left\{z_{it}; i\geq 1, t\geq 1\right\}$, is a realization of some stochastic process with almost surely finite sup-norm, for example, Gaussian process. 
Assumption (B1), along with (\ref{eq:tilde_eta}) and (\ref{eq:gamma_series}) leads to the following result:
\begin{equation}
|\bz^\prime_{mt}\bbeta_{m0}|<C,
\label{eq:B2}
\end{equation}
for some $C>0$. To see this, first let $\bbeta_0$ correspond to the true quantities $\bgamma_0$ and $\tilde\eta_0$. Then 
observe that $|\bz^\prime_{mt}\bbeta_{m0}|\leq\sum_{i=1}^m|z_{it}||\beta_{i0}|\leq\underset{t\geq 1}{\sup}~\|z_t\|\|\tilde\eta_0\|\sum_{i=1}^{\infty}|\gamma_{i0}|<C$,
since $\underset{t\geq 1}{\sup}~\|z_t\|<\infty$ by (B5), $\|\tilde\eta_0\|<\infty$ by (\ref{eq:tilde_eta}) and $\sum_{i=1}^{\infty}|\gamma_{i0}|<\infty$ 
by (\ref{eq:gamma_series}).
Condition (B1) is required for some limit calculations and boundedness of some norms associated with concentration inequalities. 

Condition (B2) 
says that the covariates at different time points, after scaling by $\sqrt{n}$, are asymptotically orthogonal.
This condition also imply the following:
\begin{align}
&\frac{1}{n}\sum_{t=1}^n\bbeta^\prime_m\bz_{m,t+k}\bz^\prime_{mt}\bbeta_m \rightarrow 0~\mbox{almost surely, and}~~
	\frac{1}{n}\sum_{t=1}^n\bbeta^\prime_{m0}\bz_{m,t+k}\bz^\prime_{mt}\bbeta_{m0} \rightarrow 0~\mbox{for any}~k> 1;\label{eq:ass2_inf}
\end{align}
To see (\ref{eq:ass2_inf}), observe that
\begin{equation}
	\frac{1}{n}\sum_{t=1}^n\bbeta^\prime_m\bz_{m,t+k}\bz'_{mt}\bbeta_m=\bbeta^\prime_m\left(\frac{\sum_{t=1}^n\bz_{m,t+k}\bz'_{mt}}{n}\right)\bbeta_m
\leq \|\bbeta_m\|^2\left\|\left(\frac{\sum_{t=1}^n\bz_{m,t+k}\bz'_{mt}}{n}\right)\right\|_{op}.
\label{eq:B1_1}
\end{equation}
In (\ref{eq:B1_1}), $\|\bbeta_m\|$ denotes the Euclidean norm of $\bbeta_m$ and for any matrix $\bA$, $\|\bA\|_{op}$ denotes the operator norm of $\bA$ given by
$\|\bA\|_{op}=\underset{\|\bu\|=1}{\sup}~\|\bA\bu\|$. By (B2), $\left\|\left(\frac{\sum_{t=1}^n\bz_{m,t+k}\bz'_{mt}}{n}\right)\right\|_{op}\rightarrow 0$ as
$n\rightarrow\infty$. Also,
\begin{equation}
\|\bbeta_m\|^2\leq\sum_{i=1}^{\infty}\gamma^2_i\tilde\beta^2_i\leq\|\tilde\eta\|^2\sum_{i=1}^{\infty}\gamma^2_i<\infty,~\mbox{almost surely},
\label{eq:B1_2}
\end{equation}
by (\ref{eq:tilde_eta}) and (\ref{eq:gamma_series}). It follows from (\ref{eq:B1_2}) that (\ref{eq:B1_1}) is almost surely finite. This and (B2) together imply
the first part of the limit \ref{eq:ass2_inf}). Since $\|\bbeta_0\|<\infty$, the second limit of \ref{eq:ass2_inf}) follows in the same way.

As shown in Section \ref{subsec:non_ultra}, $\lambda_n\rightarrow 0$ as $n\rightarrow\infty$, even if 
$\underset{t=1,\ldots,n}{\sup}~\|\bz_{mt}\|=O(n^r),~\mbox{where}~r<1$, that is, even if (B1) does not hold. 
Since we assume only as much as $\lambda_n$ is bounded above, (B3) is a reasonably mild assumption.

In (B4), (\ref{eq:ass1}) can be made to hold in practice by centering the covariates, that is, by setting $\tilde \bz_{mt}=\bz_{mt}-\bar \bz_m$, where $\bar \bz_m=\frac{1}{n}\sum_{t=1}^n\bz_{mt}$.
In (B1) (\ref{eq:ass3}) we assume that $c(\bbeta)$ and $c_{10}(\bbeta,\bbeta_0)$ remain finite for any choice of $\left\{\bbeta_{m};n=1,2,\ldots\right\}$. To see that finiteness holds, first
note that
\begin{equation}
\frac{1}{n}\sum_{t=1}^n\bbeta^\prime_m\bz_{mt}\bz'_{mt}\bbeta_m=\bbeta^\prime_m\left(\sum_{t=1}^n\bz_{mt}\bz'_{mt}\right)\bbeta_m
\leq \|\bbeta_m\|^2\left\|\left(\frac{\sum_{t=1}^n\bz_{mt}\bz'_{mt}}{n}\right)\right\|_{op}.
\label{eq:B4_1}
\end{equation}
In (\ref{eq:B4_1}), $\|\bbeta_m\|<\infty$ almost surely, by (\ref{eq:B1_2}), and $\left\|\left(\frac{\sum_{t=1}^n\bz_{mt}\bz'_{mt}}{n}\right)\right\|_{op}<\infty$ by (B3).
Hence, (\ref{eq:B1_1}) is finite.
Similarly, $\frac{1}{n}\sum_{t=1}^n\bbeta^\prime_m\bz_{mt}\bz'_{mt}\bbeta_{m0}=\bbeta^\prime_m\left(\frac{\sum_{t=1}^n\bz_{mt}\bz'_{mt}}{n}\right)\bbeta_{m0}
\leq \|\bbeta_m\|\|\bbeta_{m0}\|\left\|\left(\frac{\sum_{t=1}^n\bz_{mt}\bz'_{mt}}{n}\right)\right\|_{op}$, which is again almost surely finite due 
to (\ref{eq:tilde_eta}), (\ref{eq:gamma_series}) and (B3).
Thus, (\ref{eq:tilde_eta}) and (\ref{eq:gamma_series}) are precisely the conditions that induce sparsity within our model in the sense of controlling 
the quadratic forms involving $\bbeta_m$ and $\bbeta_{m0}$, given that (B4) holds.
Assumptions on the existence of the limits are required for conditions (S2) and (S3) of Shalizi. 
As can be observe from Section \ref{subsec:non_ultra}, 
$\frac{1}{n}\sum_{t=1}^n\bbeta^\prime_m\bz_{mt}\bz'_{mt}\bbeta_m\rightarrow 0$, almost surely as $n\rightarrow\infty$, 
if the asymptotically orthogonal covariates satisfy 
$\underset{t=1,\ldots,n}{\sup}~\|\bz_{mt}\|=O(n^r),~\mbox{where}~r<1$, that is, even if (B1) does not hold.
Hence, in this situation, the required limits of the quadratic forms exist and are zero, under very mild conditions.

Again, the limit existence assumption (B5) is required for verification of conditions (S2) and (S3) of Shalizi. 

Assumption (B6), required to satisfy condition (S5) of Shalizi, is reasonably mild. The threshold $\exp(\left(\alpha n\right)^{1/16})$ for the probabilities involving 
$\|\tilde\eta\|$ and $\|\tilde\eta'\|$
can be replaced with the order of $\sqrt{n}$ for Gaussian process priors or for independent sub-Gaussian components of $\bbeta$. 
However, note that priors such as gamma or inverse gamma for $\sigma$ do not necessarily satisfy the condition. In such cases, one can modify the prior by
replacing the tail part of the prior, after an arbitrarily large positive value, with a thin-tailed prior, such as normal. In practice, such modified priors would be
effectively the same as gamma or inverse gamma priors, and yet would satisfy the conditions of (B6).


Assumption (B7), in conjunction with boundedness of $|\gamma_i|$, for all $i$ by $L$, is a mild condition ensuring that $\mathcal G_n$ are increasing in $n$, when $n\geq n_0$, for some $n_0\geq 1$.

\subsection{High-dimensional but not ultra high-dimensional setup}
\label{subsec:non_ultra}

The setup we discussed so far deals with the so-called ultra high-dimensional problem, in the sense that $\frac{m_n}{n}\rightarrow\infty$ as $n\rightarrow\infty$. This is a challenging
problem to address and we required a prior for $\bbeta$ satisfying $\|\bbeta\|<\infty$ almost surely.
However, if we are only interested in the problem where $\frac{m_n}{n}\rightarrow 0$ as $n\rightarrow\infty$, then it is not necessary to insist on priors to ensure finiteness of $\|\bbeta\|$. 
For example, if the covariates $\bz_{mt}$ are orthogonal, 
then assuming that 
\begin{equation}
\underset{t=1,\ldots,n}{\sup}~\|\bz_{mt}\|=O(n^r),~\mbox{where}~r<1, 
\label{eq:orth}
\end{equation}
$\frac{1}{n}\sum_{t=1}^n\bz_{mt}\bz^\prime_{mt}$ has maximum eigenvalue $O(n^{r-1})$, so that (\ref{eq:B1_1})
entails 
\begin{equation}
\frac{1}{n}\sum_{t=1}^n\bbeta^\prime_m\bz_{mt}\bz'_{mt}\bbeta_m =O\left(\|\bbeta_m\|^2n^{r-1}\right).
\label{eq:B1_3}
\end{equation}
Now, if the components of $\bbeta_m$ are independent and sub-Gaussian with mean zero, then by the Hanson-Wright inequality (see, for example, \ctn{Rudelson13}) we have 
\begin{align}
&P\left(\left|\sum_{t=1}^m\beta^2_t-\sum_{t=1}^mE(\beta^2_t)\right|>n^{1-r}-\sum_{t=1}^mE(\beta^2_t)\right)\notag\\
&\qquad\leq 2\exp\left(-L_1\min\left\{\frac{\left(n^{1-r}-\sum_{t=1}^mE(\beta^2_t)\right)^2}{L^4_2m},\frac{n^{1-r}-\sum_{t=1}^mE(\beta^2_t)}{L^2_2}\right\}\right),
\label{eq:hw_finite_sum}
\end{align}
where $L_1>0$ is some constant and $L_2$ is the upper bound of the sub-Gaussian norm. Let $\tilde m=\sum_{t=1}^mE(\beta^2_t)$. 
If $\frac{n^{1-r}-\tilde m}{\sqrt{\tilde m}}\rightarrow \tilde c~(>0)$, where $\tilde c$ is finite or infinite, then (\ref{eq:hw_finite_sum}) is summable. Hence, by the Borel-Cantelli lemma,
$\sum_{t=1}^m\beta^2_t\leq n^{1-r}$ almost surely, as $n\rightarrow\infty$. It then follows from (\ref{eq:B1_3}) that $\frac{1}{n}\sum_{t=1}^n\bbeta^\prime_m\bz_{mt}\bz'_{mt}\bbeta_m<\infty$
almost surely as $n\rightarrow\infty$. 

For the non-ultra high-dimensional setup, the problem is largely simplified. Indeed, introduction of $\tilde\eta$ and $\tilde\eta^\prime$ are not required, as we can directly consider
sub-Gaussian priors for $\bbeta$ as detailed above. Consequently, in (B3), only the first two inequalities are needed and assumption (B6) is no longer required.
Since the ultra high-dimensional setup is far more challenging than the non-ultra high-dimensional setup, we consider only the former setup for our purpose, and note that the latter
setup can be dealt with using almost the same ideas but with much less effort.

Assumptions (B1)--(B6) lead to the following results that are the main ingredients in proving our posterior convergence in the ultra high-dimensional setup.
\begin{lemma}
\label{lemma:lemma1_inf}
Under (B1), (B2) and (B5), the KL-divergence rate $h(\btheta)$ exists for each $\btheta\in\bTheta$ and is given by
\begin{multline}
h(\btheta)
=\log\left(\frac{\sigma}{\sigma_0}\right)+\left(\frac{1}{2\sigma^2}-\frac{1}{2\sigma^2_0}\right)\left(\frac{\sigma^2_0}{1-\rho^2_0}
+\frac{c(\bbeta_0)}{1-\rho^2_0}\right)\\
+\left(\frac{\rho^2}{2\sigma^2}-\frac{\rho^2_0}{2\sigma^2_0}\right)\left(\frac{\sigma^2_0}{1-\rho^2_0}
+\frac{c(\bbeta_0)}{1-\rho^2_0}\right)
+\frac{c(\bbeta)}{2\sigma^2}-\frac{c(\bbeta_0)}{2\sigma^2_0}\\
-\left(\frac{\rho}{\sigma^2}-\frac{\rho_0}{\sigma^2_0}\right)\left(\frac{\rho_0\sigma^2_0}{1-\rho^2_0}+\frac{\rho_0c(\bbeta_0)}{1-\rho^2_0}\right)
-\left(\frac{c_{10}(\bbeta,\bbeta_0)}{\sigma^2}-\frac{c(\bbeta_0)}{\sigma^2_0}\right).
\label{eq:htheta_inf}
\end{multline}
\end{lemma}

\begin{theorem}
\label{theorem:theorem2_inf}
Under (B1), (B2) and (B5), the asymptotic equipartition property holds and is given by
\begin{equation*}
\underset{n\rightarrow\infty}{\lim}~\frac{1}{n}\log R_n(\btheta)=-h(\btheta).
\end{equation*}
Furthermore, the convergence is uniform on any compact subset of $\bTheta$.
\end{theorem}
Lemma \ref{lemma:lemma1_inf} and Theorem \ref{theorem:theorem2_inf} ensure that \ref{shalizi1} -- \ref{s3} hold, and \ref{s4} holds since $h(\btheta)$ is almost surely finite. 
(B6) implies that $\mathcal G_n$ increases to $\bTheta$. In Section \ref{subsec:A5} we verify \ref{s5}.

Now observe that the aim of assumption (S6) is to ensure that (see the proof of Lemma 7 of \ctn{Shalizi09}) 
for every $\varepsilon>0$ and for all $n$ sufficiently large,
\begin{equation*}
\frac{1}{n}\log\int_{\mathcal G_n}R_n(\btheta)d\pi(\btheta)\leq -h\left(\mathcal G_n\right)+\varepsilon,~\mbox{almost surely}.
\end{equation*}
Since $h\left(\mathcal G_n\right)\rightarrow h\left(\bTheta\right)$ as $n\rightarrow\infty$, it is enough to verify that for every $\varepsilon>0$ and for all $n$ sufficiently large,
\begin{equation}
\frac{1}{n}\log\int_{\mathcal G_n}R_n(\btheta)d\pi(\btheta)\leq -h\left(\bTheta\right)+\varepsilon,~\mbox{almost surely}.
\label{eq:s6_1}
\end{equation}
In this regard, first observe that 
\begin{align}
\frac{1}{n}\log\int_{\mathcal G_n}R_n(\btheta)d\pi(\btheta)&\leq\frac{1}{n}\log\left[\underset{\btheta\in\mathcal G_n}{\sup}~R_n(\btheta)\pi(\mathcal G_n)\right]\notag\\
&=\frac{1}{n}\log\left[\underset{\btheta\in\mathcal G_n}{\sup}~R_n(\btheta)\right]+\frac{1}{n}\log\pi(\mathcal G_n)\notag\\
&=\underset{\btheta\in\mathcal G_n}{\sup}~\frac{1}{n}\log R_n(\btheta)+\frac{1}{n}\log\pi(\mathcal G_n)\notag\\
&\leq \frac{1}{n}\underset{\btheta\in\mathcal G_n}{\sup}~\log R_n(\btheta),
\label{eq:s6_2}
\end{align}
where the last inequality holds since $\frac{1}{n}\log\pi(\mathcal G_n)\leq 0$. Now, letting $\mathcal S=\left\{\btheta:h(\btheta)\leq\kappa\right\}$, where
$\kappa>h\left(\bTheta\right)$ is large as desired,
\begin{align}
\underset{\btheta\in\mathcal G_n}{\sup}~\frac{1}{n}\log R_n(\btheta)
&\leq\underset{\btheta\in\bTheta}{\sup}~\frac{1}{n}\log R_n(\btheta)
=\underset{\btheta\in\mathcal S\cup\mathcal S^c}{\sup}~\frac{1}{n}\log R_n(\btheta)
\notag\\
&\leq \max\left\{\underset{\btheta\in\mathcal S}{\sup}~\frac{1}{n}\log R_n(\btheta),\underset{\btheta\in\mathcal S^c}{\sup}~\frac{1}{n}\log R_n(\btheta)\right\}.
\label{eq:s6_3}
\end{align}
From (\ref{eq:htheta_inf}) it is clear that $h(\btheta)$ is continuous in $\btheta$ and that $h(\btheta)\rightarrow\infty$ as $\|\btheta\|\rightarrow\infty$. In other words,
$h(\btheta)$ is a continuous coercive function. Hence, $\mathcal S$ is a compact set (see, for example, \ctn{Lange10}). 
Hence it easily follows that (see \ctn{Chatterjee18b}), that
\begin{equation}
\underset{\btheta\in\mathcal S}{\sup}~\frac{1}{n}\log R_n(\btheta)\rightarrow \underset{\btheta\in\mathcal S}{\sup}~-h(\btheta)=-h\left(\mathcal S\right),~\mbox{almost surely, as}~n\rightarrow\infty.
\label{eq:s6_4}
\end{equation}

We now show that 
\begin{equation}
\underset{\btheta\in\mathcal S^c}{\sup}~\frac{1}{n}\log R_n(\btheta)\leq -h\left(\bTheta\right)~\mbox{almost surely, as}~ n\rightarrow\infty.
\label{eq:s6_5}
\end{equation}
First note that
if $\underset{\btheta\in\mathcal S^c}{\sup}~\frac{1}{n}\log R_n(\btheta)> -h\left(\bTheta\right)$ infinitely often, then $\frac{1}{n}\log R_n(\btheta)> -h\left(\bTheta\right)$
for some $\btheta\in\mathcal S^c$ infinitely often. But $\frac{1}{n}\log R_n(\btheta)> -h\left(\bTheta\right)$ if and only if 
$
\frac{1}{n}\log R_n(\btheta)+h(\btheta)> h(\btheta)-h\left(\bTheta\right),~\mbox{for}~\btheta\in\mathcal S^c. 
$
Hence, if we can show that  
\begin{equation}
P\left(\left|\frac{1}{n}\log R_n(\btheta)+h(\btheta)\right|> \kappa-h\left(\bTheta\right),~\mbox{for}~\btheta\in\mathcal S^c~\mbox{infinitely often}\right)=0, 
\label{eq:s6_6}
\end{equation}
then (\ref{eq:s6_5}) will be proved. We use the Borel-Cantelli lemma to prove (\ref{eq:s6_6}). In other words, we prove that
\begin{theorem}
\label{theorem:theorem3}
Under (B5), (\ref{eq:tilde_eta}) and (\ref{eq:gamma_series}), 
\begin{equation}
\sum_{n=1}^{\infty}\int_{\mathcal S^c}P\left(\left|\frac{1}{n}\log R_n(\btheta)+h(\btheta)\right|> \kappa-h\left(\bTheta\right)\right)d\pi(\btheta)<\infty. 
\label{eq:s6_7}
\end{equation}
\end{theorem}
The proof of Theorem \ref{theorem:theorem3} heavily uses (\ref{eq:B2}), which is ensured by (B5), (\ref{eq:tilde_eta}) and (\ref{eq:gamma_series}).
Since $h(\btheta)$ is continuous, (S7) holds trivially. 

We provide detailed verification of the seven assumptions of Shalizi in the supplement, which leads to the following result:
\begin{theorem}
	\label{theorem:ar1_inf}
	Under assumptions (B1) -- (B6), the non-marginal multiple testing procedure for testing (\ref{eq:ar1_test_inf}) is consistent. 
\end{theorem}
Needless to mention, all the results on error convergence of the non-marginal method also continue to hold for this setup under (B1) -- (B6), thanks to verification of
Shalizi's conditions.

\subsection{Remark on identifiability of our model and posterior consistency}
\label{subsec:identifiability}
Note that we have modeled $\bbeta$ in terms of $\bgamma$ and $\tilde\eta$. But from the likelihood it is evident that although $\bbeta$ is identifiable, $\bgamma$ and $\tilde\eta$ are not.
But this is not an issue since our interest is in the posterior of $\bbeta$, not of $\bgamma$ or $\tilde\eta$. Indeed, Theorem 3 of Shalizi guarantees that the posterior
of the set $\{\btheta:h(\btheta)\leq h(\bTheta)+\varepsilon\}$ tends to 1 as $n\rightarrow\infty$, for any $\varepsilon>0$. We show in the supplement that $h(\bTheta)=0$ in our case.
Since $h(\btheta_0)=0$, where $\btheta_0$ is the true parameter which includes $\bbeta_0$ and lies in $\{\btheta:h(\btheta)<\varepsilon\}$ for any $\varepsilon>0$, it follows that 
the posterior of $\bbeta$ is consistent.


\section{Summary and conclusion}
\label{sec:conclusion_inf}
In this article, we have investigated asymptotic properties of the Bayesian non-marginal procedure under the general dependence structure when the number of hypotheses also tend
to infinity with the sample size. We specifically showed that our method is consistent even in this setup, and that the different Bayesian versions of the error rates converge to zero 
exponentially fast, and that the expectations of the Bayesian versions with respect to the data also tend to zero. Since our results hold for any choice of the groups, 
it follows that they hold even for singleton groups, that is, 
for marginal decision rules. The 
results associated with $\alpha$-control also continue to hold in the same spirit as the finite-dimensional setup developed in \ctn{chandra2020}. 
Interestingly, provided that Shalizi's conditions hold, almost no assumption is required on the growth rate of the number of hypotheses to establish the results
of the multiple testing procedures in high dimensions. Although in several cases, unlike the exact fixed-dimensional limits established in \ctn{chandra2020}, 
the exact high-dimensional limits associated with the error rates could not be established, 
exponential convergence to zero in high dimensions could still be achieved. Moreover, internal consistency of our results, as we make transition from fixed dimension
to high dimensions, are always ensured.

An important objective of this research is to show that the finite-dimensional time-varying variable selection problem in the autoregressive setup introduced in 
\ctn{chandra2020} admits extension to the setup where the number of covariates to be selected by our Bayesian non-marginal procedure, grows with sample size. Indeed, we have shown that under reasonable
assumptions, our asymptotic theories remain valid for this problem for both high-dimensional and ultra high-dimensional situations. Different priors for the regression coefficients  
are of course warranted, and we have discussed the classes of such relevant priors for the two different setups. As much as we are aware of, at least in the time series context,
such high-dimensional multiple hypotheses testing is not hitherto dealt with. The priors that we introduce, particularly in the ultra high-dimensional context, 
also do not seem to have been considered before. These priors, in conjunction with the equipartition property, help control sparsity of the model quite precisely. As such,
these ideas seem to be of independent interest for general high-dimensional asymptotics.

\renewcommand\thefigure{S-\arabic{figure}}
\renewcommand\thetable{S-\arabic{table}}
\renewcommand\thesection{S-\arabic{section}}
\renewcommand\thetheorem{S-\arabic{theorem}}


\begin{center}
{\bf \LARGE Supplementary Material}
\end{center}

\section{Proof of Theorem \ref{theorem:inf_dim_cons}}
\begin{proof}
	%
	From conditions 
	(\ref{eq:liminf_beta}) and (\ref{eq:limsup_beta}), it follows that there exists $n_1$ such that for all $n>n_1$
	\begin{align}
	\beta_n&>\underline{\beta}-\delta,\\
	\beta_n&<1-\delta,\text{ such that}
	\end{align}
	$\underline{\beta}-\delta>0$ and $1-\bar{\beta}>\delta$, for some $\delta>0$.
	It follows using this, (\ref{eq:ineq_start}) and (\ref{eq:ineq2}), that for $n>n_1$,
	\begin{align}
		&\sum_{i:\bd(m_n)\in \mathbb{D}_{i,m_n}^c }^{m_n} d^t_iw_{in}(\bd^t(m_n))- 
		\sum_{i:\bd(m_n)\in \mathbb{D}_{i,m_n}^c}^{m_n}d_iw_{in}(\bd(m_n))\\ 
		&\qquad\qquad> \left(1- e^{-n(J-\epsilon )}\right) 
		\sum_{i:\bd(m_n)\in \mathbb{D}_{i,m_n}^c } d^t_i-e^{-n(J-\epsilon )}\sum_{i:\bd(m_n)\in \mathbb{D}_{i,m_n}^c } d_i,~\mbox{and}\notag\\
		&\beta_n \left( \sum_{i:\bd\in \mathbb{D}_{i,m_n}^c}^{m_n} d^t_i-\sum_{i:\bd(m_n)\in \mathbb{D}_{i,m_n}^c }^{m_n} d_i\right)< (1-\delta )
		\sum_{i:\bd\in \mathbb{D}_{i,m_n}^c}^{m_n} d^t_i - (\underline\beta-\delta) \sum_{i:\bd(m_n)\in \mathbb{D}_{i,m_n}^c}^{m_n} d_i.
	\end{align}
	Now $n_1$ can be appropriately chosen such that $e^{-n(J-\epsilon )}<\min \{ \delta, \underline\beta-\delta\}$. 
	Hence, for $n>\max \{ n_0, n_1\}$,
	\begin{align*}
		&\sum_{i:\bd\in \mathbb{D}_{i,m_n}^c }^{m_n} d^t_iw_{in}(\bd^t(m_n))- \sum_{i:\bd(m_n)\in \mathbb{D}_{i,m_n}^c}^{m_n}d_iw_{in}(\bd(m_n))
		> \beta_n \left( \sum_{i:\bd(m_n)\in \mathbb{D}_{i,m_n}^c}^{m_n}d^t_i-\sum_{i:\bd(m_n)\in \mathbb{D}_{i,m_n}^c }^{m_n} d_i\right) ,\notag\\
		&\qquad~\mbox{for all}~\bd(m_n)\neq\bd^t(m_n),\text{ almost surely};\\
		\Rightarrow &\sum_{i=1}^{m_n}d^t_i(w_{in}(\bd^t(m_n))-\beta_n)>\sum_{i=1}^{m_n}d_i(w_{in}(\bd(m_n))-\beta_n) ,
		~\mbox{for all}~\bd(m_n)\neq\bd^t(m_n),\text{ almost surely};\\	
	\Rightarrow & \lim_{n\rightarrow\infty}\delta_{\mathcal {NM}}(\bd^t(m_n)|\bX_n)=1,~\mbox{almost surely}. 
	\end{align*}
	Hence, (\ref{eq:inf_dim_cons1}) holds, and by the dominated convergence theorem, (\ref{eq:inf_dim_cons2}) also follows.
\end{proof}

\section{Proof of Theorem \ref{theorem:mpBFDR_bound_inf}}
\begin{proof}
	
	\begin{align*}
		&\sum_{\bd(m_n)\neq\bzero}\frac{\sum_{i=1}^{m_n}d_i (1- w_{in}(\bd(m_n)))}{\sum_{i=1}^{m_n}d_i}
		\dnm\left(\bd(m_n)|\bX_n\right)\\
		=&\frac{\sum_{i=1}^{m_n}d_i^t (1- w_{in}(\bd^t(m_n)))}{\sum_{i=1}^{m_n}d_i^t}
		\delta_{\mathcal {NM}}\left(\bd^t(m_n)|\bX_n\right)+\sum_{\bd(m_n)\neq\bd^t(m_n)\neq\bzero}\frac{\sum_{i=1}^{m_n}d_i (1- w_{in}(\bd(m_n)))}{\sum_{i=1}^{m_n}d_i}
		\dnm\left(\bd(m_n)|\bX_n\right).
	\end{align*}
	Following Theorem \ref{theorem:inf_dim_cons}, it holds, almost surely, that there exists $N\geq 1$ such that for all $n>N$, 
	$\delta_{\mathcal {NM}}\left(\bd(m_n)|\bX_n\right)=0$ for all $\bd(m_n)\neq\bd^t(m_n)$. Therefore, for $n>N$,
	\begin{align*}
		&\sum_{\bd(m_n)\neq\bzero}\frac{\sum_{i=1}^{m_n}d_i (1- w_{in}(\bd(m_n)))}{\sum_{i=1}^{m_n}d_i}
		\delta_{\mathcal {NM}}\left(\bd(m_n)|\bX_n\right)\\
		= &\frac{\sum_{i=1}^{m_n}d_i^t (1- w_{in}(\bd^t(m_n)))}{\sum_{i=1}^{m_n}d_i^t}
		\delta_{\mathcal {NM}}\left(\bd^t(m_n)|\bX_n\right)\\
		\leq& \frac{\sum_{i=1}^{m_n}d_i^t e^ {-n(J-\epsilon ) }}{\sum_{i=1}^{m_n}d_i^t}\\
	=& e^ {-n(J-\epsilon ) }.
	\end{align*}
	Thus, (\ref{eq:inf_dim_mpBFDR}) is established. Using (\ref{eq:ineq_end}) and Corollary \ref{corr:corr1}, (\ref{eq:inf_dim_pBFDR}) follows in the same way.
\end{proof}

\subsection{Proof of Theorem \ref{theorem:lim_BFDR}}
\begin{proof}
	Note that
	\begin{align*}
		&mpBFDR\\
		&= E_{\bX_n} \left[\sum_{\bd(m_n)\in\mathbb{D}_{m_n}} \frac{\sum_{i=1}^{m_n}d_i(1- w_i(\bd(m_n)) )}{\sum_{i=1}^{m_n}d_i}
		\delta_\beta(\bd(m_n)|\bX_n)\bigg{|}\dnm(\bd(m_n)=\bzero|\bX_n)=0 \right]\\
		=& E_{\bX_n} \left[\sum_{\bd(m_n)\in\mathbb{D}_{m_n}} \frac{\sum_{i=1}^{m_n}d_i(1- w_i(\bd(m_n)) )}{\sum_{i=1}^{m_n}d_i}
		\dnm(\bd(m_n)|\bX_n)\bigg{|}\dnm(\bd(m_n)=\bzero|\bX_n)=0 \right]\\
		=& E_{\bX_n} \left[\sum_{\bd(m_n)\in\mathbb{D}_{m_n}} \frac{\sum_{i=1}^{m_n}d_i(1- w_i(\bd(m_n)) )}{\sum_{i=1}^{m_n}d_i} 
		I\left( \sum_{i=1}^{m_n}d_i>0 \right)\dnm(\bd(m_n)|\bX_n) \right] 
		\frac{1}{P_{\bX_n}\left[ \dnm(\bd(m_n)=\bzero|\bX_n)=0\right] } \\
		=&E_{\bX_n}\left[\sum_{\bd(m_n)\in\mathbb{D}_{m_n}\setminus\left\lbrace \bzero\right\rbrace} \frac{\sum_{i=1}^{m_n}d_i(1- w_i(\bd(m_n)) )}{\sum_{i=1}^{m_n}d_i} 
		\dnm(\bd(m_n)|\bX_n) \right] \frac{1}{P_{\bX_n}\left[ \dnm(\bd(m_n)=\bzero|\bX_n)=0\right] }.
	\end{align*}	
	From Theorem \ref{theorem:mpBFDR_bound_inf}, $\mfdr\rightarrow0$, as $n\rightarrow\infty$. Also we have 
	\begin{equation*}
		0\leq\sum_{\bd(m_n)\in\mathbb{D}_{m_n}\setminus\left\lbrace \bzero\right\rbrace} \frac{\sum_{i=1}^{m_n}d_i(1- w_i(\bd(m_n)) )}{\sum_{i=1}^{m_n}d_i} 
		\dnm(\bd(m_n)|\bX_n)\leq \mfdr\leq1.
	\end{equation*}
	Therefore by the dominated convergence theorem, $E_{\bX_n}\left[\sum_{\bd(m_n)\in\mathbb{D}_{m_n}\setminus\left\lbrace \bzero\right\rbrace} 
	\frac{\sum_{i=1}^{m}d_i(1- w_i(\bd(m_n)) )}{\sum_{i=1}^{m_n}d_i} 
	\dnm(\bd(m_n)|\bX_n) \right]\rightarrow0$, as $n\rightarrow\infty$. From \ref{A2} we have $\bd^t(m_n)\neq\bzero$ and from Theorem \ref{theorem:inf_dim_cons}
	we have $E_{\bX_n}[\dnm(\bd^t(m_n)|\bX_n)]\rightarrow1$. Thus $P_{\bX_n}\left[ \dnm(\bd(m_n)=\bzero|\bX_n)=0\right]\rightarrow1$, as $n\rightarrow\infty$. 
	This proves the result.
	
	It can be similarly shown that $pBFDR\rightarrow 0$, as $n\rightarrow\infty$.
\end{proof}

\section{Proof of Theorem \ref{theorem:mpBFNR_bound_inf}}
\begin{proof}
The proof follows in the same way as that of Theorem \ref{theorem:mpBFDR_bound_inf}, using (\ref{A2}) in addition.
\end{proof}

\subsection{Proof of Theorem \ref{theorem:lim_BFNR}}
\begin{proof}
The proof follows in the same way as that of Theorem \ref{theorem:lim_BFDR}, using (\ref{A2}) in addition.
\end{proof}

\subsection{Proof of Theorem \ref{theorem:mpBFDR_alpha}}
\begin{proof}
	Theorem 3.4 of \ctn{chandra2019} shows that $mpBFDR$ is non-increasing in $\beta$. Hence, for every $n>1$, 
	the maximum error that can be incurred is at $\beta=0$ where we 
	actually maximize $\sum_{i=1}^{m_n}d_iw_{in}(\bd(m_n))$.
	Let \begin{align*}
		\hat\bd(m_n) &=\argmax_{\bd(m_n)\in\mathbb{D}_{m_n}} \sum_{i=1}^{m_n}d_iw_{in}(\bd(m_n))
		= \argmax_{\bd(m_n)\in\mathbb{D}_{m_n}}\left[\sum_{i=1}^{m_{1n}}d_iw_{in}(\bd(m_n))+ \sum_{i=m_{1n}+1}^{m_n} d_iw_{in}(\bd(m_n)) \right] 
	\end{align*}	
	Since the groups in $\{G_{i,m_n}:i=1,\ldots,m_{1n}\}$ have no overlap with those in $\{G_{i,m_n}:i=m_{1n}+1,\ldots,m_n\}$, 
	$\sum_{i=1}^{m_{1n}}d_iw_{in}(\bd(m_n))$ and 
	$\sum_{i=m_{1n}+1}^{m_n}d_iw_{in}(\bd(m_n))$ can be maximized separately.

	Let us define the following notations:
	\begin{align}
		&Q_{\bd(m_n)}=\left\{i\in\{1,\ldots,m_n\}:\mbox{all elements of}~\bd_{G_{i,m_n}}~\mbox{are correct}\right\};\notag\\
		&Q_{\bd(m_n)}^{m_{1n}}=Q_{\bd(m_n)}\cap \{1,2,\ldots,m_{1n}\},~	Q_{\bd(m_n)}^{m_{1n}c}=\{1,2,\cdots,m_{1n}\}\setminus Q_{\bd(m_n)}^{m_{1n}}.\notag
	\end{align}
	
	Now,
	\begin{align*}
		&\sum_{i=1}^{m_{1n}} d_iw_{in}(\bd(m_n)) - \sum_{i=1}^{m_{1n}} d_i^tw_{in}(\bd^t(m_n))\\
		=& \left[ \sum_{i\in Q_{\bd(m_n)}^{m_{1n}}} d_iw_{in}(\bd(m_n)) - \sum_{i\in Q_{\bd(m_n)}^{m_{1n}}} d_i^tw_{in}(\bd^t(m_n)) \right] 
		+ \left[ \sum_{i\in Q_{\bd(m_n)}^{m_{1n}c}} d_iw_{in}(\bd(m_n)) 
		- \sum_{i\in Q_{\bd(m_n)}^{m_{1n}c}} d_i^tw_{in}(\bd^t(m_n)) \right]\\
		=& \sum_{i\in Q_{\bd(m_n)}^{m_{1n}c}} d_iw_{in}(\bd(m_n)) - \sum_{i\in Q_{\bd(m_n)}^{m_{1n}c}} d_i^tw_{in}(\bd^t(m_n)), 
	\end{align*}
	since for any $\bd(m_n)$, $\sum_{i\in Q_{\bd(m_n)}^{m_{1n}}} d_iw_{in}(\bd(m_n))=\sum_{i\in Q_{\bd(m_n)}^{m_{1n}}} d_i^tw_{in}(\bd^t(m_n))$ 
	by definition of $Q_{\bd(m_n)}^{m_{1n}}$.
	
	Note that $\sum_{i\in Q_{\bd(m_n)}^{m_{1n}c}} d_i^tw_{in}(\bd^t(m_n))$ can not be zero as it contradicts \ref{ass_groups} that 
	$\left\{G_{i,m_n}:i=1,\ldots,m_{1n}\right\}$ have at least one false null hypothesis. 
	
	Now, from (\ref{eq:ineq_start}) and (\ref{eq:ineq2}), we obtain for $n\geq n_0(\epsilon)$,
	\begin{align}
		&\sum_{i\in Q_{\bd(m_n)}^{m_{1n}c}} d_iw_{in}(\bd(m_n)) - \sum_{i\in Q_{\bd(m_n)}^{m_{1n}c}} d_i^tw_{in}(\bd^t(m_n))\notag\\	
		&\qquad<e^{-n(J-\epsilon)}\sum_{i\in Q_{\bd(m_n)}^{m_{1n}c}} \left(d_i+d^t_i\right)-\sum_{i\in Q_{\bd(m_n)}^{m_{1n}c}}d^t_i\notag\\
		&\qquad<2m_{1n}e^{-n(J-\epsilon)}-\sum_{i\in Q_{\bd(m_n)}^{m_{1n}c}}d^t_i.
		\label{eq:eq1}
	\end{align} 
	By our assumption (\ref{eq:infdim_limit}), $m_ne^{-n(J-\epsilon)}\rightarrow 0$ as $n\rightarrow\infty$, so that $m_{1n}e^{-n(J-\epsilon)}\rightarrow 0$
	as $n\rightarrow\infty$. Also, $\sum_{i\in Q_{\bd(m_n)}^{m_{1n}c}}d^t_i>0$. Hence, (\ref{eq:eq1}) is negative for sufficient;y large $n$.
	In other words, $\bd^t(m_n)$ maximizes $\sum_{i=1}^{m_{1n}} d_iw_{in}(\bd(m_n))$ for sufficiently large $n$.
	
	Let us now consider the term $\sum_{i=m_{1n}+1}^{m_n} d_iw_{in}(\bd(m_n))$. Note that $\sum_{i=m_{1n}+1}^{m_n} d_i^tw_{in}(\bd^t(m_n))= 0$ by \ref{ass_groups}. 
	For any finite $n$, $\sum_{i=m_{1n}+1}^{m_n} d_iw_{in}(\bd(m_n))$ is maximized for some decision configuration $\tilde\bd(m_n)$ where $\tilde d_i=1$ for at least one 
	$i\in \{m_{1n}+1,\ldots,m_n\}$. 
	In that case, $$\hat\bd^t(m_n)=(d^t_1,\ldots,d^t_{m_{1n}},\tilde d_{m_{1n}+1},\tilde d_{m_{1n}+2},\ldots,\tilde d_{m_n}),$$ so that
	for sufficiently large $n$,
	\begin{align}
		&\frac{\sum_{i=1}^{m_n} \hat d_i(1-w_{in}(\hat \bd(m_n)))}{\sum_{i=1}^{m_n}\hat d_i}
		\geq 1-\frac{\sum_{i=1}^{m_{1n}}d^t_iw_{in}(\bd^t(m_n))+(m_n-m_{1n})e^{-n(J-\epsilon)}}{\sum_{i=1}^{m_n}d^t_i+1}\notag\\
		&\qquad=\frac{1+\sum_{i=1}^{m_{1n}}d^t_i\left(1-w_{in}(\bd^t)\right)}{\sum_{i=1}^{m_n}d^t_i+1}-\frac{(m_n-m_{1n})e^{-n(J-\epsilon)}}{\sum_{i=1}^{m_n}d^t_i+1}.
		\label{eq:eq2}
	\end{align}	
       Now note that
	\begin{equation}
		0<\frac{\sum_{i=1}^{m_{1n}}d^t_i\left(1-w_{in}(\bd^t)\right)}{m_n}<e^{-n(J-\epsilon)}\frac{\sum_{i=1}^{m_{1n}}d^t_i}{m_n}
		<\frac{m_{1n}}{m_n}e^{-n(J-\epsilon)}.
		\label{eq:eq3}
	\end{equation}
	Since the right most side of (\ref{eq:eq3}) tends to zero as $n\rightarrow\infty$ due to (\ref{eq:q}), it follows that 
	$\frac{\sum_{i=1}^{m_{1n}}d^t_i\left(1-w_{in}(\bd^t)\right)}{m_n}\rightarrow 0$ as $n\rightarrow\infty$.
	Hence, dividing the numerators and denominators of the right hand side of (\ref{eq:eq2}) by $m_n$ and taking limit as $n\rightarrow\infty$
	shows that 
	\begin{equation}
		\lim_{n\rightarrow\infty}\frac{\sum_{i=1}^{m_n} \hat d_i(1-w_{in}(\hat \bd(m_n)))}{\sum_{i=1}^{m_n}\hat d_i}\geq 0.
		\label{eq:eq4}
	\end{equation}
	almost surely, for all data sequences.
	Boundedness of $\frac{\sum_{i=1}^{m_n} d_i(1-w_{in}(\bd(m_n)))}{\sum_{i=1}^{m_n}d_i}$ for all $\bd(m_n)$ and $\bX_n$ ensures uniform integrability, 
	which, in conjunction with
	the simple observation that for $\beta=0$, $$P\left(\delta_{\mathcal {NM}}(\bd(m_n)=\bzero|\bX_n)=0\right)=1$$ for all $n\geq 1$, guarantees
	that under \ref{ass_groups}, $\underset{n\rightarrow\infty}{\lim}~mpBFDR\geq 0$.
	
	Now, if $G_{m_{1n}+1},\ldots,G_{m_n}$ are all disjoint, each consisting of only one true null hypothesis, 
	then $\sum_{i=m_{1n}+1}^{m_n} d_iw_{in} (\bd(m_n))$ will be maximized by $\tilde\bd(m_n)$ where $\tilde d_i=1$ for all $i\in \{m_{1n}+1,\ldots,m_n\}$. 
	Since $d^t_i$; $i=1,\ldots,m_{1n}$ maximizes $\sum_{i=1}^{m_{1n}} d_iw_{in}(\bd(m_n))$ for large $n$, it follows that 
	$\hat\bd(m_n)=(d^t_1,\ldots,d^t_{m_{1n}},1,1,\ldots,1)$
	is the maximizer of $\sum_{i=1}^{m_n} d_iw_{in}(\bd(m_n))$ for large $n$.
	In this case,
	\begin{equation}
		\frac{\sum_{i=1}^{m_n} \hat d_i(1-w_{in}(\hat \bd(m_n)))}{\sum_{i=1}^{m_n}\hat d_i} 
		= 1-\frac{\sum_{i=1}^{m_{1n}}d^t_iw_{in}(\bd^t(m_n))+\sum_{i=m_{1n}+1}^{m_{n}}w_{in}(\bone)}{\sum_{i=1}^{m_n}d^t_i+m_n-m_{1n}}.
	\label{eq:upper_mpBFDR2}
	\end{equation}
        Now, for large enough $n$,
        \begin{equation}
	\left(1-e^{-n(J-\epsilon)}\right)\frac{\sum_{i=1}^{m_{1n}}d^t_i}{m_n}< \frac{\sum_{i=1}^{m_{1n}}d^t_iw_{in}(\bd^t(m_n))}{m_n}<\frac{\sum_{i=1}^{m_{1n}}d^t_i}{m_n}.
	\label{eq:eq5}
        \end{equation}
	Since due to (\ref{eq:p}), $\frac{\sum_{i=1}^{m_{1n}}d^t_i}{m_n}\rightarrow p$, as $n\rightarrow\infty$, it follows from (\ref{eq:eq5}) that
	\begin{equation}
		\frac{\sum_{i=1}^{m_{1n}}d^t_iw_{in}(\bd^t(m_n))}{m_n}\rightarrow p, ~\mbox{as}~n\rightarrow\infty.
		\label{eq:eq6}
	\end{equation}
	Also, since for large enough $n$,
        \begin{equation*}
		0< \frac{\sum_{i=m_{1n}+1}^{m_{n}}w_{in}(\bone)}{m_n}<\frac{(m_n-m_{1n})}{m_n}e^{-n(J-\epsilon)},
        \end{equation*}
	it follows using (\ref{eq:q}) that
	\begin{equation}
         \frac{\sum_{i=m_{1n}+1}^{m_{n}}w_{in}(\bone)}{m_n}\rightarrow 0, ~\mbox{as}~n\rightarrow\infty.
		\label{eq:eq7}
	\end{equation}
	Hence, dividing the numerator and denominator in the ratio on the right hand side of (\ref{eq:upper_mpBFDR2}) by $m_n$ and using the limits
	(\ref{eq:eq6}), (\ref{eq:eq7}) and (\ref{eq:q}) as $n\rightarrow\infty$, yields
	\begin{equation}
		\lim_{n\rightarrow\infty}\frac{\sum_{i=1}^{m_n} \hat d_i(1-w_{in}(\hat \bd(m_n)))}{\sum_{i=1}^{m_n}\hat d_i} 
		= \frac{1-q}{1+p-q}.
		\label{eq:eq8}
	\end{equation}
	Hence, in this case, the maximum $mpBFDR$ (that can be incurred at $\beta=0$) for $n\rightarrow\infty$ is given by  
	\begin{equation*}
		\lim_{n\rightarrow\infty} mpBFDR_{\beta=0}=\frac{1-q}{1+p-q}.
	\end{equation*}
	Note that this is also the maximum asymptotic $mpBFDR$ that can be incurred among all possible configurations of 
	$G_{m_{1n}+1},\ldots,G_{m_n}$. Hence, for any arbitrary configuration of groups, 
	the maximum asymptotic $mpBFDR$ that can be incurred lies in the interval $\left(0, \frac{1-q}{1+p-q} \right)$.	
\end{proof}

\subsection{Proof of Theorem \ref{corollary:beta_n0}}
\begin{proof} 
	Using the facts that $mpBFDR$ is continuous and decreasing in $\beta$ (\ctn{chandra2019}) and that $mpBFDR$ tends to $0$ (Theorem \ref{theorem:lim_BFDR}),
	the proof follows in the same way as that of Theorem 8 of \ctn{chandra2020}.
\end{proof}

\subsection{Proof of Theorem \ref{theorem:pBFDR_alpha2}}
\begin{proof}
	From \ctn{chandra2019} it is known that $mpBFDR$ and $pBFDR$ are continuous and non-increasing in $\beta$.
	If $\hat\bd(m_n)$ denotes the optimal decision configuration with respect to the additive loss function,  
	$\hat d_i=1$ for all $i$, for $\beta=0$. Thus, assuming without loss of generality that the first $m_{0n}$ null hypotheses are true, 
	\begin{equation}
		\frac{\sum_{i=1}^{m_n} \hat d_i(1-v_{in})}{\sum_{i=1}^{m_n}\hat d_i}
		=1-\frac{\sum_{i=1}^{m_{0n}}v_{in}+\sum_{i=m_{0n}+1}^{m_{n}}v_{in}}{m_n}.
		\label{eq:eq9}
	\end{equation}
	Now, $0<\frac{\sum_{i=1}^{m_{0n}}v_{in}}{m_n}<\left(1-\frac{m_{0n}}{m_n}\right)e^{-n(J-\epsilon)}$, so that 
	$\frac{\sum_{i=1}^{m_{0n}}v_{in}}{m_n}\rightarrow 0$ as $n\rightarrow\infty$.
	Also, $\left(1-e^{-n(J-\epsilon)}\right)\left(1-\frac{m_{0n}}{m_n}\right)<\frac{\sum_{i=m_{0n}+1}^{m_{n}}v_{in}}{m_n}<1-\frac{m_{0n}}{m_n}$,
	so that $\frac{\sum_{i=m_{0n}+1}^{m_{n}}v_{in}}{m_n}\rightarrow p_0$, as $n\rightarrow\infty$.
	Hence, taking limits on both sides of (\ref{eq:eq9}), we obtain
	\begin{equation*}
		\lim_{n\rightarrow\infty}\frac{\sum_{i=1}^{m_n} \hat d_i(1-v_{in})}{\sum_{i=1}^{m_n}\hat d_i}
		=p_0.
	\end{equation*}
	The remaining part of the proof follows 
	in the same way as that of Theorem \ref{corollary:beta_n0}.	
\end{proof}

\subsection{Proof of Theorem \ref{theorem:mpBFDR_alpha2}}
\begin{proof}
	The proof follows in the same way as that of Theorem 10 of \ctn{chandra2020} using the facts $mpBFDR_{\beta}>pBFDR_{\beta}$ for any multiple testing procedure,
	$\underset{n\rightarrow\infty}{\lim}~pBFDR_{\beta=0}=p_0$ (due to Theorem \ref{theorem:pBFDR_alpha2}), and that $mpBFDR$ is continuous and non-increasing in $\beta$ 
	and tends to zero as $n\rightarrow\infty$.
\end{proof}

\subsection{Proof of Theorem \ref{theorem:mpBFNR_alpha}}
\begin{proof}
	Note that by Theorem \ref{corollary:beta_n0}, there exists a sequence $\{\beta_n\}$ such that $\lim_{n\rightarrow\infty} mpBFDR_{\beta_n}=\alpha$,
	where $\alpha\in \left(0,\frac{1-q}{1+p-q}\right)$. 
	Let $\hat{\bd(m_n)}$ be the optimal decision configuration associated with the sequence $\{\beta_n\}$. 
	The proofs of Theorem \ref{theorem:mpBFDR_alpha} and \ref{corollary:beta_n0} show that 
	$\hat d_{in}=d_i^t$ for $i=1,\cdots,m_{1n}$ and $\sum_{i=m_{1n}+1}^{m_n} \hat d_{in}>0$.
	Hence, using (\ref{eq:ineq_end0}) we obtain
	\begin{align}
		&\frac{\sum_{i=1}^{m_n}(1-\hat d_{in})v_{in}} {\sum_{i=1}^{m_n}(1-\hat d_{in})} \leq \frac{\sum_{i=1}^{m_n}(1- d_i^t)v_{in} } {\sum_{i=1}^{m_n}(1-\hat d_{in})} 
		< e^{-n(J-\epsilon)}\times \frac{\sum_{i=1}^{m_n}(1- d_i^t))}{\sum_{i=1}^{m_n}(1-\hat d_{in})}\\
		\Rightarrow~&\frac{1}{n}\log\left(\fnrx \right) <-J+\epsilon+\frac{1}{n} \log\left[\sum_{i=1}^{m_n}(1- d_i^t)\right]
		-\frac{1}{n}\log\left[ \sum_{i=1}^{m_n}(1-\hat d_{in}) \right].\label{eq:eq10} 
	\end{align}
	Now,
	\begin{align}
		0\leq\frac{1}{n}\log\left[ \sum_{i=1}^{m_n}(1-d^t_{i})\right]\leq\frac{\log m_n}{n};\notag\\
		0\leq\frac{1}{n}\log\left[ \sum_{i=1}^{m_n}(1-\hat d_{in})\right]\leq\frac{\log m_n}{n}.\notag
	\end{align}
	Since $\frac{\log m_n}{n}\rightarrow 0$, as $n\rightarrow\infty$,
	\begin{align}
		&\lim_{n\rightarrow\infty}\frac{1}{n}\log\left[ \sum_{i=1}^m(1-d^t_{i})\right]=0,~\mbox{and}\label{eq:eq11}\\
		&\lim_{n\rightarrow\infty} \frac{1}{n} \log\left[\sum_{i=1}^m(1- \hat d_{in})\right] = 0.\label{eq:eq12}
	\end{align}
	As $\epsilon$ is any arbitrary positive quantity we have from (\ref{eq:eq10}), (\ref{eq:eq11}) and (\ref{eq:eq12}) that
	\begin{equation*}
	\limsup_{n\rightarrow\infty} \frac{1}{n}\log\left(\fnrx \right)\leq -J.
	\end{equation*}
\end{proof}

\section[Verification of \ref{shalizi1}-\ref{shalizi7}]{Verification of \ref{shalizi1}-\ref{shalizi7} in $AR(1)$ model with time-varying covariates and proofs of 
the relevant theorems}
\label{sec:ar1_covariates}

All the probabilities and expectations below are with respect to the true model $P$.

\subsection{Verification of \ref{shalizi1}}
\label{subsec:A1}

We obtain
\begin{align}
-\log R_n(\btheta)&=
n\log\left(\frac{\sigma}{\sigma_0}\right)+\left(\frac{1}{2\sigma^2}-\frac{1}{2\sigma^2_0}\right)\sum_{t=1}^nx^2_t
+\left(\frac{\rho^2}{2\sigma^2}-\frac{\rho^2_0}{2\sigma^2_0}\right)\sum_{t=1}^nx^2_{t-1}\notag\\
&\qquad+\frac{1}{2\sigma^2}\bbeta^\prime_m\left(\sum_{t=1}^n\bz_{mt}\bz^\prime_{mt}\right)\bbeta_m 
-\frac{1}{2\sigma^2_0}\bbeta^\prime_{m0}\left(\sum_{t=1}^n\bz_{mt}\bz^\prime_{mt}\right)\bbeta_{m0}\notag\\ 
&\qquad-\left(\frac{\rho}{\sigma^2}-\frac{\rho_0}{\sigma^2_0}\right)\sum_{t=1}^nx_tx_{t-1}
-\left(\frac{\bbeta_m}{\sigma^2}-\frac{\bbeta_{m0}}{\sigma^2_0}\right)^\prime\sum_{t=1}^n\bz_{mt}x_t\notag\\
&\qquad+\left(\frac{\rho\bbeta_m}{\sigma^2}-\frac{\rho_0\bbeta_0}{\sigma^2_0}\right)'\sum_{t=1}^n\bz_{mt}x_{t-1}.
\label{eq:RT}
\end{align}

It is easily seen that $-\log R_n(\btheta)$ is continuous in $\bX_n$ and $\btheta$.
Hence, $R_n(\btheta)$ is $\mathcal F_n\times\mathcal T$
measurable. In other words, (S1) holds.

\subsection{Proof of Lemma \ref{lemma:lemma1_inf}}
\label{subsec:A2_inf}

It is easy to see that under the true model $P$,
\begin{align}
E(x_t)&=\sum_{k=1}^t\rho^{t-k}_0\bz^\prime_{mk}\bbeta_{m0};
\label{eq:mean_true_inf}\\
E(x_{t+h}x_t)&\sim\frac{\sigma^2_0\rho^h_0}{1-\rho^2_0}+E(x_{t+h})E(x_t);~h\geq 0,
\label{eq:cov_true}
\end{align}
where for any two sequences $\{a_t\}_{t=1}^{\infty}$ and $\{b_t\}_{t=1}^{\infty}$, $a_t\sim b_t$ stands for $a_t/b_t\rightarrow 1$ as $t\rightarrow\infty$.
Hence, 
\begin{equation}
E(x^2_t)\sim \frac{\sigma^2_0}{1-\rho^2_0}+\left(\sum_{k=1}^t\rho^{t-k}_0\bz^\prime_{mk}\bbeta_{m0}\right)^2.
\label{eq:lim_mean_sq_inf}
\end{equation}

Now let 
\begin{equation}
\varrho_t=\sum_{k=1}^t\rho^{t-k}_0\bz^\prime_{mk}\bbeta_{m0} 
\label{eq:varrho_t}
\end{equation}
and for $t>t_0$, 
\begin{equation}
\tilde\varrho_t=\sum_{k=t-t_0}^t\rho^{t-k}_0\bz^\prime_{mk}\bbeta_{m0}, 
\label{eq:tilde_varrho_t}
\end{equation}
where, for any $\varepsilon>0$, $t_0$ is so large that 
\begin{equation}
\frac{C\left|\rho_0\right|^{t_0+1}}{(1-\left|\rho_0\right|^{t_0})}\leq\varepsilon.
\label{eq:small1_inf}
\end{equation}
It follows, using (\ref{eq:B2}) and (\ref{eq:small1_inf}), that for $t>t_0$,
\begin{equation}
\left|\varrho_t-\tilde\varrho_t\right|\leq \sum_{k=1}^{t-t_0-1}|\rho_0|^{t-k}\left|\bz^\prime_{mk}\bbeta_{m0}\right|\leq\frac{C|\rho_0|^{t_0+1}(1-|\rho_0|^{t-t_0+1})}{1-|\rho_0|}\leq\varepsilon.
\label{eq:diff_inequality1_inf}
\end{equation}
Hence, for $t>t_0$,
\begin{equation}
\tilde\varrho_t-\varepsilon\leq\varrho_t\leq\tilde\varrho_t+\varepsilon.
\label{eq:diff_inequality2_inf}
\end{equation}
Now,
\begin{align}
\frac{\sum_{t=1}^n\tilde\varrho_t}{n}&=\rho^{t_0}_0\left(\frac{\sum_{t=1}^n\bz_{mt}}{n}\right)^\prime\bbeta_{m0}+\rho^{t_0-1}_0\left(\frac{\sum_{t=2}^n\bz_{mt}}{n}\right)^\prime\bbeta_{m0}
+\rho^{t_0-2}_0\left(\frac{\sum_{t=3}^n\bz_{mt}}{n}\right)^\prime\bbeta_{m0}+\cdots\notag\\
&\qquad\qquad\cdots+\rho_0\left(\frac{\sum_{t=t_0}^n\bz_{mt}}{n}\right)^\prime\bbeta_{m0}+\left(\frac{\sum_{t=t_0+1}^n\bz_{mt}}{n}\right)^\prime\bbeta_{m0}\notag\\
	&\rightarrow 0,~\mbox{as}~n\rightarrow\infty,~\mbox{by virtue of (B4) (\ref{eq:ass1})}.
\label{eq:tilde_varrho_limit}
\end{align}
Similarly, it is easily seen, using (B4), that
\begin{equation}
\frac{\sum_{t=1}^n\tilde\varrho^2_t}{n}\rightarrow \left(\frac{1-\rho^{2(2t_0+1)}_0}{1-\rho^2_0}\right)c(\bbeta_0), 
~\mbox{as}~n\rightarrow\infty.
\label{eq:tilde_varrho_sq_limit}
\end{equation}
Since (\ref{eq:diff_inequality1_inf}) implies that for $t>t_0$, $\tilde \varrho^2_t+\varepsilon^2-2\varepsilon\tilde \varrho_t\leq \varrho^2_t
\leq\tilde \varrho^2_t+\varepsilon^2+2\varepsilon\tilde \varrho_t$,
it follows that
\begin{equation}
\underset{n\rightarrow\infty}{\lim}~\frac{\sum_{t=1}^n\varrho^2_t}{n}=\underset{n\rightarrow\infty}{\lim}~\frac{\sum_{t=1}^n\tilde\varrho^2_t}{n}+\varepsilon^2
=\left(\frac{1-\rho^{2(2t_0+1)}_0}{1-\rho^2_0}\right)c(\bbeta_0) 
+\varepsilon^2,
\label{eq:varrho_sq_limit}
\end{equation}
and since $\epsilon>0$ is arbitrary, it follows that
\begin{equation}
\underset{n\rightarrow\infty}{\lim}~\frac{\sum_{t=1}^n\varrho^2_t}{n}
=\frac{c(\bbeta_0)}{1-\rho^2_0}.
\label{eq:varrho_sq_limit2_inf}
\end{equation}
Hence, it also follows from (\ref{eq:mean_true_inf}), (\ref{eq:lim_mean_sq_inf}), (B4) and (\ref{eq:varrho_sq_limit2_inf}), that 

\begin{equation}
\frac{\sum_{t=1}^nE(x^2_t)}{n}\rightarrow \frac{\sigma^2_0}{1-\rho^2_0}+\frac{c(\bbeta_0)}{1-\rho^2_0},~\mbox{as}~n\rightarrow\infty
\label{eq:lim2}
\end{equation}
and
\begin{equation}
\frac{\sum_{t=1}^nE(x^2_{t-1})}{n}\rightarrow \frac{\sigma^2_0}{1-\rho^2_0}+\frac{c(\bbeta_0)}{1-\rho^2_0},~\mbox{as}~n\rightarrow\infty.
\label{eq:lim3_inf}
\end{equation}
Now note that 
\begin{equation}
x_tx_{t-1}=\rho_0x^2_{t-1}+\bz^\prime_{mt}\bbeta_0x_{t-1}+\epsilon_tx_{t-1}.
\label{eq:cov_lag1_inf}
\end{equation}
Using (\ref{eq:ass2_inf}), (\ref{eq:diff_inequality2_inf}) and arbitrariness of $\varepsilon>0$ it is again easy to see that
\begin{equation}
\frac{\sum_{t=1}^n\bz^\prime_{mt}\bbeta_{m0}E(x_{t-1})}{n}\rightarrow 0,~\mbox{as}~n\rightarrow\infty.
\label{eq:lim4_inf}
\end{equation}
Also, since for $t=1,2,\ldots,$ $E(\epsilon_tx_{t-1})=E(\epsilon_t)E(x_{t-1})$ by independence, and since $E(\epsilon_t)=0$ for $t=1,2,\ldots$, it holds that
\begin{equation}
\frac{\sum_{t=1}^nE\left(\epsilon_tx_{t-1}\right)}{n}= 0,~\mbox{for all}~n=1,2,\ldots.
\label{eq:lim5_inf}
\end{equation}
Combining (\ref{eq:cov_lag1_inf}), (\ref{eq:lim3_inf}), (\ref{eq:lim4_inf}) and (\ref{eq:lim5_inf}) we obtain
\begin{equation}
\frac{\sum_{t=1}^nE\left(x_tx_{t-1}\right)}{n}\rightarrow \frac{\rho_0\sigma^2_0}{1-\rho^2_0}+\frac{\rho_0c(\bbeta_0)}{1-\rho^2_0}.
\label{eq:lim6}
\end{equation}
Using (B4)  (\ref{eq:B2}) and arbitrariness of $\varepsilon>0$, it follows that
\begin{multline*}
h(\btheta)=\underset{n\rightarrow\infty}{\lim}~\frac{1}{n}E\left[-\log R_n(\btheta)\right]
=\log\left(\frac{\sigma}{\sigma_0}\right)
+\left(\frac{1}{2\sigma^2}-\frac{1}{2\sigma^2_0}\right)\left(\frac{\sigma^2_0}{1-\rho^2_0}
+\frac{c(\bbeta_0)}{1-\rho^2_0}\right)\\
+\left(\frac{\rho^2}{2\sigma^2}-\frac{\rho^2_0}{2\sigma^2_0}\right)\left(\frac{\sigma^2_0}{1-\rho^2_0}+\frac{c(\bbeta_0)}{1-\rho^2_0}\right)
+\frac{c(\bbeta)}{2\sigma^2}-\frac{c(\bbeta_0)}{2\sigma^2_0}\\
-\left(\frac{\rho}{\sigma^2}-\frac{\rho_0}{\sigma^2_0}\right)\left(\frac{\rho_0\sigma^2_0}{1-\rho^2_0}+\frac{\rho_0c(\bbeta_0)}{1-\rho^2_0}\right)-\left(\frac{c_{10}(\bbeta,\bbeta_0)}{\sigma^2}-\frac{c(\bbeta_0)}{\sigma^2_0}\right).
\end{multline*}
In other words, (S2) holds, with $h(\btheta)$ given by (\ref{eq:htheta_inf}).

\subsection{Proof of Theorem \ref{theorem:theorem2_inf}}
\label{subsec:A3_inf}

Note that 
\begin{equation}
x_t=\sum_{k=1}^t\rho^{t-k}_0\bz^\prime_{mk}\bbeta_{m0}+\sum_{k=1}^t\rho^{t-k}_0\epsilon_k,
\label{eq:xt}
\end{equation}
where $\tilde\epsilon_t=\sum_{k=1}^t\rho^{t-k}_0\epsilon_k$ is an asymptotically stationary Gaussian process with mean zero and covariance
\begin{equation}
cov(\tilde\epsilon_{t+h},\tilde\epsilon_t)\sim\frac{\sigma^2_0\rho^h_0}{1-\rho^2_0},~\mbox{where}~h\geq 0.
\label{eq:cov_epsilon}
\end{equation}
Then 
\begin{equation}
\frac{\sum_{t=1}^nx^2_t}{n}=\frac{\sum_{t=1}^n\varrho^2_t}{n}+\frac{\sum_{t=1}^n\tilde\epsilon^2_t}{n}+
\frac{2\sum_{t=1}^n\tilde\epsilon_t\varrho_t}{n}.
\label{eq:xt_sq_average_inf}
\end{equation}
By (\ref{eq:varrho_sq_limit2_inf}), the first term of the right hand side of (\ref{eq:xt_sq_average_inf}) converges to 
$\frac{c(\bbeta_0)}{1-\rho^2_0}$, 
as $n\rightarrow\infty$, and since $\tilde\epsilon_t$; $t=1,2,\ldots$
is also an irreducible and aperiodic Markov chain, by the ergodic theorem it follows that the second term of (\ref{eq:xt_sq_average_inf}) converges to $\sigma^2_0/(1-\rho^2_0)$
almost surely, as $n\rightarrow\infty$.
For the third term, we observe that 
\begin{equation}
|\bz^\prime_k\bbeta_0-\bz^\prime_{mk}\bbeta_{m0}|<\delta, 
\label{eq:inf_approx}
\end{equation}
for $n>n_0$, where
$n_0$, depending upon $\delta~(>0)$, is sufficiently large. Recalling from (B5) that $\hat\varrho_t=\sum_{k=1}^t\rho^{t-k}_0\bz^\prime_k\bbeta_0$, we then see that for $t>n_0$,
\begin{equation}
|\varrho_t-\hat\varrho_t|<\frac{\delta}{1-|\rho_0|}<\varepsilon,
\label{eq:hat_varrho}
\end{equation}
for $\delta<(1-|\rho_0|)\varepsilon$. 
From (\ref{eq:hat_varrho}) it follows that
\begin{align}
\underset{n\rightarrow\infty}{\lim}~\frac{2\sum_{t=1}^n\tilde\epsilon_t\varrho_t}{n}=\underset{n\rightarrow\infty}{\lim}~\frac{2\sum_{t=n_0+1}^n\tilde\epsilon_t\hat\varrho_t}{n-n_0}
\end{align}
Since by (B5) 
the limit of $\hat\varrho_t$ exists as $t\rightarrow\infty$, it follows that $\tilde\epsilon_t\hat\varrho_t$ is still an irreducible and aperiodic Markov chain with 
asymptotically stationary zero-mean Gaussian process.
%
Hence, by the ergodic theorem, the third term of (\ref{eq:xt_sq_average_inf}) converges to zero, almost surely, as $n\rightarrow\infty$. It follows that
\begin{equation}
\frac{\sum_{t=1}^nx^2_t}{n}\rightarrow\frac{\sigma^2_0}{1-\rho^2_0}+\frac{c(\bbeta_0)}{1-\rho^2_0},
\label{eq:xt_sq_average2_inf}
\end{equation}
and similarly,
\begin{equation}
\frac{\sum_{t=1}^nx^2_{t-1}}{n}\rightarrow\frac{\sigma^2_0}{1-\rho^2_0}+\frac{c(\bbeta_0)}{1-\rho^2_0}.
\label{eq:xt_sq_average3_inf}
\end{equation}

Now, since $x_t=\varrho_t+\tilde\epsilon_t$, it follows using (B2) (orthogonality) and (\ref{eq:diff_inequality2_inf}) that for 
$\tilde\bbeta_m=\bbeta_m$ or $\tilde\bbeta_m=\bbeta_{m0}$,
\begin{equation}
\underset{n\rightarrow\infty}{\lim}~\frac{\sum_{t=1}^n\tilde\bbeta^\prime_m\bz_{mt}x_t}{n}
=\underset{n\rightarrow\infty}{\lim}~\frac{\sum_{t=1}^n\tilde\bbeta^\prime_m\bz_{mt}\bz^\prime_{mt}\bbeta_{m0}}{n}
+\underset{n\rightarrow\infty}{\lim}~\frac{\sum_{t=1}^n\tilde\bbeta^\prime_m\bz_{mt}\tilde\epsilon_t}{n}.
\label{eq:z_tx_t_average_inf}
\end{equation}
By (B4), the first term on the right hand side of (\ref{eq:z_tx_t_average_inf}) is $\tilde c(\bbeta,\bbeta_0)$, where $\tilde c(\bbeta,\bbeta_0)$ is $c(\bbeta_0)$ or $c_{10}(\bbeta,\bbeta_0)$ 
accordingly as
$\tilde\bbeta_m$ is $\bbeta_{m0}$ or $\bbeta_m$. 
For the second term, due to (\ref{eq:inf_approx}), $\underset{n\rightarrow\infty}{\lim}~\frac{\sum_{t=1}^n\tilde\bbeta^\prime_m\bz_{mt}\tilde\epsilon_t}{n}
=\underset{n\rightarrow\infty}{\lim}~\frac{\sum_{t=1}^n\tilde\bbeta^\prime\bz_{t}\tilde\epsilon_t}{n}$, where $\tilde\bbeta$ is either $\bbeta$ or $\bbeta_0$.
By (B5) the limit of $\tilde\bbeta^\prime\bz_{t}$ exists as $t\rightarrow\infty$, and hence $\tilde\bbeta^\prime\bz_{t}\tilde\epsilon_t$ remains an irreducible, aperiodic
Markov chain with zero mean Gaussian stationary distribution.
Hence, by the ergodic theorem, it follows that
the second term of (\ref{eq:z_tx_t_average_inf}) is zero, almost surely. In other words, almost surely,
\begin{equation}
\frac{\sum_{t=1}^n\tilde\bbeta^\prime_m\bz_{mt}x_t}{n}\rightarrow \tilde c(\bbeta,\bbeta_0),~\mbox{as}~n\rightarrow\infty,
\label{eq:z_tx_t_average2_inf}
\end{equation}
and similar arguments show that, almost surely,
\begin{equation}
\frac{\sum_{t=1}^n\tilde\bbeta^\prime_m\bz_{mt}x_{t-1}}{n}\rightarrow 0,~\mbox{as}~n\rightarrow\infty.
\label{eq:z_tx_t_average3_inf}
\end{equation}
We now calculate the limit of $\sum_{t=1}^nx_tx_{t-1}/n$, as $n\rightarrow\infty$. By (\ref{eq:cov_lag1_inf}), 
\begin{equation}
\underset{n\rightarrow\infty}{\lim}~\frac{\sum_{t=1}^nx_tx_{t-1}}{n}=\underset{n\rightarrow\infty}{\lim}~\frac{\rho_0\sum_{t=1}^nx^2_{t-1}}{n}
+\underset{n\rightarrow\infty}{\lim}~\frac{\bbeta^\prime_{m0}\sum_{t=1}^n\bz_{mt}x_{t-1}}{n}
+\underset{n\rightarrow\infty}{\lim}~\frac{\sum_{t=1}^n\epsilon_tx_{t-1}}{n}.
\label{eq:cov_lag1_limit_inf}
\end{equation}
By (\ref{eq:xt_sq_average3_inf}), the first term on the right hand side of (\ref{eq:cov_lag1_limit_inf}) is given, almost surely, by 
$\frac{\rho_0\sigma^2_0}{1-\rho^2_0}+\frac{\rho_0c(\bbeta_0)}{1-\rho^2_0}$, 
and the second term is almost surely zero due to (\ref{eq:z_tx_t_average3_inf}).
For the third term, note that $\epsilon_tx_{t-1}=\epsilon_t\varrho_{t-1}+\epsilon_t\tilde\epsilon_{t-1}$, and hence using (\ref{eq:inf_approx}),
$\underset{n\rightarrow\infty}{\lim}~\frac{\sum_{t=1}^n\epsilon_tx_{t-1}}{n}=\underset{n\rightarrow\infty}{\lim}~\frac{\sum_{t=1}^n\epsilon_t\hat\varrho_{t-1}}{n}
+\underset{n\rightarrow\infty}{\lim}~\frac{\sum_{t=1}^n\epsilon_t\tilde\epsilon_{t-1}}{n}$.
Both $\epsilon_t\hat\varrho_{t-1}$; $t=1,2,\ldots$
and $\epsilon_t\tilde\epsilon_{t-1}$; $t=1,2,\ldots$, are sample paths of irreducible and aperiodic Markov chains having stationary distributions with mean zero.
Hence, by the ergodic theorem, the third term of (\ref{eq:cov_lag1_limit_inf}) is zero, almost surely. That is,
\begin{equation}
\underset{n\rightarrow\infty}{\lim}~\frac{\sum_{t=1}^nx_tx_{t-1}}{n}=\frac{\rho_0\sigma^2_0}{1-\rho^2_0}+\frac{\rho_0c(\bbeta_0)}{1-\rho^2_0}.
\label{eq:cov_lag1_limit2_inf}
\end{equation}

The limits (\ref{eq:xt_sq_average2_inf}), (\ref{eq:xt_sq_average3_inf}), (\ref{eq:z_tx_t_average2_inf}), (\ref{eq:z_tx_t_average3_inf}), (\ref{eq:cov_lag1_limit2_inf})
applied to $\log R_n(\btheta)$ given by (\ref{eq:RT}), shows that $\frac{\log R_n(\btheta)}{n}$ converges to $-h(\btheta)$ almost surely as $n\rightarrow\infty$.
In other words, (S3) holds.

\subsection{Verification of \ref{s4}}
\label{subsec:A4}

In the expression for $h(\btheta)$ given by (\ref{eq:htheta_inf}), note that $c(\bbeta)$ and $c_{10}(\bbeta,\bbeta_0)$ are almost surely finite. Hence, for any prior on $\sigma$ and $\rho$ such that they 
are almost surely finite, (S4) clearly holds. 
In particular, this holds for any proper priors on $\sigma$ and $\rho$.

\subsection{Verification of \ref{s5}}
\label{subsec:A5}

\subsubsection{Verification of (S5) (1)}
\label{subsubsec:A5_1}

Since $\Theta=\mathcal C'(\mathcal X)\times\mathbb R^{\infty}\times\mathbb R\times\mathbb R^+$, it is easy
to see that $h(\Theta)=0$.
Let $\bgamma_m=(\gamma_1,\ldots,\gamma_m)^\prime$, $\tilde\gamma_m=\sum_{i=1}^m|\gamma_i|$, $\btheta_m=(\tilde\eta,\bgamma_m,\rho,\sigma)$,
$\Theta_m=\mathcal C'(\mathcal X)\times\mathbb R^m\times\mathbb R\times\mathbb R^+$.
We now define
\begin{align*}
	\mathcal G_n &=\left\{\btheta_m\in\Theta_m:|\rho|\leq\exp\left(\left(\alpha n\right)^{1/16}\right),
	\tilde\gamma_m\leq\exp\left(\left(\alpha n\right)^{1/16}\right),\right.\notag\\
	&\qquad\left.\|\tilde\eta\|\leq\exp\left(\left(\alpha n\right)^{1/16}\right),\|\tilde\eta'\|\leq\exp\left(\left(\alpha n\right)^{1/16}\right),
	\exp\left(-\left(\alpha n\right)^{1/16}\right)\leq\sigma\leq\exp\left(\left(\alpha n\right)^{1/16}\right)\right\},
\end{align*}
where
$\alpha>0$. 

Since $|\gamma_i|<L<\infty$ for all $i$, it follows that $\mathcal G_n$ is increasing in $n$ for $n\geq n_0$, for some $n_0\geq 1$. To see this, note that
if $\tilde\gamma_{m_n}\leq\exp(\left(\alpha n\right)^{1/16})$, then $\tilde\gamma_{m_{n+1}}=\tilde\gamma_{m_n}+\sum_{i=m_n+1}^{m_{n+1}}|\gamma_i|
<\exp(\left(\alpha (n+1)\right)^{1/16})$
if $\sum_{i=m_n+1}^{m_{n+1}}|\gamma_i|<L(m_{n+1}-m_n)<\exp(\left(\alpha (n+1)\right)^{1/16})-\exp(\left(\alpha n\right)^{1/16})$, which holds by assumption (B7).
Since $\mathcal G_n\rightarrow\Theta$ as $n\rightarrow\infty$, there exists $n_1$ such that $\mathcal G_{n_1}$ contains $\btheta_0$. Hence, $h(\mathcal G_n)=0$
for all $n\geq n_1$. In other words, $h(\mathcal G_n)\rightarrow h(\Theta)$, as $n\rightarrow\infty$.
Now observe that 
\begin{align}
&\pi\left(\mathcal G_n\right)\notag\\
	&=\pi\left(\tilde\gamma_m\leq\exp(\left(\alpha n\right)^{1/16}),\|\tilde\eta\|\leq\exp(\left(\alpha n\right)^{1/16}),
	\|\tilde\eta'\|\leq\exp(\left(\alpha n\right)^{1/16}),\right.\notag\\
	&\qquad\left.\exp\left(-\left(\alpha n\right)^{1/16}\right)\leq\sigma\leq\exp\left(\left(\alpha n\right)^{1/16}\right)\right)\notag\\
	&\quad-\pi\left(|\rho|>\exp\left(\left(\alpha n\right)^{1/16}\right),\tilde\gamma_m\leq\exp(\left(\alpha n\right)^{1/16}),
	\|\tilde\eta\|\leq\exp(\left(\alpha n\right)^{1/16}),\|\tilde\eta'\|\leq\exp(\left(\alpha n\right)^{1/16}),\right.\notag\\
	&\qquad\left.\exp\left(-\left(\alpha n\right)^{1/16}\right)\leq\sigma\leq\exp\left(\left(\alpha n\right)^{1/16}\right)\right)\notag\\
	&\geq 1-\pi\left(|\rho|>\exp\left(\left(\alpha n\right)^{1/16}\right)\right)-\pi\left(\tilde\gamma_m>\exp\left(\left(\alpha n\right)^{1/16}\right)\right)
	-\pi\left(\|\tilde\eta\|>\exp(\left(\alpha n\right)^{1/16})\right)\notag\\
	&\qquad-\pi\left(\|\tilde\eta'\|> \exp(\left(\alpha n\right)^{1/16})\right)
	-\pi\left(\left\{\exp\left(-\left(\alpha n\right)^{1/16}\right)\leq\sigma\leq \exp\left(\left(\alpha n\right)^{1/16}\right)\right\}^c\right)\notag\\
&\geq 1-(c_{\rho}+c_{\gamma}+c_{\tilde\eta}+c_{\tilde\eta'}+c_{\sigma})\exp(-\alpha n),\notag
\end{align}
where the last step is due to (B6).

\subsubsection{Verification of (S5) (2)}
\label{subsubsec:A5_2}

First, we note that $\mathcal G_n$ is compact, which can be proved using Arzela-Ascoli lemma, in almost the same way as in \ctn{Chatterjee18b}.
Since $\mathcal G_n$ is compact for all $n\geq 1$, uniform convergence as required will be proven if we can show that $\frac{1}{n}\log R_n(\btheta)+h(\btheta)$ is  
stochastically equicontinuous almost surely in $\btheta\in\mathcal G$ for any 
$\mathcal G\in\left\{\mathcal G_n:n=1,2,\ldots\right\}$ and $\frac{1}{n}\log R_n(\btheta)+h(\btheta)\rightarrow 0$, almost surely, for all $\btheta\in\mathcal G$
(see \ctn{Newey91} for the general theory of uniform convergence in compact sets under stochastic equicontinuity). Since we have already verified pointwise convergence
of the above for all $\btheta\in\bTheta$ while verifying (S3), it remains to prove stochastic equicontinuity of $\frac{1}{n}\log R_n(\cdot)+h(\cdot)$.
Stochastic equicontinuity usually follows easily if one can prove that the function concerned is almost surely Lipschitz continuous. In our case, 
we can first verify Lipschitz continuity of $\frac{1}{n}\log R_n(\btheta)$ by showing that its first partial derivatives 
with respect to the components of $\btheta$ are almost surely bounded. With respect to $\rho$ and $\sigma$, the boundedness of the parameters in $\mathcal G$, (\ref{eq:B2}) and the limit results 
(\ref{eq:xt_sq_average2_inf}), (\ref{eq:xt_sq_average3_inf}), (\ref{eq:z_tx_t_average2_inf}), (\ref{eq:z_tx_t_average3_inf}) and (\ref{eq:cov_lag1_limit2_inf})
readily show boundedness of the partial derivatives. With respect to $\bbeta_m$, note that the derivative of 
$\frac{1}{2\sigma^2}\bbeta^\prime_m\left(\frac{\sum_{t=1}^n\bz_{mt}\bz^\prime_{mt}}{n}\right)\bbeta_m$,
a relevant expression of $\frac{1}{n}\log R_n(\btheta)$ (see (\ref{eq:RT})), is $\frac{1}{\sigma^2}\left(\frac{\sum_{t=1}^n\bz_{mt}\bz^\prime_{mt}}{n}\right)\bbeta_m$, 
whose Euclidean norm is bounded above by $\sigma^{-2}\|\left(\frac{\sum_{t=1}^n\bz_{mt}\bz^\prime_{mt}}{n}\right)\|_{op}\times\|\bbeta_m\|$.
In our case, $\|\left(\frac{\sum_{t=1}^n\bz_{mt}\bz^\prime_{mt}}{n}\right)\|_{op}\leq K<\infty$ by (B3). Moreover, $\sigma^{-2}$ is bounded in $\mathcal G$ 
and $\|\bbeta_m\|\leq\|\tilde\eta\|\times\sqrt{\sum_{i=1}^m\gamma^2_i}$, which is also bounded in $\mathcal G$. Boundedness of the partial derivatives with respect to $\bbeta_m$ of the other terms of $\frac{1}{n}\log R_n(\btheta)$ 
involving $\bbeta_m$ are easy to observe. In other words, $\frac{1}{n}\log R_n(\btheta)$ is stochastically equicontinuous. 

To see that $h(\btheta)$ is equicontinuous, first note that
in the expression (\ref{eq:htheta_inf}), except the terms involving $c(\bbeta)$ and $c_{10}(\bbeta,\bbeta_0)$, the other terms are easily seen to be Lipschitz, using boundedness of the partial derivatives. Let us
now focus on the term $\frac{c(\bbeta)}{2\sigma^2}$.
For our purpose, let us consider two different sequences $\bbeta_{1m}$ and $\bbeta_{2m}$ associated with $(\gamma_1,\tilde\eta_1)$ and $(\gamma_2,\tilde\eta_2)$, 
 respectively, such that $\bbeta^\prime_{1m}\left(\frac{\sum_{t=1}^n\bz_{mt}\bz^\prime_{mt}}{n}\right)\bbeta_{1m}\rightarrow c(\bbeta_1)$
and $\bbeta^\prime_{2m}\left(\frac{\sum_{t=1}^n\bz_{mt}\bz^\prime_{mt}}{n}\right)\bbeta_{2m}\rightarrow c(\bbeta_2)$. 
As we have already shown that 
$\bbeta^\prime_m\left(\frac{\sum_{t=1}^n\bz_{mt}\bz^\prime_{mt}}{n}\right)\bbeta_m$ is Lipschitz in $\bbeta_m$, 
we must have $\|\bbeta^\prime_{1m}\left(\frac{\sum_{t=1}^n\bz_{mt}\bz^\prime_{mt}}{n}\right)\bbeta_{1m}-\bbeta^\prime_{2m}\left(\frac{\sum_{t=1}^n\bz_{mt}\bz^\prime_{mt}}{n}\right)\bbeta_{2m}\|
\leq L\|\bbeta_{1m}-\bbeta_{2m}\|\leq L\|\gamma_1\tilde\eta_1-\gamma_2\tilde\eta_2\|$, for some Lipschitz constant $L>0$. Taking the limit of both sides as $n\rightarrow\infty$
shows that $|c(\bbeta_1)-c(\bbeta_2)|\leq L\|\gamma_1\tilde\eta_1-\gamma_2\tilde\eta_2\|$, proving that $\frac{c(\bbeta)}{2\sigma^2}$ is Lipschitz in $\eta=\gamma\tilde\eta$, when $\sigma$ is held fixed. 
The bounded partial derivative with respect to $\sigma$ also shows that $\frac{c(\bbeta)}{2\sigma^2}$ is Lipschitz in both $\eta$ and $\sigma$. Similarly, the term $\frac{c_{10}(\bbeta,\bbeta_0)}{\sigma^2}$
in (\ref{eq:htheta_inf}) is also Lipschitz continuous. 

In other words, $\frac{1}{n}\log R_n(\btheta)+h(\btheta)$ is stochastically equicontinuous almost surely in $\btheta\in\mathcal G$.
Hence, the required uniform convergence is satisfied.

\subsubsection{Verification of (S5) (3)}
\label{subsubsec:A5_3}

Continuity of $h(\btheta)$, compactness of $\mathcal G_n$ , along with its non-decreasing
nature with respect to $n$ implies that $h\left(\mathcal G_n\right)\rightarrow h\left(\bTheta\right)$, as $n\rightarrow\infty$.
Hence, (S5) holds.

\subsection{Verification of \ref{s6} and proof of Theorem \ref{theorem:theorem3}}
\label{subsec:A6}

Note that in our case,
\begin{align}
&\frac{1}{n}\log R_n(\btheta)+h(\btheta)\notag\\
&\qquad=\left(\frac{1}{2\sigma^2}-\frac{1}{2\sigma^2_0}\right)\left(\frac{\sum_{t=1}^nx^2_t}{n}-\frac{\sigma^2_0}{1-\rho^2_0}-\frac{c(\bbeta_0)}{1-\rho^2_0}\right)\notag\\
&\qquad+\left(\frac{\rho^2}{2\sigma^2}-\frac{\rho^2_0}{2\sigma^2_0}\right)\left(\frac{\sum_{t=1}^nx^2_{t-1}}{n}-\frac{\sigma^2_0}{1-\rho^2_0}-\frac{c(\bbeta_0)}{1-\rho^2_0}\right)
\notag\\
&\qquad+\frac{1}{2\sigma^2}\left(\bbeta^\prime_m\left(\frac{\sum_{t=1}^n\bz_{mt}\bz^\prime_{mt}}{n}\right)\bbeta_m-c(\bbeta)\right)
-\frac{1}{2\sigma^2_0}\left(\bbeta^\prime_{m0}\left(\frac{\sum_{t=1}^n\bz_{mt}\bz^\prime_{mt}}{n}\right)\bbeta_{m0}-c(\bbeta_0)\right)\notag\\
&\qquad-\left(\frac{\rho}{\sigma^2}-\frac{\rho_0}{\sigma^2_0}\right)\left(\frac{\rho_0\sum_{t=1}^nx^2_{t-1}}{n}
+\frac{\bbeta^\prime_{m0}\sum_{t=1}^n\bz_{mt}x_{t-1}}{n}-\frac{\rho_0\sigma^2_0}{1-\rho^2_0}
-\frac{\rho_0c(\bbeta_0)}{1-\rho^2_0}\right)\notag\\
&\qquad-\left[\left(\frac{\bbeta_m}{\sigma^2}-\frac{\bbeta_{m0}}{\sigma^2_0}\right)^\prime\left(\frac{\sum_{t=1}^n\bz_{mt}x_t}{n}\right)-\frac{c_{10}(\bbeta,\bbeta_0)}{\sigma^2}
+\frac{c(\bbeta_0)}{\sigma^2_0}\right]\notag\\
&\qquad +\left(\frac{\rho\bbeta_m}{\sigma^2}-\frac{\rho_0\bbeta_0}{\sigma^2_0}\right)^\prime\frac{\sum_{t=1}^n\bz_{mt}x_{t-1}}{n}\notag\\
&\qquad +\left(\frac{\rho}{\sigma^2}-\frac{\rho_0}{\sigma^2_0}\right)\left(\frac{\sum_{t=1}^n\epsilon_tx_{t-1}}{n}\right).
\label{eq:U}
\end{align}

Let $\kappa_1=(\kappa-h\left(\bTheta\right))/7$, $\bmu_n=E(\bx_n)$ and $\bSigma_n=Var(\bx_n)$; let $\bSigma_n=\bC_n\bC^\prime_n$ be the Cholesky decomposition. 
Also let $\by_n\sim N_n\left(\bzero_n,\bI_n\right)$, the $n$-dimensional normal distribution
with mean $\bzero_n$, the $n$-dimensional vector with all components zero and variance $\bI_n$, the $n$-dimensional identity matrix. Then
\begin{align}
&P\left(\left|\frac{1}{2\sigma^2}-\frac{1}{2\sigma^2_0}\right|
\left|\frac{\sum_{t=1}^nx^2_t}{n}-\frac{\sigma^2_0}{1-\rho^2_0}-\frac{c(\bbeta_0)}{1-\rho^2_0}\right|>\kappa_1\right)\notag\\
=&P\left(\left|\frac{1}{2\sigma^2}-\frac{1}{2\sigma^2_0}\right|
\left|\frac{\bmu_n\prime\bmu_n+2\bmu^\prime_n\bC_n\by_n+\by^\prime_n\bSigma_n\by_n}{n}-\frac{\sigma^2_0}{1-\rho^2_0}-\frac{c(\bbeta_0)}{1-\rho^2_0}\right|>\kappa_1\right)\notag\\
\leq &P\left(\left|\frac{1}{2\sigma^2}-\frac{1}{2\sigma^2_0}\right|\left|\frac{2\bmu^\prime_n\bC_n\by_n}{n}\right|>\frac{\kappa_1}{4}\right)
+P\left(\left|\frac{1}{2\sigma^2}-\frac{1}{2\sigma^2_0}\right|\left|\frac{\bmu^\prime_n\bmu_n}{n}-\frac{c(\bbeta_0)}{1-\rho^2_0}\right|>\frac{\kappa_1}{4}\right)
\label{eq:ce1}\\
&+P\left(\left|\frac{1}{2\sigma^2}-\frac{1}{2\sigma^2_0}\right|\left|\frac{\by^\prime_n\bSigma_n\by_n}{n}-tr\left(\frac{\bSigma_n}{n}\right)\right|>\frac{\kappa_1}{4}\right)
+P\left(\left|\frac{1}{2\sigma^2}-\frac{1}{2\sigma^2_0}\right|\left|tr\left(\frac{\bSigma_n}{n}\right)-\frac{\sigma^2_0}{1-\rho^2_0}\right|>\frac{\kappa_1}{4}\right).
\label{eq:ce2}
\end{align}
To deal with the first term of (\ref{eq:ce1}) first note that $2\bmu^\prime_n\bC_n\by_n$ is Lipschitz in $\by_n$, with the square of the Lipschitz constant
being $4\bmu^\prime_n\bSigma_n\bmu_n$, which is again bounded above by $K_1n$, for some constant $K_1>0$, due to (\ref{eq:B2}). It then follows using the Gaussian concentration inequality
(see, for example, \ctn{Giraud15}) that
\begin{align}
P\left(\left|\frac{1}{2\sigma^2}-\frac{1}{2\sigma^2_0}\right|\left|\frac{2\bmu^\prime_n\bC_n\by_n}{n}\right|>\frac{\kappa_1}{4}\right)
&=P\left(\left|2\bmu^\prime_n\bC_n\by_n\right|>\frac{n\kappa_1}{3}\left|\frac{1}{2\sigma^2}-\frac{1}{2\sigma^2_0}\right|^{-1}\right)\notag\\
&\leq 2\exp\left(-\frac{n\kappa^2_1}{18K_1}\left|\frac{1}{2\sigma^2}-\frac{1}{2\sigma^2_0}\right|^{-2}\right).
\label{eq:term1_bound1}
\end{align}
Now, for large enough $n$, noting that $\pi\left(\mathcal G^c_n\right)\leq\exp(-\alpha n)$ up to some positive constant, we have 
\begin{align}
&\int_{\mathcal S^c}
P\left(\left|\frac{1}{2\sigma^2}-\frac{1}{2\sigma^2_0}\right|\left|\frac{2\bmu^\prime_n\bC_n\by_n}{n}\right|>\frac{\kappa_1}{4}\right)d\pi(\btheta)\notag\\
&\leq 2\int_{\mathcal S^c}\exp\left(-\frac{n\kappa^2_1}{18K_1}\left|\frac{1}{2\sigma^2}-\frac{1}{2\sigma^2_0}\right|^{-2}\right)d\pi(\btheta)
\label{eq:term1_interm1}\\
&\leq 2\int_{\mathcal G_n}\exp\left(-\frac{n\kappa^2_1}{18K_1}\left|\frac{1}{2\sigma^2}-\frac{1}{2\sigma^2_0}\right|^{-2}\right)d\pi(\btheta)
+ 2\int_{\mathcal G^c_n}\exp\left(-\frac{n\kappa^2_1}{18K_1}\left|\frac{1}{2\sigma^2}-\frac{1}{2\sigma^2_0}\right|^{-2}\right)d\pi(\btheta)\notag\\
	&\leq 2\int_{\exp(-2\left(\alpha n\right)^{1/16})}^{\exp(2\left(\alpha n\right)^{1/16})}\exp\left(-\frac{n\kappa^2_1}{18K_1}\left|\frac{1}{2\sigma^2}-\frac{1}{2\sigma^2_0}\right|^{-2}\right)\pi(\sigma^2)d\sigma^2
+2\pi\left(\mathcal G^c_n\right)\notag\\ 
	&\leq 2\int_{\exp(-2\left(\alpha n\right)^{1/16})-\sigma^{-2}_0}^{\exp(2\left(\alpha n\right)^{1/16})-\sigma^{-2}_0}\exp\left(-C_1\kappa^2_1Tu^{-2}\right)(u+\sigma^{-2}_0)^{-2}\pi\left(\frac{1}{u+\sigma^{-2}_0}\right)du
+\tilde C\exp(-\alpha n),
\label{eq:term1_bound3}
\end{align}
for some positive constants $C_1$ and $\tilde C$.

Now, the prior $(u+\sigma^{-2}_0)^{-2}\pi\left(\frac{1}{u+\sigma^{-2}_0}\right)$ is such that large values of $u$ receive
small probabilities. Hence, if this prior is replaced by an appropriate function which has a thicker tail than the prior, then the resultant integral provides an upper bound
for the first term of (\ref{eq:term1_bound3}). We consider a function 
$\tilde\pi(u)$ which is of mixture form depending upon $n$, that is, 
we let $\tilde\pi_n(u)=c_3\sum_{r=1}^{M_n}\psi^{\zeta_{rn}}_{rn}\exp(-\psi_{rn} u^2)u^{2(\zeta_{rn}-1)}\bI_{B_n}(u)$, 
where $B_n=\left[\exp\left(-2\left(\alpha n\right)^{1/16}\right)-\sigma^{-2}_0,\exp\left(2\left(\alpha n\right)^{1/16}\right)-\sigma^{-2}_0\right]$, 
$M_n\leq\exp(\left(\alpha n\right)^{1/16})$ is the number of mixture components, $c_3>0$,  
for $r=1,\ldots,M_n$, $\frac{1}{2}<\zeta_{rn}\leq c_4n^q$, for $0<q<1/16$ and $n\geq 1$, where $c_4>0$, and $0<\psi_1\leq\psi_{rn}<c_5<\infty$, for all $r$ and $n$. 
In this case,
\begin{align}
	& \int_{\exp(-2\left(\alpha n\right)^{1/16})-\sigma^{-2}_0}^{\exp(2\left(\alpha n\right)^{1/16})-\sigma^{-2}_0}\exp\left(-C_1\kappa^2_1nu^{-2}\right)(u+\sigma^{-2}_0)^{-2}\pi\left(\frac{1}{u+\sigma^{-2}_0}\right)du\notag\\
&\leq c_3\sum_{r=1}^{M_n}\psi^{\zeta_{rn}}_{rn}
	\int_{\exp(-2\left(\alpha n\right)^{1/16})-\sigma^{-2}_0}^{\exp(2\left(\alpha n\right)^{1/16})-\sigma^{-2}_0}\exp\left[-\left(C_1\kappa^2_1nu^{-2}+\psi_{rn}u^2\right)\right]\left(u^2\right)^{\zeta_{rn}-1}du.
\label{eq:term1_bound3_1}
\end{align}
Now the $r$-th integrand of (\ref{eq:term1_bound3_1}) is minimized at 
$\tilde u^2_{rn}= \frac{\zeta_{rn}-1+\sqrt{(\zeta_{rn}-1)^2+4C_1\psi_{rn}\kappa^2_1 n}}{2\psi_{rn}}$, so that for sufficiently large $n$, 
$c_1\kappa_1\sqrt{\frac{n}{\psi_{rn}}}\leq\tilde u^2_{rn}\leq \tilde c_1\kappa_1\sqrt{\frac{n}{\psi_{rn}}}$, for some positive constants $c_1$ and $\tilde c_1$. 
Now, for sufficiently large $n$, we have $\frac{\tilde u^2_{rn}}{\log\tilde u^2_{rn}}\geq\frac{\zeta_{rn}-1}{\psi_{rn}(1-c_2)}$, for $0<c_2<1$.
Hence, for sufficiently large $n$, $C_1\kappa^2_1n\tilde u^{-2}_{rn}+\psi_{rn}\tilde u^2_{rn}-(\zeta_{rn}-1)\log(\tilde u^2_{rn})\geq c_2\psi_1\tilde u^2_{rn}
\geq C_2\kappa_1\sqrt{\psi_{rn} n}$ for some 
positive constant $C_2$. From these and (\ref{eq:term1_bound3}) it follows that
\begin{align}
&c_3\sum_{r=1}^{M_n}\psi^{\zeta_{rn}}_{rn}
	\int_{\exp(-2\left(\alpha n\right)^{1/16})-\sigma^{-2}_0}^{\exp(2\left(\alpha n\right)^{1/16})-\sigma^{-2}_0}\exp\left[-\left(C_1\kappa^2_1nu^{-2}+\psi_1u^2\right)\right]\left(u^2\right)^{\zeta_{rn}-1}du\notag\\
	&\leq c_3M_n\exp\left[-\left(C_2\kappa_1\sqrt{n\psi_1}-2\left(\alpha n\right)^{1/16}-\tilde c_5 n^q\right)\right]\notag\\
	&\leq c_3\exp\left[-\left(C_2\kappa_1\sqrt{n}-3\left(\alpha n\right)^{1/16}-\tilde c_5 n^q\right)\right].
\label{eq:term1_bound3_2}
\end{align}
for some constant $\tilde c_5$. 
Combining (\ref{eq:term1_bound3}), (\ref{eq:term1_bound3_1}) and (\ref{eq:term1_bound3_2}) we obtain
\begin{align}
&\int_{\mathcal S^c}P\left(\left|\frac{1}{2\sigma^2}-\frac{1}{2\sigma^2_0}\right|\left|\frac{2\bmu^\prime_n\bC_n\by_n}{n}\right|>\frac{\kappa_1}{4}\right)d\pi(\btheta)\notag\\
	&\leq K_2\exp\left[-\left(C_2\kappa_1\sqrt{n}-3\left(\alpha n\right)^{1/16}-\tilde c_5n^q\right)\right]+\tilde C\exp(-\alpha n).
\label{eq:term1_bound3_final}
\end{align}

For the second term of (\ref{eq:ce1}), since $\bmu_n$ is non-random, we can also view this as a set of independent realizations from any suitable
independent zero mean process with variance $\frac{c(\bbeta_0)}{1-\rho^2_0}$ on a compact set (due to (\ref{eq:B2})). In that case, by Hoeffding's inequality \ctp{Hoeffding63} we obtain 
\begin{align}
&\int_{\mathcal S^c}P\left(\left|\frac{1}{2\sigma^2}-\frac{1}{2\sigma^2_0}\right|\left|\frac{\bmu^\prime_n\bmu_n}{n}
-\frac{c(\bbeta_0)}{1-\rho^2_0}\right|>\frac{\kappa_1}{4}\right)d\pi(\btheta)\notag\\
	&\leq 2\int_{\exp(-2\left(\alpha n\right)^{1/16})}^{\exp(2(\alpha n)^{1/16})}
\exp\left(-K_3\kappa^2_1n\left|\frac{1}{2\sigma^2}-\frac{1}{2\sigma^2_0}\right|^{-2}\right)\pi(\sigma^2)d\sigma^2+\tilde C_2\exp(-\alpha n)\notag\\
	&\leq K_3\exp\left[-\left(C_3\kappa_1\sqrt{n\psi_2}-3\left(\alpha n\right)^{1/16}-\tilde c_5n^q\right)\right]+ \tilde C\exp(-\alpha n).
\label{eq:term1_bound4}
\end{align}
for some positive constants $K_3$ and $C_3$. 
The last step follows in the same way as (\ref{eq:term1_bound3_final}).

We now deal with the first term of (\ref{eq:ce2}). First note that $\|\bSigma_n\|^2_F\leq K_4n$, for some $K_4>0$, where $\|\bSigma_n\|^2_F$ is the Frobenius norm of $\bSigma_n$.
Also, any eigenvalue $\lambda$ of any matrix $\bA=(a_{ij})$ satisfies $|\lambda-a_{ii}|\leq \sum_{j\neq i}|a_{ij}|$, by the Gerschgorin's circle theorem (see, for example,
\ctn{Lange10}). In our case,
the rows of $\bSigma_n$ are summable and the diagonal elements are bounded for any $n$. The maximum row sum is attained by the middle row when $n$ is odd and the two middle rows when
$n$ is even. In other words, the maximum eigenvalue of $\bSigma_n$ remains bounded for all $n\geq 1$. That is, $\underset{n\geq 1}{\sup}~\|\bSigma_n\|_{op}<K_5$, for some positive constant $K_5$. 
Now observe that for the integral of the form $\int_{\sigma^2\in\tilde{\mathcal G_n}}\exp\left(-C_5\kappa^2_1n\left|\sigma^{-2}-\sigma^{-2}_0\right|^{-1}\right)\pi(\sigma^2)d\sigma^2$,
where $\tilde{\mathcal G_n}\subseteq\mathcal G_n$, we can obtain, using the same technique pertaining to (\ref{eq:term1_bound3_final}), that
\begin{multline}
\int_{\sigma^2\in\tilde {\mathcal G_n}}\exp\left(-C_5\kappa^2_1n\left|\sigma^{-2}-\sigma^{-2}_0\right|^{-1}\right)\pi(\sigma^2)d\sigma^2\\
	\leq C_7\exp\left[-\left(C_6\kappa_1\sqrt{n}-3\left(\alpha n\right)^{1/16}-\tilde c_5n^q\right)\right],
\label{eq:hw0}
\end{multline}
for relevant positive constants $C_6$, $\psi_3$ and $\tilde c_5$.
%
Then by the Hanson-Wright inequality, 
(\ref{eq:hw0}) and the same method for obtaining (\ref{eq:term1_bound3_final}), we obtain the following bound for 
the first term of (\ref{eq:ce2}):
\begin{align}
&\int_{\mathcal S^c}P\left(\left|\frac{1}{2\sigma^2}-\frac{1}{2\sigma^2_0}\right|\left|\frac{\by^\prime_n\bSigma_n\by_n}{n}-\tr\left(\frac{\bSigma_n}{n}\right)\right|
>\frac{\kappa_1}{4}\right)d\pi(\btheta)\notag\\
&\leq E_{\pi}\left[\exp\left[-K_6\min\left\{\frac{\frac{\kappa^2_1}{9}\left|\frac{1}{2\sigma^2}-\frac{1}{2\sigma^2_0}\right|^{-2}}{\|\frac{\bSigma_n}{n}\|^2_F},
\frac{\frac{\kappa_1}{3}\left|\frac{1}{2\sigma^2}-\frac{1}{2\sigma^2_0}\right|^{-1}}{\|\frac{\bSigma_n}{n}\|_{op}}\right\}\right]\bI_{\mathcal G_n}(\btheta)\right]+\tilde C\exp(-\alpha n)\notag\\
	&\leq K_7\exp\left[-\left(C_8\kappa_1\sqrt{n}-3\left(\alpha n\right)^{1/16}-\tilde c_5n^q\right)\right]+\tilde C\exp(-\alpha n),
\label{eq:term1_bound5}
\end{align} 
for relevant positive constants $K_7$, $C_8$, $\psi_4$ and $\tilde c_5$.

Using the same technique involving Hoeffding's bound for the second term of (\ref{eq:ce1}), it is easy to see that the second term of (\ref{eq:ce2}) satisfies the following:
\begin{align}
P\left(\left|\frac{1}{2\sigma^2}-\frac{1}{2\sigma^2_0}\right|\left|tr\left(\frac{\bSigma_n}{n}\right)-\frac{\sigma^2_0}{1-\rho^2_0}\right|>\frac{\kappa_1}{4}\right)
	&\leq \tilde K_3\exp\left[-\left(\tilde C_3\kappa_1\sqrt{n}-3\left(\alpha n\right)^{1/16}-\tilde c_5n^q\right)\right],\notag\\
&\qquad+ \tilde C\exp(-\alpha n),
\label{eq:ce1_2nd}
\end{align}
for relevant positive constants $\tilde K_3$, $\tilde C_3$, $\tilde\psi_2$ and $\tilde c_5$.

Hence, combining (\ref{eq:ce1}), (\ref{eq:ce2}), (\ref{eq:term1_bound4}), (\ref{eq:term1_bound5}) and (\ref{eq:ce1_2nd}), we obtain
\begin{align}
&E_{\pi}\left[P\left(\left|\frac{1}{2\sigma^2}-\frac{1}{2\sigma^2_0}\right|
\left|\frac{\sum_{t=1}^nx^2_t}{n}-\frac{\sigma^2_0}{1-\rho^2_0}-\frac{c(\bbeta_0)}{1-\rho^2_0}\right|>\kappa_1\right)\bI_{\mathcal S^c}(\btheta)\right]\notag\\
	&\leq K_8\exp\left[-\left(C_9\kappa_1\sqrt{n}-3\left(\alpha n\right)^{1/16}-\tilde c_5n^q\right)\right]+\tilde C\exp(-\alpha n),
\label{eq:term1_bound_final}
\end{align}
for relevant positive constants.

Let us now obtain a bound for $E_{\pi}\left[P\left(\left|\frac{\rho^2}{2\sigma^2}-\frac{\rho^2_0}{2\sigma^2_0}\right|
\left|\frac{\sum_{t=1}^sx^2_{t-1}}{s}-\frac{\sigma^2_0}{1-\rho^2_0}-\frac{c(\bbeta_0)}{1-\rho^2_0}\right|>\kappa_1\right)\bI_{\mathcal S^c}(\btheta)\right]$.
By the same way as above, we obtain, by first taking the expectation with respect to $\sigma^2\in\mathcal G_n$, the following:
\begin{align}
&E_{\pi}\left[P\left(\left|\frac{\rho^2}{2\sigma^2}-\frac{\rho^2_0}{2\sigma^2_0}\right|
\left|\frac{\sum_{t=1}^nx^2_{t-1}}{n}-\frac{\sigma^2_0}{1-\rho^2_0}-\frac{c(\bbeta_0)}{1-\rho^2_0}\right|>\kappa_1\right)\bI_{\mathcal S^c}(\btheta)\right]\notag\\
	&\leq C_{10}\int_{\rho\in\mathcal G_n}\int_{\exp(-2\left(\alpha n\right)^{1/16})}^{\exp(2\left(\alpha n\right)^{1/16})}\exp\left[-C_{11}\kappa^2_1n\left(\frac{\rho^2}{\sigma^2}-\frac{\rho^2_0}{\sigma^2_0}\right)^{-2}\right]\pi(\sigma^2)d\sigma^2\pi(\rho)d\rho
+\tilde C\exp(-\alpha n)\notag\\
	&=C_{10}\int_{\rho\in\mathcal G_n}\rho^2\int_{\rho^2\exp(-2\left(\alpha n\right)^{1/16})-\frac{\rho^2_0}{\sigma^2_0}}^{\rho^2\exp(2\left(\alpha n\right)^{1/16})-\frac{\rho^2_0}{\sigma^2_0}}\exp\left(-C_{11}\kappa^2_1nu^{-2}\right)
\left(u+\frac{\rho^2_0}{\sigma^2_0}\right)^{-2}
\pi\left(\frac{\rho^2}{u+\frac{\rho^2_0}{\sigma^2_0}}\right)du\pi(\rho)d\rho\notag\\
&\qquad\qquad\qquad\qquad\qquad\qquad\qquad\qquad\qquad\qquad\qquad\qquad\qquad\qquad+\tilde C\exp(-\alpha n),
\label{eq:term2_bound1}
\end{align}
for relevant positive constants. 
Since $\pi\left(\sigma^2>\exp(2\left(\alpha n\right)^{1/16})\right)\leq \exp(-\alpha n)$, it is evident that much the mass of 
$\left(u+\frac{\rho^2_0}{\sigma^2_0}\right)^{-2}\pi\left(\frac{\rho^2}{u+\frac{\rho^2_0}{\sigma^2_0}}\right)$ is concentrated around zero, where the function
$\exp\left(-C_{11}nu^{-2}\right)$ is small. To give greater weight to the function, we can replace 
$\left(u+\frac{\rho^2_0}{\sigma^2_0}\right)^{-2}\pi\left(\frac{\rho^2}{u+\frac{\rho^2_0}{\sigma^2_0}}\right)$
with a mixture function of the form $\tilde\pi_{\rho^2,n}(u)
=c_3\sum_{r=1}^{M_n}\rho^{2\zeta_{rn}}\psi^{\zeta_{rn}}_{rn}\exp\left(-u^2\psi_{rn}\rho^2\right)\left(u^2\right)^{(\zeta_{rn}-1)}\bI_{B_{n,\rho^2}}(u)$, 
for positive constants $0<\psi_2\leq\psi_{rn}<c_5<\infty$ and $1/2<\zeta_{rn}<c_4n^q$. Here $$B_{n,\rho^2}=\left[\rho^2\exp(-2\left(\alpha n\right)^{1/16})
-\frac{\rho^2_0}{\sigma^2_0},
\rho^2\exp(2\left(\alpha n\right)^{1/16})-\frac{\rho^2_0}{\sigma^2_0}\right].$$
As before, $0<q<1/16$ and $M_n\leq\exp\left(\left(\alpha n\right)^{1/16}\right)$.
Hence, up to some positive constant,
\begin{align}
	&\int_{\rho^2\exp(-2\left(\alpha n\right)^{1/16})-\frac{\rho^2_0}{\sigma^2_0}}^{\rho^2\exp(2\left(\alpha n\right)^{1/16})-\frac{\rho^2_0}{\sigma^2_0}}\exp\left(-C_{11}\kappa^2_1nu^{-2}\right)
\left(u+\frac{\rho^2_0}{\sigma^2_0}\right)^{-2}\pi\left(\frac{\rho^2}{u+\frac{\rho^2_0}{\sigma^2_0}}\right)du\notag\\
	&\leq \sum_{r=1}^{M_n}\rho^{2\zeta_{rn}}\psi^{\zeta_{rn}}_{rn}
	\int_{\rho^2\exp(-2\left(\alpha n\right)^{1/16})-\frac{\rho^2_0}{\sigma^2_0}}^{\rho^2\exp(2\left(\alpha n\right)^{1/16})-\frac{\rho^2_0}{\sigma^2_0}}
\exp\left[-\left(C_{11}\kappa^2_1nu^{-2}+\psi_{rn}\rho^2u^2-(\zeta_{rn}-1)\log u^2\right)\right]du.
\label{eq:root1}
\end{align}
The term within the parenthesis in the exponent of (\ref{eq:root1}) is minimized at 
$\tilde u^2_{rn}=\frac{\zeta_{rn}-1+\sqrt{(\zeta_{rn}-1)^2+4\psi_{rn}\rho^2C_{11}\kappa^2_1n}}{2\psi_{rn}\rho^2}$.
Note that $\tilde C_{01}\frac{\kappa_1}{|\rho|}\sqrt{\frac{n}{\psi_{rn}}}\leq\tilde u^2_{rn}\leq\tilde C_{11}\frac{\kappa_1}{|\rho|}\sqrt{\frac{n}{\psi_{rn}}}$, 
for large enough $n$. Hence, for large $n$, the 
term within the parenthesis in the exponent of (\ref{eq:root1}) exceeds $\psi_{rn}\tilde u^2\geq \tilde C_{02}\times|\rho|\kappa_1\sqrt{\psi_{rn} n}$, for $\tilde C_{02}>0$.
Thus, (\ref{eq:root1}) is bounded above by a constant times $\rho^{2(1+\zeta_{rn})}\exp\left(-\tilde C_{02}\times\kappa_1|\rho|\sqrt{\psi_6 n}+3\left(\alpha n\right)^{1/16}+\tilde c_5 n^q\right)$.
Combining this with (\ref{eq:term2_bound1}) 
we see that
\begin{align}
&E_{\pi}\left[P\left(\left|\frac{\rho^2}{2\sigma^2}-\frac{\rho^2_0}{2\sigma^2_0}\right|
\left|\frac{\sum_{t=1}^nx^2_{t-1}}{n}-\frac{\sigma^2_0}{1-\rho^2_0}-\frac{c(\bbeta_0)}{1-\rho^2_0}\right|>\kappa_1\right)\bI_{\mathcal S^c}(\btheta)\right]\notag\\
	&\leq\int_{\rho\in\mathcal G_n}\rho^{2(2+\zeta_{rn})}\exp\left[-\left(\tilde C_{02}\times\kappa_1|\rho|\sqrt{\psi_6 n}-3\left(\alpha n\right)^{1/4}-\tilde c_5n^q\right)\right]\pi(\rho)d\rho+\tilde C\exp(-\alpha n)\notag\\
&=\int_{\exp\left(-\left(\alpha n\right)^{1/16}\right)}^{\exp\left(\left(\alpha n\right)^{1/16}\right)}	
	\exp\left[-\left(\tilde C_{02}\times\kappa_1u^{-1}\sqrt{\psi_6 n}+2(2+\zeta_{rn})\log u
	-3\left(\alpha n\right)^{1/16}-\tilde c_5n^q\right)\right]\pi_1(u)du\notag\\
	&\qquad+\tilde C\exp(-\alpha n),\label{eq:root2}
\end{align}
where $\pi_1(u)du$ is the appropriate modification of $\pi(\rho)d\rho$ in view of the transformation $|\rho|\mapsto u^{-1}$.
Replacing $\pi_1(u)$ with a mixture function of the form $\tilde\pi_{n}(u)
=c_3\sum_{r=1}^{M_n}\psi^{\tilde\zeta_{rn}}_{rn}\exp\left(-u\psi_{rn}\right)u^{(\tilde\zeta_{rn}-1)}$, 
for positive constants $0<\psi_2\leq\psi_{rn}<\tilde c_5<\infty$ and $0<\zeta_{rn}<c_4n^q$, with $0<q<1/16$, and $M_n\leq\exp\left(\left(\alpha n\right)^{1/16}\right)$,
and applying the same techniques as before, we see from (\ref{eq:root2}) that 
\begin{align}
&E_{\pi}\left[P\left(\left|\frac{\rho^2}{2\sigma^2}-\frac{\rho^2_0}{2\sigma^2_0}\right|
\left|\frac{\sum_{t=1}^nx^2_{t-1}}{n}-\frac{\sigma^2_0}{1-\rho^2_0}-\frac{c(\bbeta_0)}{1-\rho^2_0}\right|>\kappa_1\right)\bI_{\mathcal S^c}(\btheta)\right]\notag\\
	&\leq C_{14}\exp\left(3\left(\alpha n\right)^{1/4}+\tilde c_5n^q\right)\notag\\
	&\times\sum_{t=1}^{M_n}\psi^{\tilde\zeta_{rn}}_{rn}\int_{\exp\left(-\left(\alpha n\right)^{1/4}\right)}^{\exp\left(\left(\alpha n\right)^{1/4}\right)}
	\exp\left[-\left(\tilde C_{02}\times\kappa_1u^{-1}\sqrt{\psi_6 n}+u\psi_{rn}-(\tilde\zeta_{rn}-2\zeta_{rn}-5)\log u\right)\right]du\notag\\
	&\qquad+\tilde C\exp(-\alpha n)\notag\\
	&\leq C_{14}\exp\left[-\left(C_{15}\sqrt{\kappa_1}n^{1/4}-4\left(\alpha n\right)^{1/16}-2n^q\log \tilde c_5\right)\right]+\tilde C\exp(-\alpha n),
\label{eq:term2_bound_final}
\end{align}
for relevant positive constants. 

Let us now deal with $\frac{1}{2\sigma^2}\left(\bbeta^\prime_m\left(\frac{\sum_{t=1}^n\bz_{mt}\bz^\prime_{mt}}{n}\right)\bbeta_m-c(\bbeta)\right)
=\frac{1}{2\sigma^2}\left(\frac{\sum_{t=1}^n\left(\bz^\prime_{mt}\bbeta_m\right)^2}{n}-c(\bbeta)\right)$. 
Now, again we assume as before that $\bz^\prime_{mt}\bbeta_m$; $t=1,2,\ldots,n$ is a realization from some independent zero-mean process with variance $c(\bbeta)$. Note that
$|\bz^\prime_{mt}\bbeta_m|\leq\sum_{i=1}^m|z_{it}||\beta_i|=\sum_{i=1}^m|z_{it}||\gamma_i||\tilde\eta_i|\leq \underset{t\geq 1}{\sup}~\|z_t\|\|\tilde\eta\|\sum_{i=1}^m|\gamma_i|$.
By (B1), $\underset{t\geq 1}{\sup}~\|z_t\|<\infty$. 
Let $\tilde\gamma_m=\sum_{i=1}^m|\gamma_i|$. 
Then using Hoeffding's inequality in conjunction with (\ref{eq:B2}), we obtain
\begin{equation}
P\left(\frac{1}{2\sigma^2}\left|\frac{\sum_{t=1}^n(\bz^\prime_{mt}\bbeta_m)^2}{n}-c(\bbeta)\right|>\kappa_1\right)
<2\exp\left(-\frac{n\kappa^2_1\sigma^4}{C^2\tilde\gamma^4_m\|\tilde\eta\|^4}\right).
\label{eq:term3_bound1}
\end{equation}
Then, first integrating with respect to $u=\sigma^{-2}$, then integrating with respect to $v=\|\tilde\eta\|$
and finally with respect to $w=\tilde\gamma_m$, in each case using the gamma mixture form 
$\tilde\pi_n(x)=c_3\sum_{r=1}^{M_n}\psi^{\tilde\zeta_{rn}}_{rn}\exp\left(-x\psi_{rn}\right)x^{(\tilde\zeta_{rn}-1)}$, 
for positive constants $0<\psi_2\leq\psi_{rn}<\tilde c_5<\infty$ and $0<\zeta_{rn}<c_4n^q$, with $0<q<1/16$, and $M_n\leq\exp\left(\left(\alpha n\right)^{1/16}\right)$,
we find that
\begin{multline}
E_{\pi}\left[P\left(\frac{1}{2\sigma^2}\left|\frac{\sum_{t=1}^n(\bz^\prime_{mt}\bbeta_m)^2}{n}-c(\bbeta)\right|>\kappa_1\right)\bI_{\mathcal S^c}(\btheta)\right]\\
	\leq K_9\exp\left[-\left(C_{16}\kappa^{1/4}_1\left(n\psi_7\right)^{1/8}-C_{17}\left(\alpha n\right)^{1/16}-\tilde c_5n^q\right)\right]
+\tilde C\exp(-\alpha n),
\label{eq:term3_bound_final}
\end{multline}
for relevant positive constants. 
It is also easy to see using Hoeffding's inequality using (\ref{eq:B2}) that
\begin{align}
E_{\pi}\left[P\left(\frac{1}{2\sigma^2_0}\left|\frac{\sum_{t=1}^n(\bz^\prime_{mt}\bbeta_{m0})^2}{n}-c(\bbeta_0)\right|>\kappa_1\right)\bI_{\mathcal S^c}(\btheta)\right]
&\leq \tilde K_9\exp\left[-\left(\tilde C_{16}\kappa^2_1 n\right)\right],
\label{eq:term3_bound_final2}
\end{align}
for relevant constants.
%

We next consider $P\left(\left|\frac{\rho}{\sigma^2}-\frac{\rho_0}{\sigma^2_0}\right|\left|\frac{\rho_0\sum_{t=1}^nx^2_{t-1}}{n}+\frac{\bbeta^\prime_{m0}\sum_{t=1}^n\bz_{mt}x_{t-1}}{n}
-\frac{\rho_0\sigma^2_0}{1-\rho^2_0}-\frac{\rho_0c(\bbeta_0)}{1-\rho^2_0}\right|>\kappa_1\right)$. Note that
\begin{align}
&P\left(\left|\frac{\rho}{\sigma^2}-\frac{\rho_0}{\sigma^2_0}\right|\left|\frac{\rho_0\sum_{t=1}^nx^2_{t-1}}{n}+\frac{\bbeta^\prime_{m0}\sum_{t=1}^n\bz_{mt}x_{t-1}}{n}
-\frac{\rho_0\sigma^2_0}{1-\rho^2_0}-\frac{\rho_0c(\bbeta_0)}{1-\rho^2_0}\right|>\kappa_1\right)\notag\\
	&\qquad\qquad\leq P\left(\left|\frac{\rho}{\sigma^2}-\frac{\rho_0}{\sigma^2_0}\right|\left|\frac{\sum_{t=1}^nx^2_{t-1}}{n}
	-\frac{\sigma^2_0}{1-\rho^2_0}-\frac{c(\bbeta_0)}{1-\rho^2_0}\right|>\frac{\kappa_1}{2\rho_0}\right)\label{eq:term4_bound1}\\
	&\qquad\qquad\qquad\qquad+P\left(\left|\frac{\rho}{\sigma^2}-\frac{\rho_0}{\sigma^2_0}\right|\left|\frac{\bbeta^\prime_{m0}\sum_{t=1}^n\bz_{mt}x_{t-1}}{n}\right|>\frac{\kappa_1}{2}\right).
\label{eq:term4_bound2}
\end{align}
Note that the expectation of \eqref{eq:term4_bound1} admits the same upper bound as (\ref{eq:term2_bound_final}).
To deal with (\ref{eq:term4_bound2}) we let $\tilde x_t=(\bz^\prime_t\bbeta_0)x_{t-1}$ and $\tilde\bx_n=(\tilde x_1,\ldots,\tilde x_n)^\prime$. Then
$\tilde\bx_n\sim N_n\left(\tilde\bmu_n,\tilde\bSigma_n\right)$, where $\tilde\bmu_n$ and $\tilde\bSigma_n=\tilde\bC_n\tilde\bC^\prime_n$ are appropriate modifications of 
$\bmu_n$ and $\bSigma_n=\bC_n\bC^\prime_n$ associated with (\ref{eq:term1_bound1}). Note that $\tilde\bx_n=\tilde\bmu_n+\tilde\bC_n\by_n$, where
$\by_n\sim N_n\left(\bzero_n,\bI_n\right)$. Using (\ref{eq:B2}) we obtain the same form of the bound for (\ref{eq:term4_bound2})
as (\ref{eq:term1_bound1}). That is, we have
\begin{align}
&P\left(\left|\frac{\rho}{\sigma^2}-\frac{\rho_0}{\sigma^2_0}\right|\left|\frac{\bbeta^\prime_{m0}\sum_{t=1}^n\bz_{mt}x_{t-1}}{n}\right|>\frac{\kappa_1}{2}\right)\notag\\
&\leq P\left(\left|\bone^\prime_n\tilde\bC_n\by_n\right|>\frac{n\kappa_1}{4}\left|\frac{\rho}{\sigma^2}-\frac{\rho_0}{\sigma^2_0}\right|^{-1}\right)
+P\left(\left|\bmu^\prime_n\bone_n\right|>\frac{n\kappa_1}{4}\left|\frac{\rho}{\sigma^2}-\frac{\rho_0}{\sigma^2_0}\right|^{-1}\right)\notag\\
&\leq 2\exp\left(-K_{10}\kappa^2_1n\left|\frac{\rho}{\sigma^2}-\frac{\rho_0}{\sigma^2_0}\right|^{-2}\right)
+P\left(\left|\bmu^\prime_n\bone_n\right|>\frac{n\kappa_1}{4}\left|\frac{\rho}{\sigma^2}-\frac{\rho_0}{\sigma^2_0}\right|^{-1}\right),
\label{eq:term4_bound3}
\end{align}
where $K_{10}$ is some positive constant. Using the same method as before 
again we obtain a bound for the expectation
of the first part of (\ref{eq:term4_bound3}) of similar form as 
$\exp\left[-\left(\tilde C_{16}\sqrt{\kappa_1}n^{1/4}-\tilde C_{17}\left(\alpha n\right)^{1/16}-\tilde\alpha_5n^q\right)\right]+\tilde C\exp(-\alpha n)$,
for relevant positive constants. As before, here $0<q<1/16$.
For the second part of (\ref{eq:term4_bound3}) we apply the method involving Hoeffding's inequality as before, and obtain a bound of the above-mentioned form.  
Hence combining the bounds for the expectations of \eqref{eq:term3_bound1} and (\ref{eq:term4_bound2}) we see that
\begin{align}
&E_{\pi}\left[P\left(\left|\frac{\rho}{\sigma^2}-\frac{\rho_0}{\sigma^2_0}\right|\left|\frac{\rho_0\sum_{t=1}^nx^2_{t-1}}{n}+\frac{\bbeta^\prime_{m0}\sum_{t=1}^n\bz_{mt}x_{t-1}}{n}
-\frac{\rho_0\sigma^2_0}{1-\rho^2_0}-\frac{\rho_0c(\bbeta_0)}{1-\rho^2_0}\right|>\kappa_1\right)\bI_{\mathcal S^c}(\btheta)\right]\notag\\
	&\leq K_{12}\exp\left[-\left(C_{18}\sqrt{\kappa_1}n^{1/4}-C_{19}\left(\alpha n\right)^{1/16}-\tilde\alpha_5n^q\right)\right]
	+\tilde C\exp(-\alpha n),
\label{eq:term4_bound_final}
\end{align}
for relevant positive constants.

Now let us bound the probability 
$P\left(\left|\left(\frac{\bbeta_m}{\sigma^2}-\frac{\bbeta_{m0}}{\sigma^2_0}\right)^\prime\left(\frac{\sum_{t=1}^n\bz_{mt}x_t}{n}\right)-\frac{c_{10}(\bbeta,\bbeta_0)}{\sigma^2}
+\frac{c(\bbeta_0)}{\sigma^2_0}\right|
>\kappa_1\right)$.
Observe that
\begin{align}
&P\left(\left|\left(\frac{\bbeta_m}{\sigma^2}-\frac{\bbeta_{m0}}{\sigma^2_0}\right)^\prime\left(\frac{\sum_{t=1}^n\bz_{mt}x_t}{n}\right)
-\frac{c_{10}(\bbeta,\bbeta_0)}{\sigma^2}+\frac{c(\bbeta_0)}{\sigma^2_0}\right|
>\kappa_1\right)
\notag\\
&\leq P\left(\left|\frac{\sum_{t=1}^n(\bz^\prime_{mt}\bbeta_m)x_t}{n}-c_{10}(\bbeta,\bbeta_0)\right|>\frac{\kappa_1\sigma^2}{2}\right)
+ P\left(\left|\frac{\sum_{t=1}^n(\bz^\prime_{mt}\bbeta_{m0})x_t}{n}-c(\bbeta_0)\right|>\frac{\kappa_1\sigma^2_0}{2}\right).
\label{eq:term5_bound1}
\end{align}
Using the Gaussian concentration inequality as before it is easily seen that 
\begin{align}
&E_{\pi}\left[P\left(\left|\frac{\sum_{t=1}^n(\bz^\prime_{mt}\bbeta_m)x_t}{n}-c_{10}(\bbeta,\bbeta_0)\right|>\frac{\kappa_1\sigma^2}{2}\right)\bI_{\mathcal S^c}(\btheta)\right]\notag\\
	&\leq 2\int_{\bgamma_m,\tilde\eta\in\mathcal G_n}\int_{\exp(-2\left(\alpha n\right)^{1/16})}^{\exp(2\left(\alpha n\right)^{1/16})}
\exp\left(-\frac{K_{13}\kappa^2_1n\sigma^4}{\|\bbeta\|^2}\right)d\pi(\bbeta,\sigma^2)+\tilde C\exp(-\alpha n)\notag\\
	&\leq C_{20}\exp\left[-\left(C_{21}\sqrt{\kappa_1}n^{1/4}-C_{22}\left(\alpha n\right)^{1/16}-\tilde c_5n^q\right)\right]+\tilde C\exp(-\alpha n),
\label{eq:term5_bound2}
\end{align}
for relevant positive constants.

The Gaussian concentration inequality also ensures that the second term of (\ref{eq:term5_bound1}) is bounded above by $2\exp(-K_{13}\kappa^2_1n)$, for some $K_{13}>0$.
Combining this with (\ref{eq:term5_bound1}) and (\ref{eq:term5_bound2}) we obtain
\begin{align}
&E_{\pi}\left[P\left(\left|\left(\frac{\bbeta_m}{\sigma^2}-\frac{\bbeta_{m0}}{\sigma^2_0}\right)^\prime
\left(\frac{\sum_{t=1}^n\bz_{mt}x_t}{n}\right)-\frac{c_{10}(\bbeta,\bbeta_0)}{\sigma^2}+\frac{c(\bbeta_0)}{\sigma^2_0}\right|
>\kappa_1\right)\bI_{\mathcal S^c}(\btheta)\right]\notag\\
	&\leq K_{14}\exp\left[-\left(C_{23}\sqrt{\kappa_1}n^{1/4}-C_{24}\left(\alpha n\right)^{1/16}-\tilde c_5n^q\right)\right]
+\tilde C\exp(-\alpha n)+2\exp(-K_{13}\kappa^2_1n),
\label{eq:term5_bound_final}
\end{align}
for relevant positive constants. Note that, here $0<q<1/16$. 

For $P\left(\left|\left(\frac{\rho\bbeta_m}{\sigma^2}-\frac{\rho_0\bbeta_{m0}}{\sigma^2_0}\right)^\prime\left(\frac{\sum_{t=1}^n\bz_{mt}x_{t-1}}{n}\right)\right|>\kappa_1\right)$, we note that
\begin{align}
&P\left(\left|\left(\frac{\rho\bbeta_m}{\sigma^2}-\frac{\rho_0\bbeta_{m0}}{\sigma^2_0}\right)^\prime\left(\frac{\sum_{t=1}^n\bz_{mt}x_{t-1}}{n}\right)\right|>\kappa_1\right)\notag\\
&\leq P\left(\left|\frac{\sum_{t=1}^n(\bz^\prime_{mt}\bbeta_m)x_{t-1}}{n}\right|>\frac{\kappa_1\sigma^2}{2\rho}\right)
+ P\left(\left|\frac{\sum_{t=1}^n(\bz^\prime_{mt}\bbeta_{m0})x_{t-1}}{n}\right|>\frac{\kappa_1\sigma^2_0}{2\rho_0}\right).
\label{eq:term6_bound1}
\end{align}
For the first term of (\ref{eq:term6_bound1}) we apply the Gaussian concentration inequality followed by taking expectations with respect to $\sigma^2$, $|\rho|$,  
$|\tilde \gamma_m|$ and $\|\tilde\eta\|$.
This yields the bound $$K_{15}\exp\left[-\left(C_{25}\kappa^{1/8}_1n^{1/16}-C_{26}\left(\alpha n\right)^{1/16}-n^q\log \tilde c_5\right)\right]+\tilde C\exp(-\alpha n),$$ 
for relevant positive constants.
The bound for the second term is given by $2\exp(-K_{16}\kappa^2_1n)$. Together we thus obtain
\begin{align}
&E_{\pi}\left[P\left(\left|\left(\frac{\rho\bbeta_m}{\sigma^2}-\frac{\rho_0\bbeta_{m0}}{\sigma^2_0}\right)^\prime\left(\frac{\sum_{t=1}^n\bz_{mt}x_{t-1}}{n}\right)\right|>\delta_1\right)
\bI_{\mathcal G_n}(\btheta)\right]\notag\\
	&\leq \tilde K_{16}\exp\left[-\left(C_{26}\kappa^{1/8}_1n^{1/16}-C_{27}\left(\alpha n\right)^{1/16}-n^q\log \tilde c_5\right)\right]+2\exp(-K_{16}\kappa^2_1n).
\label{eq:term6_bound_final}
\end{align}

We now deal with the last term $P\left(\left|\left(\frac{\rho}{\sigma^2}-\frac{\rho_0}{\sigma^2_0}\right)\left(\frac{\sum_{t=1}^n\epsilon_tx_{t-1}}{n}\right)\right|>\kappa_1\right)$.
Recall that $\bx_n=\bmu_n+\bC_n\by_n$, where $\bC_n\bC^\prime_n=\bSigma_n$ and $\by_n\sim N_n\left(\bmu_n,\bI_n\right)$. Let $\bepsilon_{n-1}=(\epsilon_2,\ldots,\epsilon_n)^\prime$.
Then $\sum_{t=1}^n\epsilon_tx_{t-1}=\bepsilon^\prime_{n-1}\bx_{n-1}=\sigma_0\left(\by^\prime_n\bmu_n+\by^\prime_{n-1}\bC_{n-1}\by_{n-1}\right)$. Application of the Gaussian concentration 
inequality and the Hanson-Wright inequality we find that
\begin{align}
&P\left(\left|\left(\frac{\rho}{\sigma^2}-\frac{\rho_0}{\sigma^2_0}\right)\left(\frac{\sum_{t=1}^n\epsilon_tx_{t-1}}{n}\right)\right|>\kappa_1\right)\notag\\
&\leq P\left(\frac{\left|\by^\prime_n\bmu_n\right|}{n}>\frac{\kappa_1}{\sigma_0}\left|\frac{\rho}{\sigma^2}-\frac{\rho_0}{\sigma^2_0}\right|^{-1}\right)
+ P\left(\frac{\by^\prime_{n-1}\bC_{n-1}\by_{n-1}}{n}>\frac{\kappa_1}{\sigma_0}\left|\frac{\rho}{\sigma^2}-\frac{\rho_0}{\sigma^2_0}\right|^{-1}\right)\notag\\
&\leq K_{17}\exp\left(-K_{18}\kappa^2_1n\left|\frac{\rho}{\sigma^2}-\frac{\rho_0}{\sigma^2_0}\right|^{-2}\right),
\label{eq:term7_bound1}
\end{align}
for some positive constants $K_{17}$ and $K_{18}$.
Taking expectation of (\ref{eq:term7_bound1}) with respect to $\pi$ we obtain as before 
\begin{multline}
E_{\pi}\left[P\left(\left|\left(\frac{\rho}{\sigma^2}-\frac{\rho_0}{\sigma^2_0}\right)\left(\frac{\sum_{t=1}^n\epsilon_tx_{t-1}}{n}\right)\right|>\kappa_1\right)
\bI_{\mathcal S^c}(\btheta)\right] \\
	\leq K_{19}\exp\left[-\left(K_{20}\sqrt{\kappa_1}n^{1/4}-K_{21}\left(\alpha n\right)^{1/16}-\tilde c_5n^q\right)\right]+\tilde C\exp(-\alpha n),
\label{eq:term7_bound_final}
\end{multline}
for relevant positive constants. Recall that $0<q<1/16$.

Combining (\ref{eq:term1_bound_final}), (\ref{eq:term2_bound_final}), (\ref{eq:term3_bound_final}), (\ref{eq:term4_bound_final}), (\ref{eq:term5_bound_final}), (\ref{eq:term6_bound_final})
and (\ref{eq:term7_bound_final}), 
we see that 
\begin{equation*}
\sum_{n=1}^{\infty}E_{\pi}\left[P\left(\left|\frac{1}{n}\log R_n(\btheta)+h(\btheta)\right|>\delta\right)\bI_{\mathcal S^c}(\btheta)\right]<\infty.
\end{equation*}
This verifies (\ref{eq:s6_7}) and hence \ref{s6}.

\subsection{Verification of \ref{shalizi7}}
\label{subsec:A7}
Since $\mathcal G_n\rightarrow\bTheta$ as $n\rightarrow\infty$, it follows that for any set $A$ with $\pi(A)>0$,
$\mathcal G_n\cap A\rightarrow\bTheta\cap A=A$, as $n\rightarrow\infty$.
In our case, $\mathcal G_n$, and hence $\mathcal G_n\cap A$, are decreasing in $n$, so that $h\left(\mathcal G_n\cap A\right)$ must be
non-increasing in $n$. Moreover, for any $n\geq 1$, $\mathcal G_n\cap A\subseteq A$, so that $h\left(\mathcal G_n\cap A\right)\geq h(A)$, for all $n\geq 1$.
Hence, continuity of $h$ implies that $h\left(\mathcal G_n\cap A\right)\rightarrow h(A)$, as $n\rightarrow\infty$, and (S7) is satisfied.

Thus (S1)--(S7) are satisfied, so that Shalizi's result stated in the main manuscript holds. It follows that all our asymptotic results of our main manuscript
apply to this multiple testing problem.

\bibliographystyle{natbib}
\bibliography{irmcmc}

\newpage

\end{document}